\documentclass[a4paper,twoside]{article} 
\usepackage{amsfonts,amssymb,mathrsfs,amsmath}
\usepackage[dvips]{color}
\usepackage{graphicx}
\title{Rigid Polyboxes and Keller's Conjecture}
\date{}
\author{Andrzej P. Kisielewicz\\
\\
{\small Wydzia{\l} Matematyki, Informatyki i Ekonometrii, Uniwersytet Zielonog\'orski}\\
{\small ul. Z. Szafrana 4a, 65-516 Zielona G\'ora, Poland}\\
{\small A.Kisielewicz@wmie.uz.zgora.pl}\\
}
\numberwithin{equation}{section}
\newtheorem{pr}{\sc Proposition}
\newtheorem{ex}[pr]{\sc Example}
\newtheorem{lemat}[pr]{\sc Lemma}
\newtheorem{st}[pr]{\sc Statement}
\newtheorem{tw}[pr]{\sc Theorem}
\newtheorem{wn}[pr]{\sc Corollary}
\newtheorem{df}{\sc Definition}

\newtheorem{uw}{\sc Remark}
\newtheorem{uwi}[uw]{\sc Remarks}

\newtheorem{nap}{\sc Example }

\newtheorem{nps}[nap]{\sc Examples}

\def\ka #1{\mathscr{#1}}
\def\kal #1 #2{\mathscr{#1}^{#2}}
\def\proof{\noindent \textit{Proof.\,\,\,}}

\def\zet{\mathbb{Z}}

\def\er{\mathbb{R}}

\def\Pud #1{\operatorname{Box}(#1)}

\def\iver #1{\mbox{\tt [} #1 \mbox{\tt]}}

\begin{document}

\numberwithin{pr}{section}
\numberwithin{uw}{section}
\maketitle
\begin{abstract}
A cube tiling of $\er^d$ is a  family of pairwise disjoint cubes $[0,1)^d+T=\{[0,1)^d+t:t\in T\}$ such that $\bigcup_{t\in T}([0,1)^d+t)=\er^d$. Two cubes $[0,1)^d+t$, $[0,1)^d+s$ are called a twin pair if their closures have a complete facet in common, that is if $|t_j-s_j|=1$ for some $j\in [d]=\{1,\ldots, d\}$ and $t_i=s_i$ for every $i\in [d]\setminus \{j\}$. 
In $1930$, Keller conjectured that in every cube tiling of $\er^d$ there is a twin pair. Keller's conjecture is true for dimensions $d\leq 6$ and false for all dimensions $d\geq 8$. For $d=7$ the conjecture is still open. Let $x\in \er^d$, $i\in [d]$, and let $L(T,x,i)$ be the set of all $i$th coordinates $t_i$ of vectors $t\in T$ such that $([0,1)^d+t)\cap ([0,1]^d+x)\neq \emptyset$ and $t_i\leq x_i$. Let $r^-(T)=\min_{x\in \er^d}\; \max_{1\leq i\leq d}|L(T,x,i)|$ and $r^+(T)=\max_{x\in \er^d}\; \max_{1\leq i\leq d}|L(T,x,i)|$. It is known that  Keller's conjecture is true in dimension seven for cube tilings $[0,1)^7+T$ for which $r^-(T)\leq 2$. In the present paper we show that it is also true for $d=7$ if $r^+(T)\geq 6$. Thus, if $[0,1)^d+T$ is a counterexample to Keller's conjecture in dimension seven, then $r^-(T),r^+(T)\in \{3,4,5\}$. 


\end{abstract}
\section{Introduction}
A {\it cube tiling} of $\er^d$ is a  family of pairwise disjoint cubes $[0,1)^d+T=\{[0,1)^d+t:t\in T\}$ such that $\bigcup_{t\in T}([0,1)^d+t)=\er^d$. Two cubes $[0,1)^d+t$, $[0,1)^d+s$ are called a {\it twin pair} if their closures have a complete facet in common, that is if $|t_j-s_j|=1$
for some $j\in [d]=\{1,\ldots, d\}$ and $t_i=s_i$ for every $i\in [d]\setminus \{j\}$. In $1907$, Minkowski \cite{Min} conjectured that in every {\it lattice} cube tiling of $\er^d$, i.e. when $T$ is a lattice in $\er^d$, there is a twin pair, and in $1930$, Keller \cite{Ke1} generalized this conjecture to any cube tiling of $\er^d$. Minkowski's conjecture was confirmed by Haj\'os \cite{H} in $1941$. In $1940$, Perron \cite{P} proved that Keller's conjecture is true for all dimensions $d\leq 6$. In 1986, Szab\'o \cite{Sz2} showed that if there is a counterexample to Keller's conjecture in dimension $d$, then there is a counterexample two-periodic cube tiling  $[0,1)^n+T$ of $\er^n$, where $T\subset (1/2)\zet^n$ and $d\leq n$. Moreover, Corr\'adi and Szab\'o \cite{CS2} reduced Keller's conjecture for $T\subset (1/2)\zet^d$  to a problem in graph theory. They defined {\it a $d$-dimensional Keller graph} whose vertices are all strings from the set  $\{0,1,2,3\}^d$. Two vertices are adjacent if they differ in at least two positions, but  in at one position the difference is two modulo four. Thus, Keller's cube tiling conjecture says that a maximum clique in a $d$-dimensional Keller graph has less than $2^d$ vertices.  
The results of  Corr\'adi and Szab\'o  inspired Lagarias and Shor \cite{LS1} who, in 1992, constructed a cube tiling of $\er^{10}$ which does not contain a twin pair  and thereby refuted  Keller's cube tiling conjecture. Finally, in $2002$, Mackey \cite{M} gave a counterexample to Keller's conjecture in dimension eight, which also shows that this conjecture is false in dimension nine. For $d=7$ Keller's conjecture is still open.  

Let  $[0,1)^d+T$ be a cube tiling, $x\in \er^d$ and $i\in [d]$, and let $L(T,x,i)$ be the set of all $i$th coordinates $t_i$ of vectors $t\in T$ such that $([0,1)^d+t)\cap ([0,1]^d+x)\neq \emptyset$ and $t_i\leq x_i$ (Figure 1). It is known that $1\leq |L(T,x,i)|\leq 2^{d-1}$ for every $x\in \er^d$ and every $i\in [d]$ (compare Sections 2.2 and 2.3).
% for every $x\in \er^d$ and $i\in [d]$.  

\smallskip
{\center
\includegraphics[width=5cm]{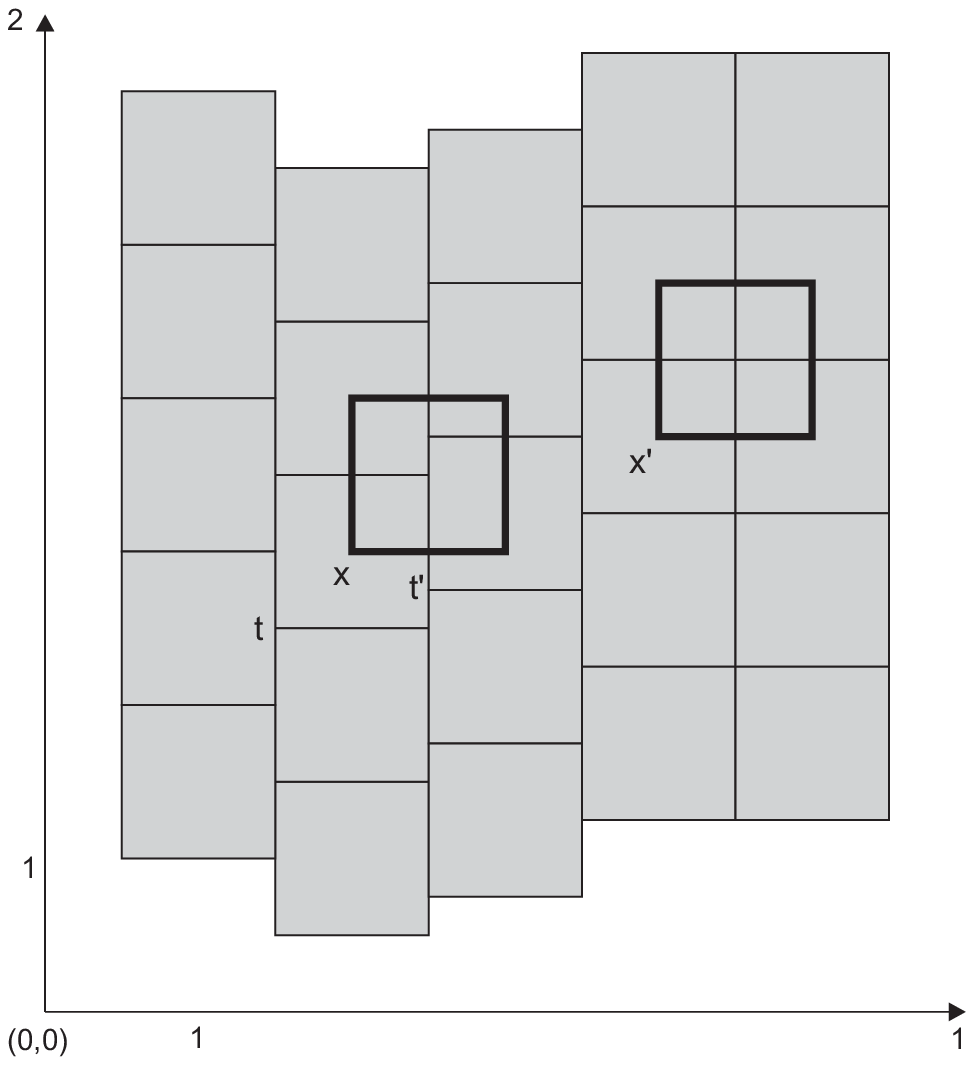}\\
}

\noindent{\footnotesize 
Fig. 1. A portion of a cube tiling $[0,1)^2+T$ of $\er^2$. The number of elements in $L(T,x,i)$ depends on the position of $x\in \er^2$. For $x=(2,3)$, we have $L(T,x,1)=\{3/2\}(=\{t_1\})$ and $L(T,x,2)=\{5/2,11/4\}(=\{t_2,t_2'\})$, while for $x'=(4,15/4)$, we have $L(T,x',1)=\{7/2\}$ and $L(T,x',2)=\{13/4\}$. This portion of the tiling  $[0,1)^2+T$ shows that $r^-(T)=1$ and $r^+(T)=2$.
  }

\medskip
%To every cube tiling  $[0,1)^n+T$ of $\er^n$ we assign the number $r(T)$ which is the maximum value of $|L(T,x,i)|$ over all $x\in \er^n$ and all $i\in [d]$. 
%$\begin{equation}
%$\label{re}
%r(T)=
Let 
\begin{equation}
\label{re}
r^-(T)=\min_{x\in \er^d}\; \max_{1\leq i\leq d}|L(T,x,i)|\;\;\;\; {\rm and}\;\;\;\; r^+(T)=\max_{x\in \er^d}\; \max_{1\leq i\leq d}|L(T,x,i)|.
\end{equation}
%and
%\begin{equation}
%\label{rem}
%M(T)=\max_{1\leq i\leq d}\max_{x\in \er^d}|L(T,x,i)|
%\end{equation}
%We have $r^-(T),r^+(T)\in [2^{d-1}]$.

In 2010, Debroni et al. \cite{De} computed that the maximum clique in the $7$-dimensional Keller graph has $124$ vertices, which implies that Keller's conjecture is true for all cube tilings $[0,1)^7+T$ of $\er^7$ with  $T\subset (1/2)\zet^7$ or equivalently, $T\subset a+\zet^7\cup b+\zet^7$, where fixed $a,b\in [0,1)^7$ are such that $a_i\neq b_i$ for every $i\in [7]$. Observe now that the condition $r^-(T)\leq 2$ means that  there is $x\in \er^7$ such that the set $L(T,x,i)$ contains at most two elements for every $i\in [7]$ or, equivalently, there is $x\in \er^7$ such that the set $T_1\subset T$ consisting of all $t$ for which $([0,1)^7+t)\cap ([0,1]^7+x)\neq \emptyset$ is a subset of the set $a+\zet^7\cup b+\zet^7$. Thus, it is easy to show that the result of Debroni et al. proves also that the conjecture is true for cube tilings of $\er^7$ for which $r^-(T)\leq 2$. Indeed, if there is no twin pair in the set $\{[0,1)^7+t:t\in T_1\}$, then extending this family to the two-periodic tiling $[0,1)^7+T$ of $\er^7$, where $T=T_1+2\zet^7$, we obtain a cube tiling with $T\subset a+\zet^7\cup b+\zet^7$ without twin pairs, which contradicts the result of Debroni et al.  

This paper is motivated by Keller's conjecture in dimension seven, the last unsolved case of this conjecture. We prove that 

\begin{tw}
\label{keli}
Keller's conjecture is true for cube tilings $[0,1)^7+T$ of $\er^7$ for which $r^+(T)\geq 6$.
\end{tw}
 
It follows from this theorem and the result of Debroni et al. that if $[0,1)^7+T$ is a counterexample to Keller's conjecture in dimension seven, then $r^-(T),r^+(T)\in \{3,4,5\}$ (Corollary \ref{ke2}). 

Similarly to Perron's approach in \cite{P} (see also \cite{LP1}), we based our methods on a deeper examination of the local structure of cube tilings. The main difference between Perron's approach and the one presented here is that we use the notion of a rigid system of boxes, which was introduced in \cite{LS2}, and widely examined in \cite{KP}. Roughly speaking, Perron's methods are combinatorial, while ours are strongly geometric.   

Works on Minkowski's and Keller's conjectures revealed a number of interesting problems concerning the structure of cube tilings. Lagarias and Shor \cite{LS2} formulated a new problem on the structure of cube tilings of $\er^d$: Let $K_d$ be the largest integer such that every cube tiling of $\er^d$ contains two cubes that have a common face of dimension $K_d$. What is the upper bound on $K_d$? Since Keller's conjecture is true for $d\leq 6$, we have $K_d=d-1$ for $d\leq 6$. Generally, in \cite{LS2} it was shown that $K_d\leq d-(1/3)\sqrt{d}$ for every $d$. Moreover, in that paper the authors considered subsets of $\er^d$ which can be represent as a union of disjoint unit cubes, satisfying certain additional condition, only in the one manner. This is, mentioned above, a rigid system of boxes (cubes in this case) which was considered in \cite{LS2} in order to improve the upper bound on $K_d$. In the present paper the rigidity of a system of boxes is one of the crucial tools in examining the structure of systems of boxes. 
%For example, in \cite{Kis1} we showed how a cube tiling code designed in \cite{LS2} can be used to obtain an interesting partitions and matchings of a $d$-dimensional cube. 

The maximum clique in the Keller graph in dimension six contains $60$ vertices which was computed by David Applegate (see Section 1 in \cite{De}). This proves that Keller's conjecture is true for all cube tilings $[0,1)^6+T$ of $\er^6$ such that $r^-(T)\leq 2$. On the other hand, the results obtained in the presented paper
%as we will see in the last section of the paper,  
(mainly Theorem \ref{12}) %which is a base for proving Theorem \ref{keli} 
prove that  Keller's conjecture is true for all cube tilings $[0,1)^6+T$ of $\er^6$ such that $r^+(T)\geq 3$.  Thus, these two results give a new proof of Keller's conjecture in dimensions $d\leq 6$. We will present this proof in the final section of the paper. 

A new approach to Minkowski's conjecture can be found in Kolountzakis's paper \cite{Ko1} and in \cite{KP2}. The Hajos proof of Minkowski's conjecture stimulates the development work on the factorization of abelian groups. These topics are examined in Szabo's book \cite{Sz1}. A fine survey of  tilings of $\er^d$ by clusters of unit cubes and Minkowski's conjecture is Stein's and Szabo's book \cite{SS}. 

It is not widely known that Keller made two conjectures on twin pairs in a cube tiling of $\er^d$. The second conjecture, posed in \cite{Ke2}, says that every cube tiling of $\er^d$ contains {\it a column} of unit cubes, i.e. a family of the form $\{[0,1)^d+t+ne_i: n\in \zet\}$, where $e_i$ is the $i$-th element of the standard basis of $\er^d$. This conjecture has been proved for $d\leq 6$ by \L  ysakowska and Przes\l awski in \cite{LP1}. Furthermore, these authors, in \cite{LP2}, described the meta-structure of cube tilings of $\er^3$ and non-extensible  systems of unit cubes. Such systems were also examined, among other things, by  Dutour Sikiri\'c and Itoh in \cite{DI}.

The present paper is organized as follows. In Section 2 we give basic notions concerning the systems of boxes and abstract words. These issues were developed in \cite{GKP,KP}. Since they are not widely known, we present them in detail. In Section 3 we describe the structure of  systems of abstract words. Next, in Section 4, we prove the fundamental for our purposes theorem on systems of words (Theorem \ref{12}). 
This theorem implies immediately Theorem \ref{keli}, which will be shown in the final Section 5. 
Actually, we prove in that section a general theorem on the existence of twin pairs in cube tilings of $\er^d$.
At the end of the paper we give an interpretation of our results for cliques in a $d$-dimensional Keller graph.

%Moreover, we extend the definition of a $d$-dimensional Keller graph and give an interpretation of our results in the terms of that graph. 
%Related topics dealing with the structure of polybox codes and cube tilings can be found in \cite{DI,SS}.
%and stimulated the development work on the other partitions than those on the unit cubes of the spaoxesce $\er^d$.  
%In \cite{GKP,KP} systems of boxundamentales were studied. 

\section{Basic notions}

In this section we present the basic notions on dichotomous boxes and words (details can be found in \cite{GKP,KP}). We start with systems of boxes.

In the whole paper, if $\ka X$ is a family of sets, then $\bigcup \ka X=\bigcup_{A\in \ka X}A$. Moreover, if $Y$ is a set, then a {\it partition of} $Y$ is a family $\ka Y$ of its pairwise disjoint subsets such that $\bigcup \ka Y=Y$.  

\subsection{Dichotomous boxes and polyboxes}

Let $X_1,\ldots ,X_d$ be non-empty sets with $|X_i|\geq 2$ for every $i\in [d]$. The set $X=X_1\times\cdots \times X_d$ is called a $d$-{\it box}.
A non-empty set $K \subseteq X$ is called a \textit{ box} if $K=K_1\times\cdots \times K_d$ and
$K_i\subseteq X_i$ for each $i\in [d]$. By  $\Pud {X}$ we denote the set of all boxes in $X$. 

\vspace{-0mm}
{\center
\includegraphics[width=12cm]{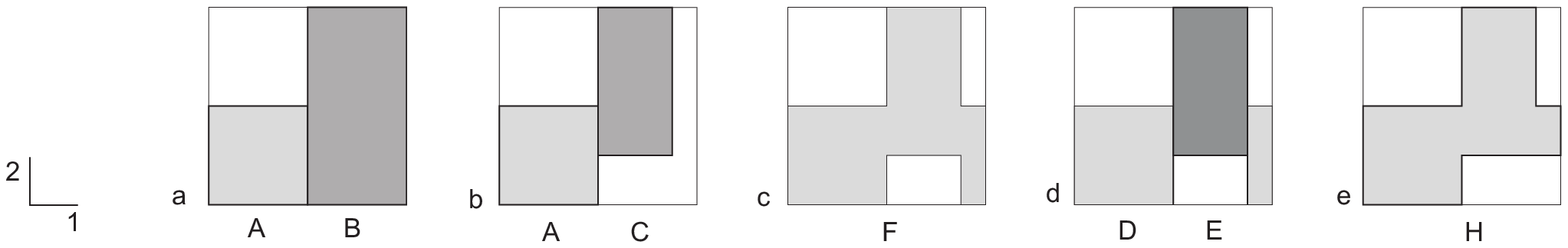}\\
}

\smallskip
\noindent{\footnotesize Fig. 2. The box $A\subset [0,1]^2$ is proper, and $B$ is not. The boxes $A$ and $B$ are dichotomous, while $A$ and $C$ are not. The set $F\subset [0,1]^2$ is a polybox, and $\ka F=\{D,E\}$ is suit for it. Moreover, $F$ is rigid. The set $H$ is not a polybox. 
}

\medskip 
%The set of all boxes in $X$ is denoted by $\Pud X$. 
The box  $K$ is said to be \textit{ proper} if $K_i\neq X_i$ for each $i\in [d]$.  
Two boxes $K$ and $G$ in $X$ are called \textit{dichotomous}
if there is $i\in [d]$ such that $K_i=X_i\setminus G_i$. A \textit{suit} is any collection of pairwise
dichotomous boxes. A suit is \textit{proper} if it consists of proper boxes. 
A non-empty set $F\subseteq X$ is said to be a \textit{ polybox} if
there is a suit $\ka F$ for $F$, i.e. if $\bigcup \ka F=F$. In other words, $F$ is a polybox if it has a partition into pairwise dichotomous boxes. 
A polybox $F$ is {\it rigid} if it has exactly one suit, that is if $\ka F$ and $\ka G$ are suits for a rigid polybox, then $\ka F=\ka G$. 
(Figures $2c$ and $3d,e$; the polyboxes  $\bigcup \ka F^{3,A}$ and $\bigcup \ka F^{3,A'}$ in Figure 5 are not rigid.) 

\vspace{-2mm}
{\center
\includegraphics[width=14cm]{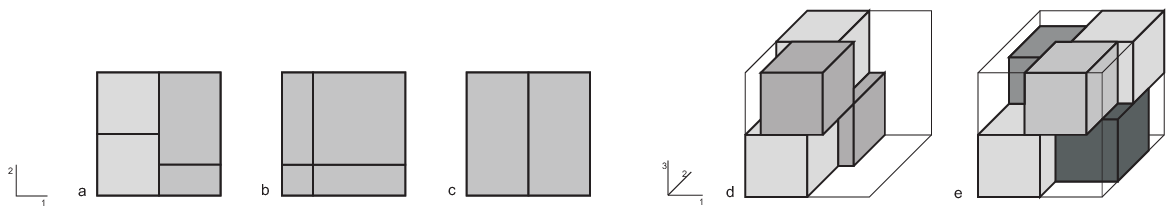}\\
}

\smallskip
\noindent{\footnotesize Fig. 3. Partition $a$ is a minimal partition of $[0,1]^2$, and Partition $b$ is a minimal partition which is a simple partition. Partition $c$ is a simple partition which is not a minimal partition. The systems of boxes $d, e$ are two suits for rigid polyboxes in the $3$-box $[0,1]^3$. %The two brightest boxes are $[0,1)^3$ and $[1,2]^3$. The remaining three boxes are $[3/4,7/4)\times [0,1)\times [1,2]$ (the light box), $[0,1)\times [1,2]\times [3/4,7/4)$ (the darker box) and $[1,2]\times [3/4,7/4)\times [0,1)$ (the darkest box).
}

\medskip
The important property of proper suits is that, for every proper suits $\ka F$ and $\ka G$ for a polybox $F$, we have $|\ka F|=|\ka G|$ (see the suits  $\ka F^{3,A}$ and $\ka F^{3,A'}$ in Figure 5). Thus, we can define a {\it box number} $|F|_0=$ the number of boxes in any proper suit for the polybox $F$ (compare (2.4) and Theorem 2.4 in \cite{KP}). In Figure 5 we have $|\bigcup \ka F^{3,A}|_0=3$. Obviously, the above property is not true for suits which are not proper (see Figure $7a$).
 A proper suit for a $d$-box $X$ is called a {\it minimal partition } of $X$ (Figures 3 and 5). In \cite{GKP} we showed that
 \begin{tw}
 \label{min}
 A suit  $\ka F$ is a minimal partition of a $d$-box $X$ if and only if $|\ka F|=2^d$.
 \end{tw}
 
 %Every minimal partition of a $d$-box has $2^d$ boxes. 

A family $\ka C\subset \Pud X$ is called a {\it simple partition} of $X$ if for every $K, G\in \ka C$ and every $i\in [d]$ we have $K_i=G_i$ or, if $G_i\neq X_i$, $K_i=X_i\setminus G_i$ and $\ka C$ is a suit for $X$ (Figures $3b,c$). 

Two boxes $K,G\subset X$ are said to be a {\it twin pair} if $K_j=X_j\setminus G_j$ for some $j\in [d]$ and $K_i=G_i$ for every $i\in [d]\setminus\{j\}$. Alternatively, two dichotomous boxes $K,G$ are a twin pair if $K\cup G$ is a box. (In Figure 3, Partitions $a,b,c$ contains twin pairs, while the suits in Figure $3d,e$ do not contain a twin pair).
Observe that the suit for a rigid polybox cannot contain a twin pair.

%It should be explicitly emphasized that 
%We shall base our considerations on the methods and notions developed there. 

%One of such partition is a minimal partition whose properties were studied in \cite{GKP,KP}(compare also \cite{ABHK}). Here is an example of such a partition and its connection with a cube tiling.       

\subsection{The structure of minimal partitions}
In order to sketch our approach to the problem of the existence of twin pairs in a cube tiling of $\er^d$, we describe the structure of a minimal partition.
A graph-theoretic description of this structure can be found in \cite{CS1,La} (see also \cite{LS2}).   
  
Let $X$ be a $d$-box. A set $l_i=\{x_1\}\times \cdots \times \{x_{i-1}\}\times X_i \times \{x_{i+1}\}\times \cdots \times \{x_d\}$, where $x_j\in X_j$ for $j\in [d]\setminus\{i\}$, is called a {\it line } in $X$. A set $F\subseteq X$ is called an {\it $i$-cylinder} (Figure 4) if for every line $l_i$ 
one has 
$$
l_i\cap F=l_i \; \; \; {\rm or}\; \; \; l_i\cap F=\emptyset.
$$

\vspace{-0mm}
{\center
\includegraphics[width=8cm]{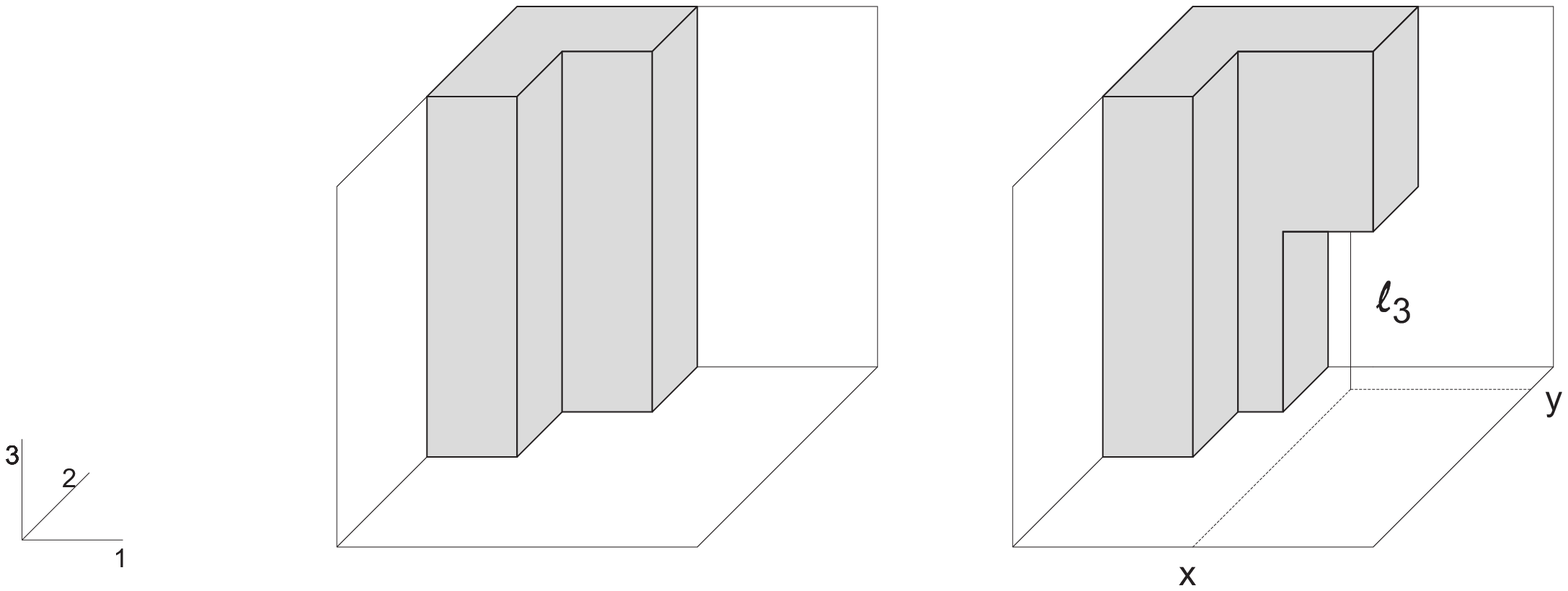}\\
}

\medskip
\noindent{\footnotesize Fig. 4. The set on the left is a $3$-cylinder in $X=[0,1]^3$, and the set on the right is not because the line $l_3=\{x\}\times \{y\}\times [0,1]$ has a non-empty intersection with this set but $l_3$ is not entire contained in it. 
}

\medskip
Let $\ka F$ be a minimal partition, and 
%It is clear that if the minimal partition $\ka F=\{([0,1]^d+x)\cap ([0,1)^d+t)\neq\emptyset:t\in T\}$ contains a twin pair, then the tiling   $[0,1)^d+T$ also contains a twin pair. Thus, it is sufficient to focus on twin pairs in $\ka F$. We can indicate sets of boxes in $\ka F$, in which twin pairs, if they exist, can be found.
let $B\subset X_i$ be a set such that there is a box $K\in \ka F$ with $K_i\in \{B, B^c\}$, where $B^c=X_i\setminus B$.  Let
$$
\ka F^{i,B}=\{K\in \ka F: K_i=B\}\quad {\rm and} \quad \ka F^{i,B^c}=\{K\in \ka F: K_i=B^c\}.
$$
 Since boxes in $\ka F$ are pairwise dichotomous, the set $\bigcup (\ka F^{i,B}\cup \ka F^{i,B^c})$ is an $i$-cylinder, and the set of boxes $\ka F^{i,B}\cup \ka F^{i,B^c}$ is a suit for it. As, by Theorem \ref{min}, $|\ka F|=2^d$, it follows that the boxes in $\ka F$ can form at most $2^{d-1}$ pairwise disjoint $i$-cylinders. More precisely, for every $i\in [d]$ there are sets $B^1,\ldots ,B^{k_i}\subset X_i$ such that $B^n\not\in \{B^m,(B^m)^c\}$ for every $n,m\in [k_i], n\neq m$, and 
\begin{equation}
\label{rep}
\ka F=\ka F^{i,B^1}\cup \ka F^{i,(B^1)^c}\cup \cdots \cup \ka F^{i,B^{k_i}}\cup \ka F^{i,(B^{k_i})^c}. 
\end{equation}
The boxes in $\ka F$ are proper, and hence $|\ka F^{i,B^n}\cup \ka F^{i,(B^n)^c}|\geq 2$. Thus, $k_i\leq 2^{d-1}$ for every $i\in [d]$. Observe that that is why $1\leq |L(T,x,i)|\leq 2^{d-1}$ for every cube tiling $[0,1)^d+T$, $x\in \er^d$ and $i\in [d]$. 

If $K$ is a box in $X$ and  $\ka G$ is a family of boxes, then let  
$$
K_{i^c}=K_1\times\cdots \times K_{i-1}\times K_{i+1}\times\cdots \times K_d\quad {\rm and} \quad \ka G_{i^c}=\{K_{i^c}\colon K\in \ka G\}.
$$
Since $\bigcup (\ka F^{i,B}\cup \ka F^{i,B^c})$ is an $i$-cylinder, the sets of boxes $\ka F^{i,B}_{i^c}$ and $\ka F^{i,B^c}_{i^c}$, where $\ka F^{i,B}_{i^c}=(\ka F^{i,B})_{i^c}$, are two suits for the polybox $\bigcup \ka F^{i,B}_{i^c}=\bigcup \ka F^{i,B^c}_{i^c}$, which is a polybox in the $(d-1)$-box $X_{i^c}$ (Figure 3). The sets $\ka F^{i,B}_{i^c}$ and $\ka F^{i,B^c}_{i^c}$ are proper suits for the polybox $\bigcup \ka F^{i,B}_{i^c}$ and therefore  $|\ka F^{i,B}_{i^c}|=|\ka F^{i,B^c}_{i^c}|$ 

{\center
\includegraphics[width=9cm]{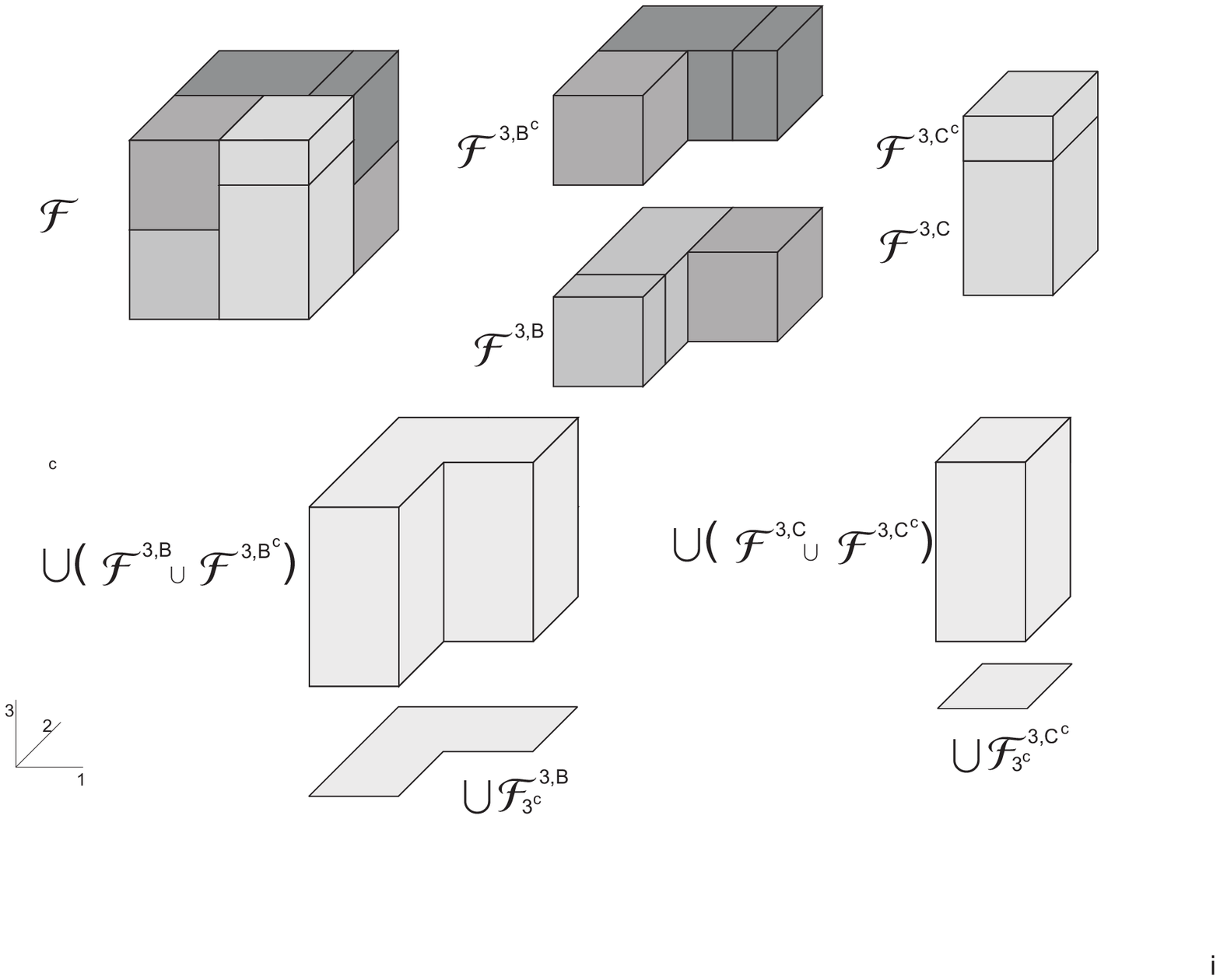}\\
}

\vspace{-6mm}
\noindent{\footnotesize 
Fig. 5. The minimal partition $\ka F=\ka F^{3,B}\cup \ka F^{3,B^c}\cup \ka F^{3,C}\cup \ka F^{3,C^c}$ of the $3$-box $X=[0,1]^3$ ($B=[0,1/2)$, $C=[0,3/4)$). The set $\bigcup (\ka F^{3,B}\cup \ka F^{3,B^c})$ is a $3$-cylinder and $\ka F^{3,B}\cup \ka F^{3,B^c}$ is a proper suit for it. 
  }

\medskip
Now, if $K,G\in \ka F$ are a twin pair, then there is a suit $\ka F^{i,B}\cup \ka F^{i,B^c}\subset \ka F$ for some $i$-cylinder  such that $K,G\in  \ka F^{i,B}\cup \ka F^{i,B^c}$. Thus, $K,G\in \ka F^{i,B}$ or $K,G\in \ka F^{i,B^c}$ or  $K\in \ka F^{i,B}$ and $G\in \ka F^{i,B^c}$. In the third case $K_{i^c}=G_{i^c}\in \ka F^{i,B}_{i^c}\cap \ka F^{i,B^c}_{i^c}$. So, if $\ka F^{i,B}_{i^c}\cap \ka F^{i,B^c}_{i^c}\neq \emptyset$, then there is a twin pair in $\ka F$. Obviously, there is exactly one $i\in [d]$ such that $K\in \ka F^{i,B}$ and $G\in \ka F^{i,B^c}$ (thus, $K,G\in \ka F^{j,B^j}$ or $K,G\in \ka F^{j,(B^j)^c}$ for every $j\in [d]\setminus\{i\}$). 

Now we can ask on the maximal positive integer $n$ such that if  $\ka F^{i,B}$ and $\ka F^{i,B^c}$ do not contain a twin pair and  $|\ka F^{i,B}|\leq n$, then   $\ka F^{i,B}_{i^c}=\ka F^{i,B^c}_{i^c}$ (which means that there is a twin pair $K,G$ such that  $K\in \ka F^{i,B}$ and $G\in \ka F^{i,B^c}$). This is rather obvious that if $m=|\ka F^{i,B}_{i^c}|$ is very small and $\ka F^{i,B}_{i^c}$ is a suit without twin pairs, then the polybox $\bigcup \ka F^{i,B}_{i^c}$ has exactly one suit which is equal to $\ka F^{i,B}_{i^c}$ (see Figure $2d$ for $m=2$ and Figure $3d$ for $m=4$). We will show in Section 4 that $n=11$, from where Theorem \ref{keli} will easily follow. 

\medskip
%In \ref{CS1} Corradi and Szabo consider the following graph approach to Keller's conjecture for cube tilings $[0,1)^d+T$ of $\er^d$, where $T\subset a+\zet^7\cup b+\zet^7$ (compare Section 1). 
We now discuss a connection between the structure of minimal partitions and a graph theoretical approach proposed by Lawrence \cite{La}, Corr\'adi and Szab\'o \cite{CS1}.  

\smallskip
Let $G_1,\ldots ,G_d$ be bipartite graphs with the same vertex set $V$, where $|V|=2^d$,  such that each of them is a disjoint union of regular complete bipartite graphs. Let $G_1\cup \cdots \cup G_d$ be the complete graph with the vertex set $V$.

We show that every minimal partition $\ka F$ induces in a natural way the graphs $G_1,\ldots ,G_d$ with the vertex set $V=\ka F$. 

Since $\ka F$ is a minimal partition, for every $i\in [d]$ the partition $\ka F$ has a representation of the form (\ref{rep}). Two vertices $K,G\in V$ are adjacent in  $G_i$ if there is $j\in [k_i]$ such that $K\in \ka F^{i, B^j}$ and  $G\in \ka F^{i, (B^j)^c}$ or $G\in \ka F^{i, B^j}$ and  $K\in \ka F^{i, (B^j)^c}$. Obviously, $G_i$ is a bipartite graph of $V$. Since for every $j\in [k_i]$ and every $K\in \ka F^{i,B^j}$ the box $K$ is dichotomous to every $G\in \ka F^{i,(B^j)^c}$ and vice versa, the graph $G_i$ decomposes into precisely $k_i$ pairwise disjoint regular complete bipartite graphs (recall that  $|\ka F^{i, B^j}|=|\ka F^{i, (B^j)^c}|$) each of which has the vertex set $\ka F^{i, B^j}\cup \ka F^{i, (B^j)^c}$ for $j\in [k_i]$ with the bipartition  $\{\ka F^{i, B^j},\ka F^{i, (B^j)^c}\}$. Since every two boxes in the minimal partition $\ka F$ are dichotomous, $G_1\cup \cdots \cup G_d$ is a complete graph with the vertex set $V=\ka F$ 

\smallskip
Conversely, every family $G_1,\ldots ,G_d$ of bipartite graphs on $V$ with the above properties induces a minimal partition of a $d$-box $X$. To show this, for every $i\in [d]$ let $X_i$ be an arbitrary set which has at least $k_i$ proper subsets $B^j$ such that $B^n\not\in \{B^m,(B^m)^c\}$ for every $n,m\in [k_i],n\neq m$, where $k_i$ is the number of all regular complete  bipartite graphs in $G_i$. If $G^1_i,\ldots ,G^{k_i}_i$ are the regular complete bipartite graphs in $G_i$, then let $\{(G^j_i)^+,(G^j_i)^-\}$ be the bipartition of $G^j_i$ for $j\in [k_i]$. For every $i\in [d]$ let $f_i:\{(G^j_i)^+,(G^j_i)^-:j\in [k_i]\}\rightarrow 2^X$ be an injection such that $f_i((G^j_i)^+)=X_i\setminus f_i((G^j_i)^-)$ 
for every $j\in [k_i]$. For every $v\in V$ there are unique $(\varepsilon_1,\ldots , \varepsilon_d)\in \{+,-\}^d$ and $(j_1,\ldots ,j_d)\in [k_1]\times\cdots \times [k_d]$ such that $v\in (G^{j_i}_i)^{\varepsilon_i}$ for every $i\in [d]$. Define a box $K_v$ in $X$ by the formula 
$$
K_v=f_1((G^{j_1}_1)^{\varepsilon_1})\times \cdots \times f_d((G^{j_d}_d)^{\varepsilon_d}).
$$ 
By the definition of $K_v$ for every  $v\in V$ the box $K_v$ is proper in $X$. Since $G_1\cup\cdots \cup G_d$ is a complete graph, for every two vertices $v,w\in V$ the boxes $K_v$ and $K_w$ are dichotomous, and since $|V|=2^d$, by Theorem \ref{min}, the family $\ka F_V=\{K_v:v\in V\}$ is a minimal partition. 

\smallskip
In the above graph approach Keller's conjecture says (in brackets the notation for dichotomous boxes is given):  There is $i\in [d]$ and an edge $e$ in $G_i$ (there are   $K\in \ka F^{i,B}$ and $G\in \ka F^{i,B^c}$) such that the graph $G_j\cup e$ contains a triangle for every $j\in [d]\setminus\{i\}$ ($K,G\in \ka F^{j,B^j}$ or $K,G\in \ka F^{j,(B^j)^c}$ for $j\in [d]\setminus\{i\}$) (\cite{CS1}).

\subsection{Cube tilings and dichotomous boxes}

Every two cubes $[0,1)^d+t$ and $[0,1)^d+p$ in an arbitrary cube tiling $[0,1)^d + T$ of $\er^d$ satisfy {\it Keller's condition}: There is $i\in [d]$ such that $t_i-p_i\in \zet\setminus\{0\}$, where $t_i$ and $p_i$ are $i$th coordinates of the vectors $t$ and $p$ (\cite{Ke1}). For any cube $[0,1]^d+x$, where $x=(x_1,...,x_d)\in \er^d$, the family $\ka F_x=\{([0,1)^d+t)\cap ([0,1]^d+x)\neq\emptyset:t\in T\}$ is a partition of the cube $[0,1]^d+x$, in which, because of Keller's condition, every two boxes $K,G\in \ka F_x$ are dichotomous, that is, there is $i\in [d]$ such that $K_i$ and $G_i$ are disjoint and $K_i\cup G_i=[0,1]+x_i$. Moreover, since cubes in cube tilings are half-open, every box $K\in \ka F_x$ is proper, and consequently the family $\ka F_x$ is a minimal partition. The structure of the partition $\ka F_x$ reflects the local structure of the cube tiling $[0,1)^d+T$. Obviously, a cube tiling $[0,1)^d+T$ contains a twin pair if and only if the partition $\ka F_x$ contains a twin pair for some $x\in \er^d$ (\cite{La,P}) (see Figure 1). %Knowing the structure of a minimal partition given in the previous section we can give an interpretations of the numbers $r^-(T)$ and $r^+(T)$. %Let $[0,1)^d+T$ be a cube tiling of $\er^d$, and   The number $r^+(T)$ is 
Observe also that if $\ka F_x=\ka F_x^{i,B^1}\cup \ka F_x^{i,(B^1)^c}\cup \cdots \cup \ka F_x^{i,B^{k_i(x)}}\cup \ka F_x^{i,(B^{k_i(x)})^c}$, then $|L(T,x,i)|=k_i(x)$ (compare (\ref{rep})).  

%\begin{tw}
%\label{corr}
%Let $G_1,\ldots ,G_d$ be bipartitie graphs of the same vertex set $V$ such that each of them is a disjoint union of at most two regular complete bipartite graphs. If $G_1\cup \cdots \cup G_d$ is the complete graph of vertex set $V$, then there is $i\in [d]$ and an edge $e$ of $G_i$ such that $e\cup G_j$ contains a triangle for every $j\in [d]\setminus \{i\}$.
%\end{tw}   

\vspace{-2mm}
\subsection{Dichotomous words and polybox codes}

The results in  the present paper are formulated and proved  in full generality. Suits have the form of systems of abstract words.   
We collect below  basic notions concerning words (details can be found in \cite{KP}).

A set $S$ of arbitrary objects will be called an \textit{alphabet}, and the elements of $S$ will be called \textit{letters}. A permutation $s\mapsto s'$ of the alphabet $S$ such that $s''=(s')'=s$ and $s'\neq  s$ is said to be a \textit{complementation}. We add a special letter $\ast$ to the set $S$ and the set $S\cup \{*\}$ is denoted by $*S$. We set $*'=*$ and the letter $*$ is the only letter with this property (compare also the beginning of the next section). 
%(At the beginning of the next section we give another example of the differences of the letter $*$ from all the letters in $S$.) 
Each sequence of letters $s_1\ldots s_d$ from the set $*S$ is called a \textit{word}. The set of all words of length $d$ is denoted by $(*S)^d$, and by $S^d$ we denote the set of all words $s_1\ldots s_d$ such that $s_i\neq *$ for every $i\in [d]$.  Two words $u=u_1\ldots u_d$ and $v=v_1\ldots v_d$ are \textit{dichotomous} if there is $j\in[d]$ such that  $u_j\neq *$ and $u'_j=v_j$. If $V\subset (*S)^d$ consists of pairwise dichotomous words, then we call it a \textit{polybox code} (or \textit{polybox genome}). (In the next section we give examples of polybox codes and their relationships with suits.)
A pair of words $u,v\in (*S)^d$ is \textit{a twin pair } if there is $j\in[d]$ such that $u'_j=v_j$, where $u_j\neq *$ and  $u_i=v_i$ for every $i\in [d]\setminus\{j\}$.

If $A\subseteq [d]$ and $A=\{i_1<\cdots <i_n\}$, then $u_{A}=u_{i_1}\ldots u_{i_n}$ and  $V_{A}=\{v_{A}:v\in V\}$ for $V\subset (*S)^d$. If $\{i\}^c=[d]\setminus \{i\}$, then we write $u_{i^c}$ and $V_{i^c}$ instead of $u_{\{i\}^c}$ and $V_{\{i\}^c}$, respectively. If $V\subset (*S)^d$, $s\in *S$ and $i\in [d]$, then let $V^{i,s}=\{v\in V:v_i=s\}$. The representation
\begin{equation}
\label{rep1}
V=V^{i,l_1}\cup V^{i,l_1'}\cup \ldots \cup V^{i,l_{k_i}}\cup V^{i,l_{k_i}'},  
\end{equation}
where  $l_j,l_j'\in *S$ and $V^{i,l_j}\cup V^{i,l_j'}\neq\emptyset$ for $j\in [k_i]$, will be called a {\it distribution of words in} $V$.  

%Let $S=\{a,a',b,b'\}$ and let $V\subset S^d$ be a polybox codes. 
Let us discuss briefly a connection between dichotomous words and adjacent vertices in a $d$-dimensional Keller graph (see Section 1). Recall that two vertices $v$ and $w$ in the $d$-dimensional Keller graph on the vertex set $\{0,1,2,3\}^d$ are adjacent if there are $i,j\in [d]$, $i\neq j$, such that $v_i\neq w_i, v_j\neq w_j$ and $|v_i-w_i|=2$ or $|v_j-w_j|=2$. Define a complementation on the alphabet $\{0,1,2,3\}$ by $0'=2$ and $1'=3$. Thus, two vertices $v,w\in \{0,1,2,3\}^d$ are adjacent in the Keller graph if and only if the words $v,w$ are dichotomous and they do not form a twin pair. %Since $S=\{0,1,2,3\}$ is now an alphabet with a complementation we can consider a polybox code $V\subset S^d$. Of course, we allow that $V$ contains a twin pair $u,q$ and  when this words are treated as a vertices in a Keller graph, then $u$ and $q$ are not adjacent. %This is one of the reasones why we do not consider our polybox codes as a subgraphs of Keller graph. 
In the paper we consider polybox codes $V$ whose words are written down in an alphabet $S$ which has more then four letters and therefore the elements of $V$ cannot be considered as vertices of the Keller graph. But when $V\subset S^d$, where $S=\{a,a',b,b'\}$, the reader who is familiar with the Keller graphs may assume that $0=a,2=a',1=b$ and $3=b'$.

\subsection{Realizations of polybox codes}
Let $X=X_1\times \cdots \times X_d$ be a $d$-box.
Suppose that for each $i\in[d]$ a mapping $f_i\colon *S\to \Pud {X_i}$ is such that $f_i(s)\neq X_i$ for $s\in S$, $f_i(s')=X_i\setminus f_i(s)$ for $s\neq *$ and $f_i(*)=X_i$. (This is why the letter $*$ is special.) We define the mapping $f\colon (*S)^d\to \Pud  X$ by
$$
f(s_1\ldots s_d)=f_1(s_1)\times\cdots\times f_d(s_d).
$$ 

About such defined function $f$ we will say that it \textit{preserves dichotomies}. Thus, if $v\in (*S)^d$ and $v_i=*$ for some $i\in [d]$, then $f(v)$ is not a proper box, while for every $v\in S^d$ the box $f(v)$ is proper.  

If $V\subset (*S)^d$, then the set of boxes $f(V)=\{f(v)\colon v\in V\}$ is said to be a \textit{realization} of the set of words $V$. For example, the family $\{A,B\}$ in Figure $2a$ is a realization $f(V)$ of the code $V=\{aa,a'*\}$, where $f_1(a)=[0,1/2),f_2(a)=[0,1/2)$ and $f_2(*)=[0,1]$ (for more sophisticated examples of realizations of polybox codes see Figures 6,7 and Example \ref{p}). 

Clearly, if $V$ is a polybox code, then $f(V)$ is a suit for the polybox $\bigcup f(V)$. The realization is said to be \textit{exact} if for each pair of words  
$v,w \in V$,  if  $v_i\not\in \{w_i, w'_i\}$, then $f_i(v_i)\not\in\{f_i(w_i), X_i\setminus f_i(w_i)\}$ (Figure 6). 

A polybox code $V\subset (*S)^d$ is called  a {\it partition code } if any realization $f(V)$ of $V$ is a suit for a $d$-box $X$. For example, $W=\{a*,a'*\}$ is a partition code (compare Figure $3c$ for $f(W)$, where $f_1(a)=[0,1/2)$), while $V=\{aa,a'*\}$ is not a partition code, where $a\in S$ (Figure $2a$ for $f(V)$). Observe that, if  $V\subset S^d$ is  a partition code, then $f(V)$ is a minimal partition. Indeed, since $v\in S^d$ for every $v\in V$, the box $f(v)$ is proper and thus $f(V)$ is a proper suit for $X$. Moreover, if a partition code $V\subset S^d$ has a distribution of words of the form (\ref{rep1}), and $\ka F$ is an exact realization of $V$, then for every $j\in [k_i]$ the set $\ka F^{i,B^j}\cup \ka F^{i,(B^j)^c}$ is an exact realization of the polybox code $V^{i,l_j}\cup V^{i,l_j'}$, where   $\ka F^{i,B^j}\cup \ka F^{i,(B^j)^c}$, for $j\in [k_i]$,  are such as in (\ref{rep}). 

A partition code $C\subset (*S)^d$ is said to be  {\it simple} if 
for every $v,w\in C$ and every $i\in [d]$ we have $v_i=w_i$ or $v_i=w'_i$. For examples, the codes  $W=\{a*,a'*\}$ and $U=\{ab,a'b,ab',a'b'\}$, where $a,b\in S$, are simple, while $V=\{aa,aa',a'b,a'b'\}$ is not a simple code (see Figures $6a,b,c$ for $f(V)$). 
%If for every $i\in [d]$ and every $v=v_1\cdots v_d\in C$ we have $v_i\neq *$ then $C$ is called  {\it a proper} simple partition genome. 

We will exploit some abstract but very useful realization of polybox codes. This sort of realization was invented in \cite {ABHK}, where it was the crucial tool in proving the main theorem of that paper. 
   
Let $S$ be an alphabet with a complementation, and let
$$
 ES=\{B\subset S\colon |\{s,s'\}\cap B|=1, \text{whenever $s\in S$}\},
$$
$$ 
E s=\{B\in ES\colon s\in B\}\;\; {\rm and}\;\; E*=ES.
$$   
Let $V\subset (*S)^{d}$ be a polybox code, and let $v\in V$. The {\it equicomplementary} realization of the word $v$ is the box   
$$
\breve{v}=Ev_1\times \cdots \times Ev_d
$$
in the $d$-box $(ES)^d =ES\times \cdots \times ES.$ The equicomplementary realization of the code $V$ is the family
$$
E(V)=\{\breve{v}:v\in V\}.
$$ 
If $S$ is finite, $s_1,\ldots , s_n\in S$ and $s_i\not\in\{s_j, s'_j\}$ for every $i\neq j$, then
\begin{equation}
\label{dkostki}
|E s_1 \cap\dots\cap E s_n|=(1/2^{n})|ES|.
\end{equation}
In the paper we will assume that $S$ is finite, unless it will be explicitly stated otherwise. (Since we are working with a finite number of codes $V^1,\ldots ,V^k\subset S^d$, 
%and each of them has at most $2^d$ words, 
 if $S$ has infinite number of elements,  the set $S_1\subset S$ consisting of all the letters that appear in the words from $V^1\cup \cdots \cup V^k$ is finite because $|V^i|\leq 2^d$ for $i\in [k]$. Thus, in that case we can restrict ourselves to the finite alphabet $S_1$ and use (\ref{dkostki}) with $S_1$ insted of $S$.)

%of pairwise dichotomous words of length $d$, they always have at most $2^d$ elements. Thus,   

The value of the realization $E(V)$, where $V\subset S^d$,  lies in the above equality.  In particular, boxes in $E(V)$ are of the same size: $|E v_i|=(1/2)|ES|$ for every $i\in [d]$ and consequently  $|\breve{v}|=(1/2^d)|ES|^d$ for $v\in E(V)$. Thus, two boxes $\breve v, \breve w\subset (ES)^d$ are dichotomous if and only if $\breve v\cap \breve w=\emptyset.$

Moreover, from (\ref{dkostki}) we obtain the following important lemma. 

\vspace{-1mm}
\begin{lemat}
\label{=c}
Let $w,u,v\in S^d$, and let $\ka D$ be a simple partition of the $d$-box $\breve w$. If boxes $\breve w\cap \breve u$ and  $\breve w\cap \breve v$ belong to $\ka D$, then there is a simple partition code $C\subset S^d$ such that $u,v\in C$. In particular, if $\breve w\cap \breve u$ and  $\breve w\cap \breve v$ form  a twin pair, then $u$ and $v$ are a twin pair. 
\end{lemat}

\proof Assume on the contrary that $u$ and $v$ do not belong to any simple partition code $C\subset S^d$. Then there is $j\in [d]$ such that $u_j\not\in \{v_j,v_j'\}$. Since $\breve w\cap \breve u, \breve w\cap \breve v\in \ka D$, for every $i\in [d]$ we have $Ew_i\cap Eu_i=Ew_i\cap Ev_i$ or, if $Ew_i\neq Ev_i$, $Ew_i\cap Eu_i=Ew_i\setminus Ew_i\cap Ev_i$. 

If $Ew_j\cap Eu_j=Ew_j\cap Ev_j$, then   $Ew_j\cap Eu_j\cap Ev_j=Ew_j\cap Ev_j$.  On the other hand, since $u_j\not\in \{v_j,v_j'\}$ and $Ew_j\cap Eu_j=Ew_j\cap Ev_j$,  it follows that $Ew_j\neq Ev_j$ and $Ew_j\neq Eu_j$. Therefore, by (\ref{dkostki}),  $|Ew_j\cap Eu_j\cap Ev_j|=(1/8)|ES|$, while $|Ew_j\cap  Ev_j|=(1/4)|ES|$, a contradiction. 
%Thus,  $Ew_j\cap Eu_j\neq Ew_j\cap Ev_j$

Let now $Ew_j\cap Eu_j=Ew_j\setminus Ew_j\cap Ev_j$.  Clearly, $Ew_j\cap Ev_j\neq\emptyset$ and $Ew_j\cap Eu_j\neq\emptyset$, and then $Ew_j\neq Ev_j$ and $Ew_j\neq Eu_j$. Since   $u_j\not\in \{v_j,v_j'\}$, again by (\ref{dkostki}),  $|Ew_j\cap Eu_j\cap Ev_j|=(1/8)|ES|$. Thus, $Ew_j\cap Eu_j$ and $Ew_j\cap Ev_j$ are not disjoint, a contradiction.

%In the second case the non-empty sets $Ew_j\cap Eu_j$ and $Ew_j\cap Ev_j$ are disjoint, which is impossible, since, again by (\ref{dkostki}),  $|Ew_j\cap Eu_j\cap Ev_j|=(1/8)|ES|$. Thus, $v,u\in C$ for some simple partition code $C\subset S^d$.

If  $\breve w\cap \breve u$ and  $\breve w\cap \breve v$ form  a twin pair, then $Ew_j\cap Eu_j=Ew_j\setminus Ew_j\cap Ev_j$, $Ew_j\neq Ev_j$, for exactly one $j\in [d]$ and $Ew_i\cap Eu_i=Ew_i\cap Ev_i$ for every $i\in [d]\setminus\{j\}$. Therefore, by the first part of the lemma, $u,v\in C$ for some simple partition code $C\subset S^d$. Thus, $u_j=v_j'$ and $u_i=v_i$ for every $i\in [d]\setminus\{j\}$. 
\hfill{$\square$}

\medskip
In a general case the above lemma is not true; for a given three boxes $K,G$ and $H$ in a $d$-box $X$ such that $K$ and $G$ are dichotomous and $K\cap H$, $G\cap H$ belong to a simple partition of $H$, it can happen that $K$ and $G$ are not members of the same simple partition of $X$. For example, let $X=[0,4]\times [0,4]$, $K=[2,4]\times [2,4]$, $G=[2,3)\times [0,2)$ and $H=[1,3)\times [1,3)$. Then  the boxes $K\cap H=[2,3)\times [2,3)$,  $G\cap H=[2,3)\times [1,2)$ belong to a simple partition $\ka C=\{[1,2)\times [1,2), [1,2)\times [2,3), [2,3)\times [1,2), [2,3)\times [2,3)\}$ of $H$. But $K$ and $G$ do not belong to  the same simple partition of $X$ because $K_1\not \in \{G_1,X_1\setminus G_1\}$.    

\medskip
In our proofs we will need a realization with `good' properties such as these given in Lemma \ref{=c}. The realization $E(V)\subseteq (ES)^d$ of $V$ is the best possible for our purposes what will be shown in Section 2.8. 

\smallskip
Let $V\subset (*S)^d$ be  a polybox code, and let $f(V)$ be an exact realization of $V$. The set $f(V)$ is a suit (a set of pairwise dichotomous boxes), while $V$ describes the structure of it. The code $V$ has infinitely many exact realizations which may be very different from each other.   
For example, the partitions in Figures $6a,b,c$ are pairwise different but they are all the  exact realizations of the polybox code $V=\{aa,aa',a'b,a'b'\}$. The differences can even be related to the number of dimensions of a specific partition; the sets in Partition $6c$ are $3$-dimensional, but this partition can be regarded as 2-dimensional minimal partition with the same structure as Partitions $6a$ and $6b$.

{\center
\includegraphics[width=12cm]{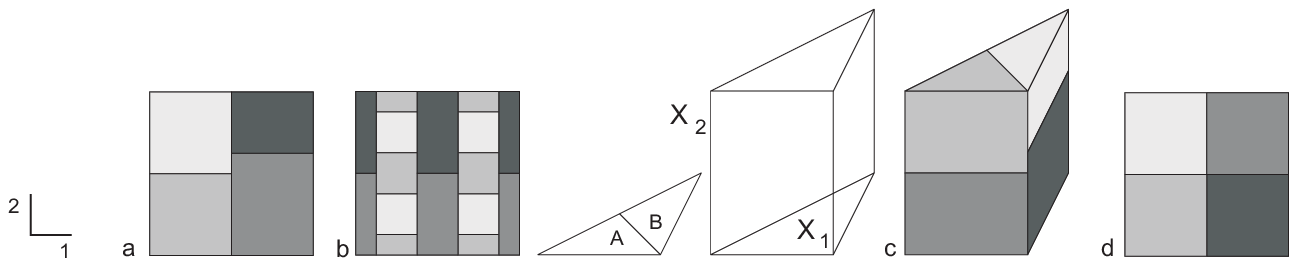}\\
}

\medskip
\noindent{\footnotesize Fig. 6. Partitions $a-c$ are exact realization of the code $V=\{aa,aa',a'b,a'b'\}$. Partition $d$ is also a realization of $V$ but not exact. Partition $c$ is a minimal partition of the $2$-box $X=X_1\times X_2$, where $X_1$ is the triangle $A\cup B$, and $X_2=[0,1]$. The realizations $a,b,d$ are partitions of $2$-box $[0,1]^2$; in the first realization we have $f_1(a)=f_2(a)=[0,1/2)$ and $f_2(b)=[0,5/8)$; in the second one we have $f_1(a)=[1/8,3/8)\cup (5/8,7/8), f_2(a)=[0,1/2)$ and  $f_2(b)=[1/8,3/8)\cup (5/8,7/8)$; in the third case $f_1(a)=A$ (and then $f_1(a')=B$) , $f_2(a)=[0,1/2)$ and $f_2(b)=[0,5/8)$; in the last partition we have $f_i(a)=f_i(b)=[0,1/2)$ for $i=1,2$.    
  }
 
\subsection{Equivalent and rigid polybox codes}
Let $V,W\subset (*S)^d$ be polybox codes, and let $v\in (*S)^d$. We say that $v$ is \textit{covered} by $W$, and write $v\sqsubseteq W$, if $f(v)\subseteq \bigcup f(W)$ for every mapping $f$ that preserves dichotomies. For example, the word $v=bb$ is not covered by the code $W=\{aa,a'*\}$ as there is  a realization $f(W)$ such that $f(bb)\not\subset \bigcup f(W)$. For example, if $X=[0,1]^2$, $f_1(a)=f_2(a)=[0,1/2)$ (clearly, $f_2(*)=[0,1]$) and $f_1(b)=[1/3,2/3),f_2(b)=[0,2/3)$, then $f(bb)=[1/3,2/3)\times [0,2/3)\not\subset [0,1/2)^2\cup [1/2,1]\times [0,1]$ (obviously, one can find a realization $g(W)$ such that $g(v)\subset \bigcup g(W)$). It can be easily checked that for every $l\in *S$ the word $w=la$ is covered by $W$ (see also Example \ref{p}). 

If  $v\sqsubseteq W$ for every $v\in V$, then we write $V \sqsubseteq W$.

Polybox codes $V,W\subset (*S)^d$ are said to be {\it equivalent} if $V \sqsubseteq W$ and  $W \sqsubseteq V$ (Figure 7).  Thus, $V$ and $W$ are equivalent if and only if $\bigcup f(V)=\bigcup f(W)$ for every mapping $f$ that preserves dichotomies.

A polybox code $V\subset S^d$ is called {\it rigid} if there is no code $W\subset S^d$ which is equivalent to $V$ and $V\neq W$ (Example \ref{char}). Thus, a polybox code $V\subset S^d$ is rigid if and only if for every exact realization $f(V)$ of $V$ the polybox $\bigcup f(V)$ is rigid.
Observe that, rigid polybox codes cannot contain a twin pair. 

The following two results describe one of the most important property of polybox codes.
% are  is strongly connected to the notion of the rigidity of polyboxes and polybox codes.

Let $V\subset (*S)^d$, $|V|\geq 2$, be a partition code. %and let $v\in V$ be a word with the smallest number of the letter $*$. 
It follows from [Lemma 8.1, \cite{KP}] (see also  [Lemma 2.1, \cite{LS2}]) that there is a simple partition code $C\subset (*S)^d$ and two words $v,w\in V$ such that $v,w\in V\cap C$ and
\begin{equation}
\label{0indd}
|\{i\in [d]: v_i=w'_i,\;v_i\neq *\}|\equiv 1\;\; ({\rm mod}\; 2).
\end{equation}
Similarly, it follows from [Theorem 10.6, \cite{KP}] that, if $v\in S^d$ and $W\subseteq S^d$ is a polybox code such that $v\sqsubseteq W, v\not\in W$, then there is  a simple partition code $C\subset S^d$ and there are two words $w,u\in W\cap C$ such that 
\begin{equation}
\label{00indd}
|\{i\in [d]: w_i=u'_i\}|\equiv 1\;\; ({\rm mod}\; 2).
\end{equation}
(see Figure 12 for (\ref{0indd}) and Examples \ref{char} and \ref{p} for (\ref{00indd})).
Observe that it follows from the above that if a polybox codes $W\subseteq S^d$ is not rigid (which in particular means that $v\sqsubseteq W$ and $v\not\in W$ for some $v\in S^d$), then $W$ has to contain the above described words $w$ and $u$. 

\vspace{-0mm}
{\center
\includegraphics[width=12cm]{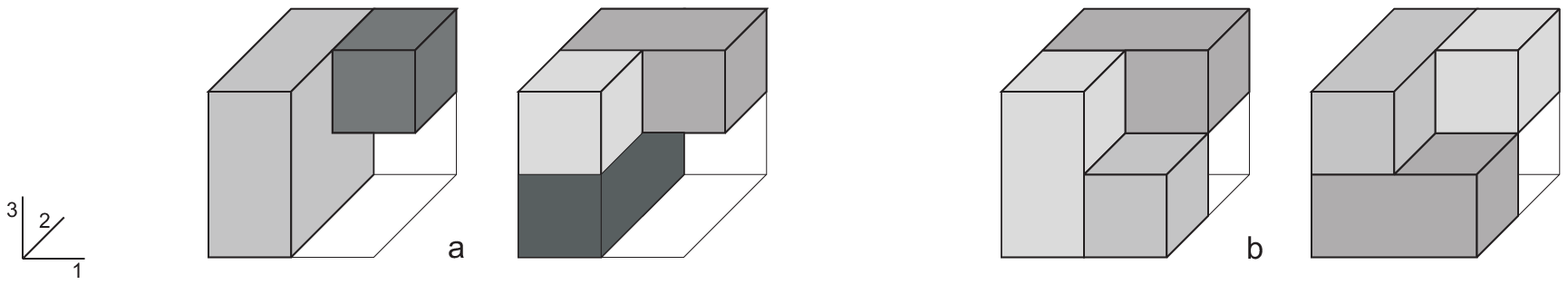}\\
}

\medskip
\noindent{\footnotesize Fig. 7. Two examples of realizations $f(V), f(W)$ of equivalent polybox codes in $X=[0,1]^3$. In $(a)$ we have the realizations of equivalent polybox codes $V$ and $W$, where $V=\{l_1**,l_1'l'_2l'_3\}$ (on the left in $(a)$) and $W=\{l_1*l_3,l_1l_2l'_3,*l'_2l'_3\}$ (on the right in $(a)$). In $(b)$ we have the realizations of equivalent polybox codes $U$ and $Q$, where $U=\{l_1l_2*,*l_2'l'_3,l_1'l_2l_3\}$ (on the left in $(b)$) and $Q=\{*l_2l_3,l_1*l'_3,l_1'l'_2l'_3\}$ (on the right in $(b)$). In these realizations we have $f_i(l_i)=[0,1/2)$ and $f_i(*)=[0,1]$ for $i=1,2,3$.
}

\medskip
The definition of the relation $\sqsubseteq$ is rather cumbersome tool to decide whether $w\sqsubseteq V$. Below we give very useful and easy test, especially in the computations, to check whether  $w\sqsubseteq V$. 

\medskip
Let  
$g\colon S^d\times S^d\to \zet$ be defined by the formula
\begin{equation}
\label{g}
g(v,w)=\prod^d_{i=1}(2[v_i=w_i]+[w_i\not\in\{v_i, v'_i\}]),
\end{equation}
where $\iver {p}=1$ if the sentence $p$ is true and $\iver {p}=0$ if it is false.

Let $w\in S^d$, and let $V\subset S^d$ be a polybox code. Then 

\begin{equation}
\label{2d}
\breve{w}\subseteq \bigcup E(V)\Leftrightarrow w \sqsubseteq V \Leftrightarrow \sum_{v\in V} g(v,w)= 2^d.
\end{equation}

It follows from the definition of equivalent polybox codes $V,W\subseteq S^d$ and (\ref{2d}) that $V$ and $W$  are equivalent if and only if $\bigcup E(V)=\bigcup E(W)$. Another characterization of equivalent polybox codes $V,W\subseteq S^d$ which stems from $(\ref{2d})$ is the following: Polybox codes $V,W\subset S^d$  are equivalent if and only if $\sum_{v\in V} g(v,w)= 2^d$  for every $w\in W$ and $\sum_{w\in W} g(w,v)= 2^d$ for every $v\in V$. 
\begin{ex}
\label{char}
{\rm Let $V=\{aaaa,a'a'a'a,baa'a,a'baa,aa'ba,bbba'\}$. If $l\not \in \{a,a'\}$, then  
$\sum_{v\in V} g(v,bbbl)= 1+1+2+2+2+8=2^4$, and if $l=a$, then $\sum_{v\in V} g(v,bbbl)= 2+2+4+4+4+0=2^4$. Therefore, for every $l\in S$, by (\ref{2d}), $bbbl\sqsubseteq V$. (Thus, $V$ has to contain two words $v,u$ described in (\ref{00indd}). These are $aaaa$ and $a'a'a'a$.) In particular, for every $l\in S$ the twin pair $bbbl,bbbl'$ is covered by $V$. Obviously, again by (\ref{2d}), this means that for every realization $f(V)$ the boxes $f(bbbl)$ and $f(bbbl')$ are contained in $\bigcup f(V)$. It can be checked, that $ \sum_{v\in V} g(v,w)< 2^d$ for every $w\in S^d$ such that $w\not\in \bigcup_{l\in S}\{bbbl,bbbl'\}$. Thus, every such word $w$ is not covered by $V$. Occasionally, we will denote this fact by $w\not \sqsubseteq V$. Moreover, it follows from the above that $V$ is rigid.  
}
\end{ex}
  
%Clearly, by (\ref{2d}), $V$ is rigid if and only if the polybox $\bigcup E(V)$ is rigid.

%Recall that, a polybox $F$ is rigid if it has only one suit. Thus, the exact realization of a rigid polybox code is the sole suit for a rigid polybox. 

\medskip 
Let $s_*=*\cdots *\in (*S)^d$ and let $\bar{g}(\cdot,s_*)\colon (*S)^d\to \zet$ be defined as follows:
$$
\bar{g}(v,s_*)=\prod^d_{i=1}(2[v_i=*]+[v_i\neq*]).
$$

\begin{lemat}
\label{2dd}
Let $V\subset (*S)^d$ be a polybox code. The code $V$ is a partition code if and only if $\sum_{v\in V} \bar{g}(v,s_*)= 2^d$. 
\end{lemat}
\proof By (\ref{dkostki}), the definition of $\bar{g}$ and the equality $E*=ES$, we have $\bar{g}(v,s_*)|ES|^d/2^d=|\breve{v}|$.
If  $V$ is a partition code, then $\sum_{v\in V} |\breve{v}|=|ES|^d$, and thus  $\sum_{v\in V} \bar{g}(v,s_*)=2^d$. If  $\sum_{v\in V} \bar{g}(v,s_*)=2^d$, then $\sum_{v\in V} |\breve{v}|=|ES|^d$, which means that $V$ is a partition code.
\hfill{$\square$}

\begin{wn}
\label{uV}
Let $V\subset S^d$ be a polybox code and let $u\in S^d$. For every $v\in V$ let $\bar{v}\in (*S)^d$ be defined  in the following way: If  $v_i\neq u_i$, then $\bar{v}_i=v_i$, and if  $v_i=u_i$, then $\bar{v}_i=*$. Let $\breve{u}\cap \breve{v}\neq\emptyset$ for every $v\in V$. If $u\sqsubseteq V$, then $\bar{V}=\{\bar{v}:v\in V\}$ is a partition code. 
\end{wn}
\proof By (\ref{dkostki}),  $|\breve{v}\cap \breve{u}|=(1/2^k)(1/2^d)|ES|^d$, where $k=|\{i:v_i\neq u_i\}|$, and by the definition of the function $\bar{g}$, we have $1/2^k=(1/2^d)\bar{g}(\bar{v},s_*)$. The set $\{\breve{u}\cap \breve{v}:v\in V\}$ is a suit for the $d$-box $\breve{u}$. Thus,  $\sum_{v\in V} |\breve{v}\cap \breve{u}|=|\breve{u}|$. Therefore,  $\sum_{\bar{v}\in \bar{V}}(1/2^d)\bar{g}(\bar{v},s_*)(1/2^d)|ES|^d=(1/2^d)|ES|^d$, which gives $\sum_{\bar{v}\in \bar{V}}\bar{g}(\bar{v},s_*)=2^d$. By Lemma \ref{2dd}, $\bar{V}$ is a partition code.
\hfill{$\square$}

\medskip
Let $X=X_1\times \cdots \times X_d$ be a $d$-box. For every $i\in [d]$ let $S_i$ be the set of all pairs $(A,i)$, where $A$ is a proper subset of $X_i$.
With the set of all proper boxes in $X$ we associate the set of words $S^d$, where $S=\bigcup_{i=1}^dS_i$. We define a complementation $(A,i)\mapsto (A,i)'$ on $S$ by the formula $(A,i)'=(A^c,i)$, where $A^c=X_i\setminus A$. %For example, $[0,1]_i'=(1,2]$ for $X_i=[0,2]$, while $[0,1]_j'=(1,3]$ for $X_j=[0,3]$.
Obviously, if  $X_i$ is infinite for some $i\in [d]$, then $S$ is infinite.

Observe now that if  $\ka F$ is a suit for a polybox $F\subseteq X$, then the set of words $V=\{\bar{A}_1\ldots \bar{A}_d\colon  A_1\times \cdots \times A_d\in \ka F\}\subset (*S)^d$, where $\bar{A}_i=(A_i,i)$ if $(A_i,i)\in S$ and $\bar{A}_i=*$ if $A_i=X_i$ for $i\in [d]$, is a polybox code, and  the suit $\ka F$ is an exact realization of $V$. (Of course, this is one of the many ways of receiving polybox codes for suits.)
\begin{lemat}
\label{kodd}
Let $X$ be a $d$-box, $F\subseteq X$ be a polybox, and let $\ka F$ and $\ka G$ be suits for $F$. 
The polybox codes $V=\{\bar{K}_1\ldots \bar{K}_d: K_1\times \cdots \times K_d\in \ka F\}$ and  $W=\{\bar{G}_1\ldots \bar{G}_d: G_1\times \cdots \times G_d\in \ka G\}$ are equivalent.   
\end{lemat}   
\proof For $G=G_1\times \cdots \times G_d$ and $k\in [d]$ let $G^k=E\bar{G}_1\times \cdots \times E\bar{G}_k\times G_{k+1}\times \cdots \times G_d$. %We assume first that $\ka F$ and $\ka G$ are proper suits for the polybox $F$. 
We show that for every $G\in \ka G$ we have  
\begin{equation}
\label{si}
G^1\subseteq \bigcup_{K\in \ka F}K^1. 
\end{equation}
To do this, let $G\in \ka G$ be fixed.
Observe that for every $K\in \ka F$ such that  $G\cap K\neq\emptyset$, where  $G_1\neq K_1$, $K_1\neq X_1$, and every $x_{1^c}\in G_{1^c}\cap K_{1^c}$ there is $H\in \ka F$ such that $H_1=X_1\setminus K_1$ and $x_{1^c}\in H_{1^c}$. To see this assume first that $G_1\subset K_1$. Then there is  $y\in K$ with $y_1\in K_1\setminus G_1$ and $y_{1^c}=x_{1^c}$. The point $y$ has to belong to a box $Q\in \ka G$ such that $Q_1=X_1\setminus G_1$ because boxes in $\ka G$ are pairwise dichotomous. Then there is a point $z\in Q$ such that $z_1\in X_1\setminus K_1$ and $z_{1^c}=y_{1^c}$. Since $z_{1^c}\in K_{1^c}$ and $\ka F$ is a suit, there is $H\in \ka F$ such that $z\in H$ and $H_1=X_1\setminus K_1$. 

The case $G_1\not \subset K_1$ is considered in the very similar way: For every $x_{1^c}\in G_{1^c}\cap K_{1^c}$ we may choose $x\in G$ such that $x_1\in X_1\setminus K_1$. Since $\ka F$ is a suit, there is $H\in \ka F$ such that $x\in H$ and $H_1=X_1\setminus K_1$.  

\smallskip
1. If $\tilde{x}\in G^1$ is such that $\tilde{x}_{1^c}\in K_{1^c}$ for a box $K\in \ka F$ with $K_1=G_1$, then $\bar{x}\in K^1$ because $E\bar{G}_1=E\bar{K}_1$.

2.  If $\tilde{x}\in G^1$ is such that $\tilde{x}_{1^c}\in K_{1^c}$ for a box $K\in \ka F$ with $K_1=X_1$, then $\bar{x}\in K^1$ because $E\bar{G}_1\subseteq E\bar{K}_1=ES$.

3. If $\tilde{x}\in G^1$ is such that $\tilde{x}_{1^c}\in K_{1^c}\cap H_{1^c}$ for some $K,H\in \ka F$ with $H_1=X_1\setminus K_1$, then $\bar{x}\in K^1\cup H^1$ because $E\bar{G}_1\subseteq E\bar{K}_1\cup E\bar{H}_1=ES$.

Since each point in $G^1$ has one of the properties mentioned in the points $1$ $2$ and $3$, the proof of (\ref{si}) is completed. 

Let $\ka F^1=\{K^1\colon K\in \ka F\}$ and   $\ka G^1=\{G^1\colon G\in \ka G\}$. It follows from (\ref{si}) that  $\bigcup \ka G^1\subseteq \bigcup \ka F^1$. In the same way we show that $\bigcup \ka F^1\subseteq \bigcup \ka G^1$. Thus, $\bigcup \ka G^1=\bigcup \ka F^1$.

Let $1<k<d$ and assume that $\bigcup \ka F^k=\bigcup \ka G^k$, where $\ka F^k=\{K^k:K\in \ka F\}$ and   $\ka G^k=\{G^k:G\in \ka G\}$. Along the same lines as for $k=1$ we show that 
$$
G^{k+1}\subseteq \bigcup_{K\in \ka F}K^{k+1}
$$
from where we infer that $\bigcup \ka F^{k+1}=\bigcup \ka G^{k+1}$.
By induction, $\bigcup \ka F^d=\bigcup \ka G^d$. Since $\bigcup \ka F^d=\bigcup E(V)$ and $\bigcup \ka G^d=\bigcup E(W)$, by (\ref{2d}), $V\sqsubseteq W$ and $W\sqsubseteq V$. Thus, $V$ and $W$ are equivalent. 

%Therefore, $\bigcup f(V)= \bigcup f(W)$ for every $f$ preserving dichotomies (compare the beginning of this section). 

%\smallskip
%To prove the general case, for every $G\in \ka G$ let $A(G)=\{i\in [d]\colon G_i\neq X_i\}$. If $G$ is not proper box, then let $\ka G_G$ be any proper suit for the box $G$ such that $H_{A(G)}=G_{A(G)}$ for every $H\in \ka G_G$; if $G$ is proper, then $\ka G_G=\{G\}$. In the same way we define suits $\ka F_K$ for a box $K\in \ka F$. Since $\ka F$ and $\ka G$ are suits for $F$, it follows that $\tilde{\ka G}=\bigcup_{G\in \ka G}\ka G_G$ and $\tilde{\ka F}=\bigcup_{K\in \ka F}\ka F_K$ are proper suits for the polybox $F$.
%and $\bigcup_{G\in \ka G}\bigcup \ka G_G=\bigcup_{K\in \ka F}\bigcup \ka F_K$. 
%Let $W_{G}=\{\bar{H}_1\ldots \bar{H}_d\colon H\in \ka G_G\}$,  $V_{K}=\{\bar{H}_1\ldots \bar{H}_d\colon H\in \ka F_K\}$, and let $\tilde{W}=\bigcup_{G\in \ka G}W_G$,  $\tilde{V}=\bigcup_{K\in \ka F}V_K$.  By the first part of the proof we have 
%$$
%\bigcup f(\tilde{V})=\bigcup f(\tilde{W})
%$$
%for every $f$ that preserves dichotomies. Note that for every  $w=w(G)=\bar{G}_1\ldots \bar{G}_d\in W$ and  $v=v(K)=\bar{K}_1\ldots \bar{K}_d\in V$ we have 
%$$
%f(w)=\bigcup f(W_{G})\quad {\rm and} \quad f(v)=\bigcup f(V_{K}).
%$$

%Thus, $\bigcup f(V)=\bigcup f(\tilde{V})$ and $\bigcup f(W)=\bigcup f(\tilde{W})$ and consequently $\bigcup f(V)=\bigcup f(W)$ for every $f$ that preserves dichotomies. Hence, $V$ and $W$ are equivalent. 
\hfill{$\square$}

\subsection{Key result in proving Theorem \ref{keli}}

The proof of Theorem \ref{keli} is based on the following theorem. 

\vspace{-1mm}
\begin{tw}
\label{12}
If $V,W\subset S^d$ are equivalent polybox codes which do not contain twin pairs and $W\cap V=\emptyset$, then
$|V|\geq 12$.
\end{tw}

\vspace{-2mm}
Recall that $V$ and $W$ are equivalent if and only if $\bigcup f(V)=\bigcup f(W)$ for every mapping $f$ that preserves dichotomies (see Section 2.5.). Thus, the simplest reformulation of the above theorem in the terms of dichotomous boxes is the following: For every $d$-box $X$ and every two proper suits $\ka F$ and $\ka G$ which do not contain a twin pair, if $\bigcup \ka F=\bigcup \ka G$ and $\ka F\cap \ka G=\emptyset$, then $|\ka F|\geq 12$.   Equivalently: For every $d$-box $X$ and every polybox $F\subset X$, if there is a proper suit $\ka F$ for $F$ which does not contain a twin pair and $|F|_0\leq 11$ (that is $|\ka F|\leq 11$), then every proper suit $\ka G$ for $F$ contains a twin pair or $\ka G=\ka F$ (Corollary \ref{Sf}).
This result was  announced in Section 2.2. In Section 4 we give a short proof of it.

Let  $S=\{0,1,2,3\}$, and let $0'=2$ and $1'=3$ (compare Section 2.4.). Since now $S=\{0,1,2,3\}$ is an alphabet with a complementation, we can speak about equivalent cliques in a $d$-dimensional Keller graph in the sense of the definition from Section 2.6. 
Theorem \ref{12} for cliques in a $d$-dimensional Keller graph reads as follows: Every two equivalent cliques in a $d$-dimensional Keller graph with at most 11 vertices are equal. We prove it in the last section of the paper (Corollary \ref{clik11}) where we also discuss some issues related to the Keller graph. 

It should be emphasized that Theorem \ref{12}  applies to polybox codes whose words belong to $S^d$, where $S$ is an arbitrary finite alphabet with a complementation, while vertices in a $d$-dimensional Keller graph are elements of $\{0,1,2,3\}^d$.

\subsection{Geometry of dichotomous boxes}

In this section we describe the main techniques which are used in the paper, and which are based on the properties of the realization $E(V)$ (see (\ref{dkostki})). 

Usually we will consider two disjoint and equivalent polybox codes $V\subset (*S)^d$ and  $W\subset (*S)^d$ 
%(in Sections 4, 5 and 6 we will consider only the case  $V,W\subset S^d$). 
Recall, that polybox codes $V$ and $W$ are equivalent if and only if $\bigcup E(V)=\bigcup E(W)$, where $\bigcup E(V)=\bigcup_{\breve{v}\in E(V)}\breve{v}$.   Moreover, let us recall that if $i\in [d]$ and  $v\in (*S)^d$, then $v_i$ denotes the letter standing in $v$ at the $i$-th position, while $v_{i^c}\in (*S)^{d-1}$ is the word that arises from $v$ by skipping  the letter $v_i$ in the word $v$. If $V\subset (*S)^d$, then $V_{i^c}=\{v_{i^c}:v\in V\}$.

Our goal is to reveal a structure of $V$ and $W$ or estimate the cardinality of $V$ (if $V,W\subset S^d$ and  $V,W$ are equivalent, then $|V|=|W|$).  Below we describe the most important techniques applied in the paper.

\smallskip
$\bullet$ {\it The structure of $V$ from the suit for $\breve{w}$}. Let $w\sqsubseteq V$. Then $\breve{w}\subseteq \bigcup E(V)$ and the set of boxes $\ka F_w=\{\breve{w}\cap \breve{v}:v\in V\}$ is a suit for $\breve{w}$. In Example  \ref{p} we show what kind of information can be obtained from the structure of $\ka F_w$.

{\center
\includegraphics[width=10cm]{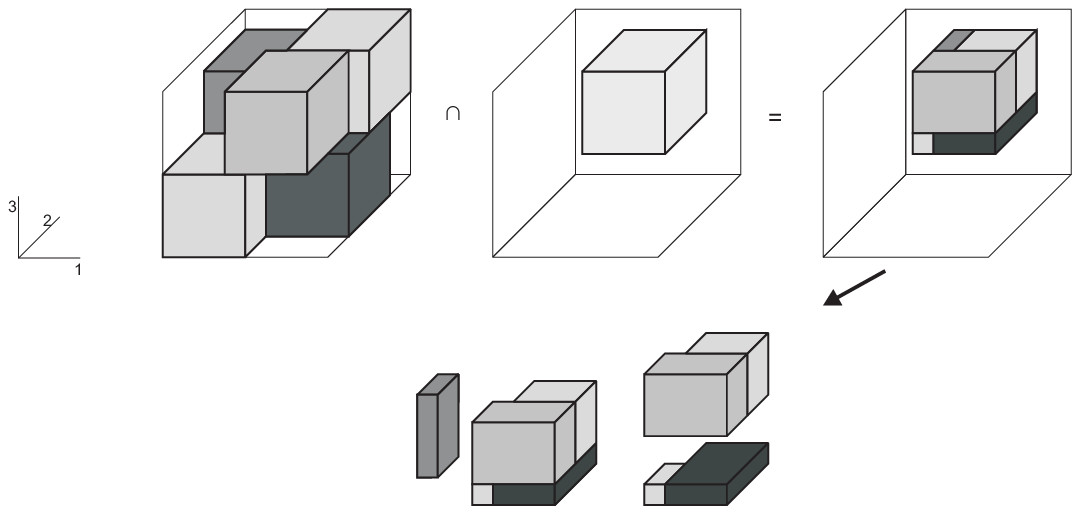}\\
}

\bigskip
\noindent{\footnotesize Fig. 8. On the top: The light box (in the middle) is contained in the sum of five pairwise dichotomous boxes (the boxes on the left). These boxes determine a partition of the light box into pairwise dichotomous boxes (the partition on the right). On the bottom: The boxes in this partition are arranged into $3$-cylinders.      
  }

\begin{ex}
\label{p}
{\rm In Figure 8 the five boxes on the left are a realization of the polybox code $V=\{aaa,a'a'a',baa',a'ba, aa'b\}$, and the box in the middle is a realization of the word $w=bbb$. Since $w\sqsubseteq V$, we have $\breve{w}\subset \bigcup E(V)$. Thus, the $3$-box $\breve{w}$ is divided into pairwise dichotomous boxes $\breve{w}\cap \breve{v}$ for $v\in V$, and the set  $\bigcup (\{\breve{w}\cap \breve{v}: v\in Q\} \cup \{\breve{w}\cap \breve{v}: v\in P\})$, where  $P=\{v\in V^{3,a}: \breve{w}\cap \breve{v}\neq\emptyset\}$ and $Q=\{v\in V^{3,a'}: \breve{w}\cap \breve{v}\neq\emptyset\}$, is a $3$-cylinder in the box $\breve{w}$. Therefore, $\bigcup \{(\breve{w}\cap \breve{v})_3: v\in Q\}=\bigcup \{(\breve{w}\cap \breve{v})_3: v\in P\}$. Thus, the polybox $\bigcup \{(\breve{w}\cap \breve{v})_{3^c}: v\in Q\}$ is divided twice into pairwise dichotomous boxes without twin pairs and  $|Q|=|P|=2$. In Lemma \ref{l1} we will show that these two information allow us predict the structure of $Q_{3^c}$ and $P_{3^c}$: $Q_{3^c}=\{ba,a'a'\}$ and $P_{3^c}=\{aa,a'b\}$.

}
\end{ex}

$\bullet$  {\it The  structure of $W$ from the distribution of words in $V$}.  Below, in (\textbf{P}), (\textbf{V}), (\textbf{C}) and (\textbf{Co}) we show how to use an information on a  distribution of words in $V$ of the form (\ref{rep1}) to say something about distribution of words in $W$.

\vspace{0mm}
{\center
\includegraphics[width=11cm]{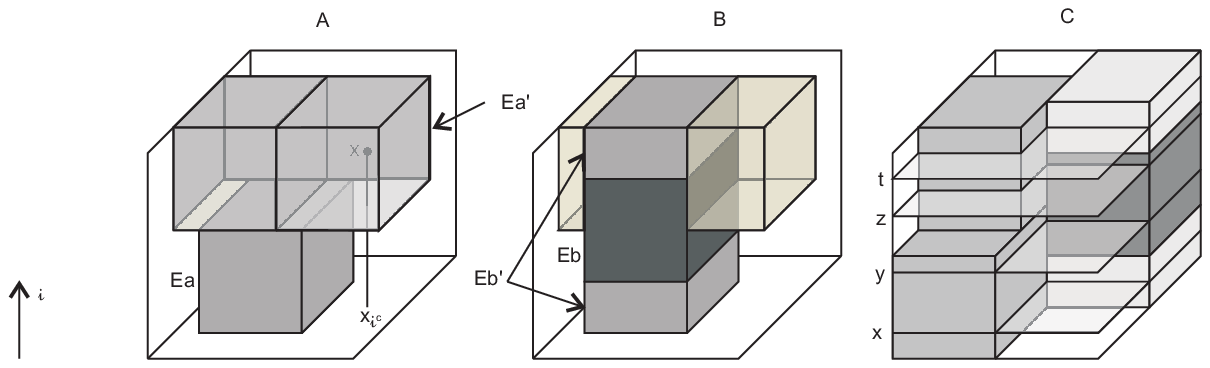}\\
}

\medskip
\noindent{\footnotesize 
Fig. 9. A scheme of realizations $E(V)$ (A), $E(W)$ (B) and $E(U)$ (C), where $V=V^{i,a}\cup V^{i,a'}$, $W=W^{i,a'}\cup W^{i,b}\cup W^{i,b'}$ and  $U=U^{i,a}\cup U^{i,a'}\cup U^{i,b}\cup U^{i,b'}$. We have $\bigcup E(V)=\bigcup E(W)$.    
  }

Let $V,W,U\subset S^d$ be polybox codes and assume that $V$ and $W$ are equivalent. Recall that $V^{i,l}_{i^c}=(V^{i,l})_{i^c}$ for every $i\in [d]$ and $l\in S$.

\medskip
\noindent
(\textbf{P}): {\bf Projections}. Suppose that there is $x\in \bigcup E(V^{i,l'})$ such that $x_{i^c}\not\in \bigcup E(V^{i,l}_{i^c})$ (see Figure 9A, where $l=a$). Since words in $V$ and $W$ are dichotomous, $\bigcup E(V)=\bigcup E(W)$ and (\ref{dkostki}), the point $x$ can be covered only by a box $\breve{w}\in E(W)$ such that $w_i=l'$. Thus, $W^{i,l'}\neq\emptyset$. In particular, if $W^{i,l}=\emptyset$ and $W^{i,l'}\neq\emptyset$, then $W^{i,l'}\sqsubseteq V^{i,l'}$ (Figure 9A).

\smallskip
\noindent
(\textbf{S}): {\bf Slices}. By (\ref{dkostki}), for every $r\in El\cap Es$, $l\not\in \{s,s'\}$ the set $\pi^i_r=ES\times \cdots \times ES\times \{r\}\times ES\times \cdots \times ES$, where $\{r\}$ stands at the $i$-th position, slices the sets $\bigcup E(U^{i,l})$ and $\bigcup E(U^{i,s})$ simultaneously (Figure 9C, where $r\in \{x,y,z,t\}$).
  
\smallskip
\noindent
(\textbf{V}): {\bf Volumes}. Let $|V^{i,l}|=n$ and $|V^{i,l'}|=m$, and let $n< m$. Since all boxes $\breve{u}$, $u\in S^d$, are of the same size and $n<m$, by ({\bf P}),
%there are words in $w\in W$ such that $w\in W^{i,a'}$ and 
$|W^{i,l'}|\geq m-n$.

\smallskip
\noindent
(\textbf{C}): {\bf Cylinders}. Suppose that $V^{i,l}\cup V^{i,l'}=\emptyset$ and $W^{i,l}\cup W^{i,l'}\neq\emptyset$ for some $l\in S$. Then $\bigcup E(W^{i,l}_{i^c})=\bigcup E(W^{i,l'}_{i^c})$, and hence the set $\bigcup E(W^{i,l}\cup W^{i,l'})$ in an $i$-cylinder in $(ES)^d$ (compare Figure 9A and 9B, where $l=b$). Indeed, if the set $\bigcup E(W^{i,l}\cup W^{i,l'})$ is not an $i$-cylinder, then $\bigcup E(W^{i,l}_{i^c})\neq\bigcup E(W^{i,l'}_{i^c})$.
%$$
%(\bigcup E((W^{i,b}_{i^c})\setminus \bigcup E((W^{i,b'}_{i^c}))\cup (\bigcup E((W^{i,b'}_{i^c})\setminus \bigcup E((W^{i,b}_{i^c}))\neq\emptyset.
%$$
By (\textbf{P}), $V^{i,l}\cup V^{i,l'}\neq\emptyset$, which is not true. Observe that, by (\ref{2d}),  the codes $W^{i,l}_{i^c}$ and $W^{i,l'}_{i^c}$ are equivalent.

%Note that $(W^{i,l}_{i^c}\sqsubseteq (V^{i,a}_{i^c}$ and $(W^{i,l}_{i^c}\sqsubseteq (V^{i,a'}_{i^c}$ for $l\in \{b,b'\}$.    

\smallskip
\noindent
(\textbf{Co}): {\bf Coverings}. Suppose that polybox codes $V^{i,l}_{i^c}$ and $W^{i,l}_{i^c}\cup W^{i,s_1}_{i^c}\cup \cdots \cup W^{i,s_k}$ are equivalent, where $s_n\not\in \{l,l',s_j,s_j'\}$ for every $n,j\in [k],n\neq j$. Then 
$$
W^{i,s_1}_{i^c}\cup \cdots \cup W^{i,s_k}\sqsubseteq V^{i,l'}_{i^c}\quad   {\rm and}\quad  V^{i,l}_{i^c}\sqsubseteq  W^{i,l}_{i^c}\cup W^{i,s_1'}_{i^c}\cup \cdots \cup W^{i,s_k'}_{i^c}.
$$
(In Figure 9A and 9B we have $V^{i,a}_{i^c}=W^{i,b}_{i^c}$, where $l=a,s_j=b$ for $j\in [k]$ and $W^{i,a}=\emptyset$). Indeed, since boxes in $E(V)$ are pairwise dichotomous and $\breve{w}_{i^c}\subseteq \bigcup E(V^{i,l}_{i^c})$ for every $w\in W^{i,s_1}_{i^c}\cup \cdots \cup W^{i,s_k}$,  each point $x\in \breve{w}\setminus \bigcup E(V^{i,l})$ has to be covered by the set $\bigcup E(V^{i,l'})$. Therefore,  $\breve{w}_{i^c}\subseteq \bigcup E(V^{i,l'}_{i^c})$ for every $w\in W^{i,s_1}_{i^c}\cup \cdots \cup W^{i,s_k}$, and consequently, by (\ref{2d}),  $w_{i^c}\sqsubseteq V^{i,l'}_{i^c}$ for $w\in W^{i,s_1}_{i^c}\cup \cdots \cup W^{i,s_k}$. Thus, $W^{i,s_1}_{i^c}\cup \cdots \cup W^{i,s_k}\sqsubseteq V^{i,l'}_{i^c}$. In the similar manner we show that  $V^{i,l}_{i^c}\sqsubseteq  W^{i,l}_{i^c}\cup W^{i,s_1'}_{i^c}\cup \cdots \cup W^{i,s_k'}_{i^c}$.

\bigskip
For fixed  $x\in ES$ and $i\in [d]$ let 
$$
\pi^i_x=ES\times \cdots \times ES\times \{x\}\times ES\times \cdots \times ES,
$$
where $\{x\}$ stands at the $i$th position.

\medskip
$\bullet$  {\it The structure of $V$ and $W$ from slices of the sets $\bigcup E(V)$ and $\bigcup E(W)$ by a set $\pi^i_x$}. Since $V\subset (*S)^d$ is a polybox code, the slice $\pi^i_x\cap \bigcup E(V)$ is a "flat" polybox in $(ES)^d$, that is boxes which are contained in this polybox have the factor $\{x\}$ at the $i$th position.

Therefore, we define a polybox $(\pi^i_x\cap \bigcup E(V))_{i^c}$ in the $(d-1)$-box $(ES)^{d-1}$:
$$
(\pi^i_x\cap \bigcup E(V))_{i^c}=\bigcup \{\breve{v}_{i^c}: v\in V\; {\rm and}\; \pi^i_x\cap \breve{v}\neq\emptyset\}.  
$$
The polybox $(\pi^i_x\cap \bigcup E(V))_{i^c}$ does not depend on a particular choice of a polybox code because if $W$ is an equivalent polybox code to $V$, then $\bigcup E(V)=\bigcup E(W)$, and hence $(\pi^i_x\cap \bigcup E(V))_{i^c}=(\pi^i_x\cap \bigcup E(W))_{i^c}$. 

%One of the main proof technique which is used in this paper derives from geometric tomography. 

We will slice a polybox  $\bigcup E(V)$ by the set $\pi^i_x$ for various $x\in ES$ (see (\textbf{S})). In particular, we will pay attention whether the polybox code $\{v_{i^c}:v\in V\;{\rm and}\; \pi^i_x\cap \breve{v}\neq\emptyset\}$ is rigid (Figure 10) because its rigidity will allow us to estimate the number of words in $V$ and $W$. %(To make the pictures more readable, in Figure 6 we use a realization $f(V)$ in $[0,1]^3$, as $E(V)$ is a discrete set. Clearly, adequate definitions for $f$ are analogous to that for $E$.)  

\vspace{-0mm}
{\center
\includegraphics[width=11cm]{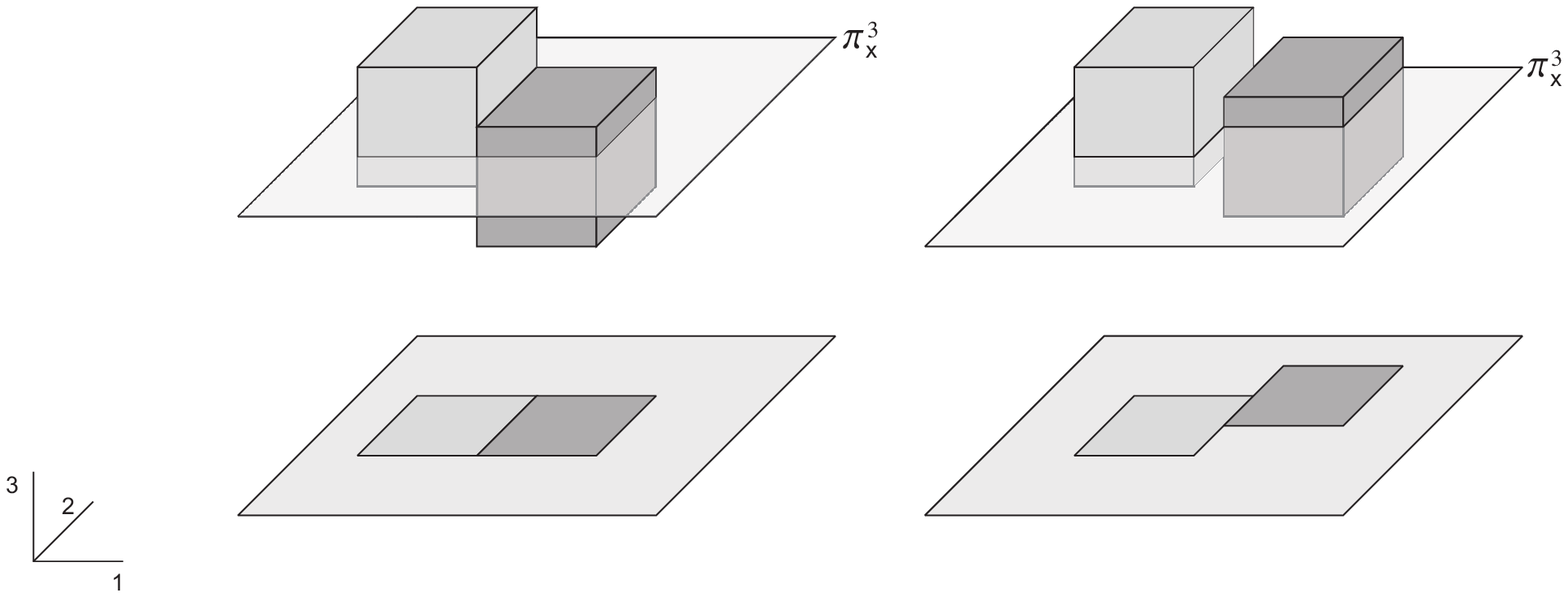}\\
}

\noindent{\footnotesize Fig. 10. The realizations $f(V)$ and $f(U)$ of two polybox codes $V=\{v,w\}$, where $v=aaa$ (the light box) and $w=a'ab$ (boxes on the left) and $U=\{v,u\}$, where $u=a'bb$ (boxes on the right) in the $3$-box $X=[0,1]^3$ ($X$ is not pictured).%, where $f_1(a)=[0,1/2),f_2(a)=[1/4,3/4),f_2(b)=[0,1/2],f_3(a)=[1/4,3/4)$ and $f_3(b)=[0,1/2)$. 
The polyboxes $\bigcup f(V)$ and $\bigcup f(U)$ are sliced by the set $\pi_x^3=[0,1]^2\times \{x\}$, where $x\in f_3(a)\cap f_3(b)$. Below each slice we have polyboxes in $X_{3^c}=[0,1]^2$: $(\pi^3_x\cap \bigcup f(V))_{3^c}=(f(v))_{3^c}\cup (f(w))_{3^c}$ (on the left) and $(\pi^3_x\cap \bigcup f(U))_{3^c}=(f(v))_{3^c}\cup (f(u))_{3^c}$ (on the right). The polybox code $\{v_{i^c},u_{i^c}\}$ is rigid, while $\{v_{i^c},w_{i^c}\}$  is not because it is a  twin pair. }

\medskip
Recall that the box number $|F|_0$ is the number of boxes in any proper suit for the polybox $F$. From (\ref{dkostki}) we deduce the following lemma which is useful in estimation of the number of  words in a polybox code by slices.
\begin{lemat}
\label{tech}
Let $V\subseteq S^d$ be a polybox code. Assume that there are letters $l_1,l_2\in S$, $l_1\not\in \{l_2,l_2'\}$, and $i\in [d]$ such that  $|(\pi^i_x\cap \bigcup E(V))_{i^c}|_0\geq m$ for every $x\in El_1\cap El_2$ and $|(\pi^i_y\cap \bigcup E(V))_{i^c}|_0\geq n$ for every $y\in El_1'\cap El_2'$ . Then $|V|\geq m+n$.
\end{lemat}

\proof Let $A\subseteq S$ be such that $V=\bigcup_{l\in A} V^{i,l}$ and $V^{i,l}\neq \emptyset$ for every $l\in A$. Let us divide the set $S$ into two disjoint sets $S^0$ and $S^1$ such that $S^1=\{l':l\in S^0\}$. Since $l_1\not\in \{l_2,l_2'\}$, we can assume that $l_1,l_2\in S_0$. We also divide the set $A$ into two disjoint sets $B$ and $C$ such that $B\subseteq S^0$ and $C\subseteq S^1$. 
%Let $B_1=B\cup C'$ and $C_1=C\cup B'$, where $B'=\{l':l\in B\}$ and $C'=\{l':l\in C\}$.
By (\ref{dkostki}), the sets $\bigcap_{l\in B}El\cap El_1\cap El_2$ and  $\bigcap_{l\in C}El\cap El_1'\cap El_2'$ are nonempty. Let $x\in \bigcap_{l\in B}El\cap El_1\cap El_2$ and  $y\in \bigcap_{l\in C}El\cap El_1'\cap El_2'$. Then $|(\pi^i_x\cap \bigcup E(V))_{i^c}|_0=\sum_{l\in B}|V^{i,l}|\geq m$ and $|(\pi^i_y\cap \bigcup E(V))_{i^c}|_0=\sum_{l\in C}|V^{i,l}|\geq n$. Since $|V|=\sum_{l\in B}|V^{i,l}|+\sum_{l\in C}|V^{i,l}|$, it follows that $|V|\geq n+m.$ 
\hfill{$\square$}

\medskip
$\bullet$ {\it The structure of $V$ from  the equality  $\sum_{v\in V} g(v,w)= 2^d$}. %Since $|\breve{v}|=1/2^n|(ES)^n|$, if $|V^{i,l}|\neq |V^{i,l'}|$, then we can estimate the cardinality of $W^{i,l}$ or $W^{i,l'}$ (see (V)). 
Let $V,W\subset S^d$ be equivalent and disjoint polybox codes. Then for every $w\in W$ we have $w\sqsubseteq V$ and $w\not\in V$. By (\ref{2d}),  $\sum_{v\in V} g(v,w)= 2^d$, where $g(v,w)\in \{0,1,2,\ldots ,2^{d-1}\}$ for every $v\in V$. Assume that $w=b\ldots b$ and let $\{v^1,\ldots ,v^k\}\subseteq V$ be such that $\breve{w}\cap \breve{v}^i\neq\emptyset$ for every $i\in [k]$ and $w\sqsubseteq \{v^1,\ldots ,v^k\}$. The solutions of the system of the equations $\sum_{i=1}^dx_i2^{d-i}=2^d, \sum_{i=1}^dx_i=k$, where $x_i$ are non-negative integers for $i\in [k]$, show the frequency of the letter $b$ in the words from the set $\{v^1,\ldots ,v^k\}$. We explain this on the following example. Recall first that, $g(v,w)= 2^{d-i}$ if and only if $v_j=b$ for every $j\in I\subset [d]$, $|I|=d-i$ and $v_j\not \in \{b,b'\}$ for $j\in [d]\setminus I$. In the example we assume that $d=3$, $w=bbb$ and $k=5$. The above system has two solutions: $x_1=0,x_2=3,x_3=2$ and $x_1=1,x_2=0,x_3=4$. It follows from the first solution that in the set $\{v^1,\ldots ,v^5\}$ there are exactly three words such that each of them contains exactly one letter $b$ and two words which have no letter $b$ or, by the second solution,  in the set $\{v^1,\ldots ,v^5\}$ there is exactly one word with two letters $b$ and the rest four words have no letter $b$. This observation is quite useful in the computations as it restricts the number of words which have to be considered during computations.
In the same way we can use the equality from Lemma \ref{2dd}.  

\subsection{Graph of siblings on a polybox code}

In Section 2.8 we described slices of a polybox $\bigcup E(V)$ by the sets $\pi^i_x$. Observe that if the set of boxes $\{v_{i^c}:v\in V\;{\rm and}\; \pi^i_x\cap \breve{v}\neq\emptyset\}$ contains a twin pair, say $v_{i^c}$ and $w_{i^c}$, and $V$ does not contain a twin pair, then $v_i\not\in \{w_i,w_i'\}$ (see Figure 10, the picture on the left). As we will see the number of such pairs $v,w$ in $V$ can help in estimation of the number of words in polybox codes $V$. Therefore, we now define a graph on a polybox code $V$.

Two words $v,u\in S^d$ such that $v_i\not\in \{u_i,u_i'\}$ for some $i\in [d]$, and  the pair $u_{i^c}$, $v_{i^c}$ is a twin pair are called $i$-{\it siblings} (in Figure 10 the boxes on the left are a realization of $i$-siblings for $i=3$). 

Let $V\subseteq S^d$ be a polybox code.  {\it A graph of siblings on V} is a graph $G=(V,\ka E)$ 
in which two vertices $u,v\in V$ are adjacent if they are $i$-siblings for some $i\in [d]$. We colour each edge in $\ka E$ with the colours from the set $[d]$: An edge $e\in \ka E$ has a colour $i\in [d]$ if its endpoints are $i$-siblings. The graph $G$ is simple and, if $V$ does not contain a twin pair, $d(v)\leq d$ for every $v\in V$, where $d(v)$ denotes the number of neighbors of $v$. (To show that  $d(v)\leq d$, suppose on the contrary that it is not true. Then there are two vertices $u$ and $w$ which are adjacent to $v$ such that $v_j'=u_j$ and $v_j'=w_j$. 
Since $v,u$ are $k$-siblings, and $v,w$ are $n$-siblings for some $k,n\in [d]\setminus \{j\}$, we have $v_k\not\in \{u_k,u_k'\}$, $v_n\not \in \{w_n,w_n'\}$ and $v_{\{k,j\}^c}=u_{\{k,j\}^c}$, $v_{\{n,j\}^c}=w_{\{n,j\}^c}$.
The vertices $u$ and $w$ are dichotomous, and therefore it must be $u_k'=w_n$ and $k=n$, which means that $u$, $w$ are a twin pair, a contradiction.)
Similarly it is easy to see that the graph $G$ does not contain triangles.

\begin{lemat}
\label{sybb}  
Let $G=(V,\ka E)$ be a graph of siblings on a polybox code $V\subset S^d$, $u$ and $v$ be adjacent vertices,  and let $d(u)=n$ and $d(v)=m$. If $n+m=2d$, then there are $i\in [d]$ and $l\in S$ such that   
%We shall show that if $n=d$ and $m=d$, then there are $j\in [d]$ and a letter $l\in S$ such that 
\begin{equation}
\label{nmm3}
|V^{i,l}\cup V^{i,l'}|\geq 2d-2, 
\end{equation}
and if $n+m\leq 2d-1$, then
\begin{equation}
\label{nm3}
|V^{i,l}\cup V^{i,l'}|\geq n+m-1 
\end{equation}
for some $i\in [d]$ and some $l\in S$. 
\end{lemat}
\proof By $N(u)$ and $N(v)$ we denote the set of all neighbors of $u$ and $v$, respectively. Assume without loss of generality that  $u=ll\ldots l$ and $v=sl\ldots l$, where $s\not\in \{l,l'\}$, and $u_i=l$, $v_i=l'$ for some $i\in \{2,\ldots, d\}$. Arrange the words from the set $N(u)\cup N(v)$ into a matrix $A$ such that words  from $N(u)\cup N(v)$ form  the rows of $A$.
%where $u$ and $v$ are its first two rows. 
By the definition of $i$-siblings, for every $w\in  (N(u)\cup N(v))\setminus \{u,v\}$ at most one letter $w_i$, where $i\in \{2,\ldots, d\}$, can be different from $l$ and $l'$. Therefore, if for every  $i\in \{2,\ldots, d\}$ there are at least three letters in the $i$th column of $A$ which are different from $l$ and $l'$, then there are at least $3(d-1)+2$ words in  the set $N(u)\cup N(v)$, which is impossible. Thus, there is $i\in  \{2,\ldots, d\}$ such that there are at least $2d-2$ words in $N(u)\cup N(v)\subset V$ with the letter $l$ or $l'$ at the $i$th position, and then $|V^{i,l}\cup V^{i,l'}|\geq 2d-2$. The proof of (\ref{nm3}) is very similar.
\hfill{$\square$} 

\medskip     
By $d(G)$ we denote the average degree of a graph $G$, and $N(S)$ denotes the set of all neighbors of vertices $v\in S$. 
In the sequel we will need the following lemma, which is probably known.  

\begin{lemat} 
\label{graph}
Let $G=(V,\ka E)$ be a simple graph, and let 
$$
m=\max\{d(v)+d(u):v,u\in V\;{\rm and}\; v,u\;{\rm are}\;{\rm adjacent}\}.
$$
Then  $d(G)\leq m/2$.
\end{lemat}

\proof Let $V_1\subset V$ be the set of all vertices $v$ such that $d(v)>m/2$, and let $\ka E_1\subset \ka E$ be the set of all edges which are incident with vertices from $V_1$. Since there is no edge with endpoints in the set $V_1$, the graph $G_1=(V,\ka E_1)$ is a bipartite with the bipartition $\{V_1,V\setminus V_1\}$. We will show that the graph $G_1$ contains a matching of $V_1$. To do this, let $S\subset V_1$. The number of edges in $\ka E_1$ which are incident with vertices from $S$ is greater than $|S|m/2$. On the other hand the  number of edges in $\ka E_1$, which are incident with vertices from $N(S)\subset V\setminus V_1$, is at most $|N(S)|m/2$. Each edge from $\ka E_1$ is incident with $S$ if and only if it is incident with $N(S)$. Therefore, $|N(S)|m/2>|S|m/2$, and thus $|N(S)|>|S|$. By
the marriage theorem, there is a matching of the set $V_1$. Let $V_2\subset V\setminus V_1$ be the set of endpoints of edges from the matching of $V_1$. Then $|V_1|=|V_2|$ and consequently
$$
d(G)=\frac{\sum_{v\in V_1}d(v)+ \sum_{v\in V_2}d(v)+\sum_{v\in V\setminus(V_1\cup V_2)}d(v)}{|V|}
$$    
$$
\leq \frac{|V_1|m+(|V|-2|V_1|)\frac{m}{2}}{|V|}=\frac{m}{2}.
$$    
\hfill{$\square$}

\medskip
Notice that if $u,v$ are $i$-siblings in a polybox code $V$ such that $u_i=l$ and  $v_i=s$, then for every $x\in El\cap Es$ the set $\{w_{i^c}\colon w\in V\;{\rm and}\; \pi^i_x\cap \breve{w}\neq\emptyset\}$ contains the twin pair $u_{i^c},v_{i^c}$ (see Figure 10).

\section{Small polybox codes without twin pairs.}

%To prove this theorem we need several auxiliary results. First, we will describe one polybox code and two partition codes without twin pairs which contain a few words. Next, we will prove the special cases in which  Theorem \ref{12} holds. 
%Finally, we will describe polybox and its two suits which, together with previously obtained results, allow us to prove Theorem \ref{12}.   
In graph theory it is very important to know the structure of all small graphs, that is, the graphs  with a small number of vertices. Similarly in the case of polybox codes it is very useful to know the structure of codes with a few words. In this section we describe first the structures of two equivalent polybox codes without twin pairs having two words each, and next we give the structures of two partition codes without twin pairs with five and six words.

\medskip
{\center
\includegraphics[width=6.6cm]{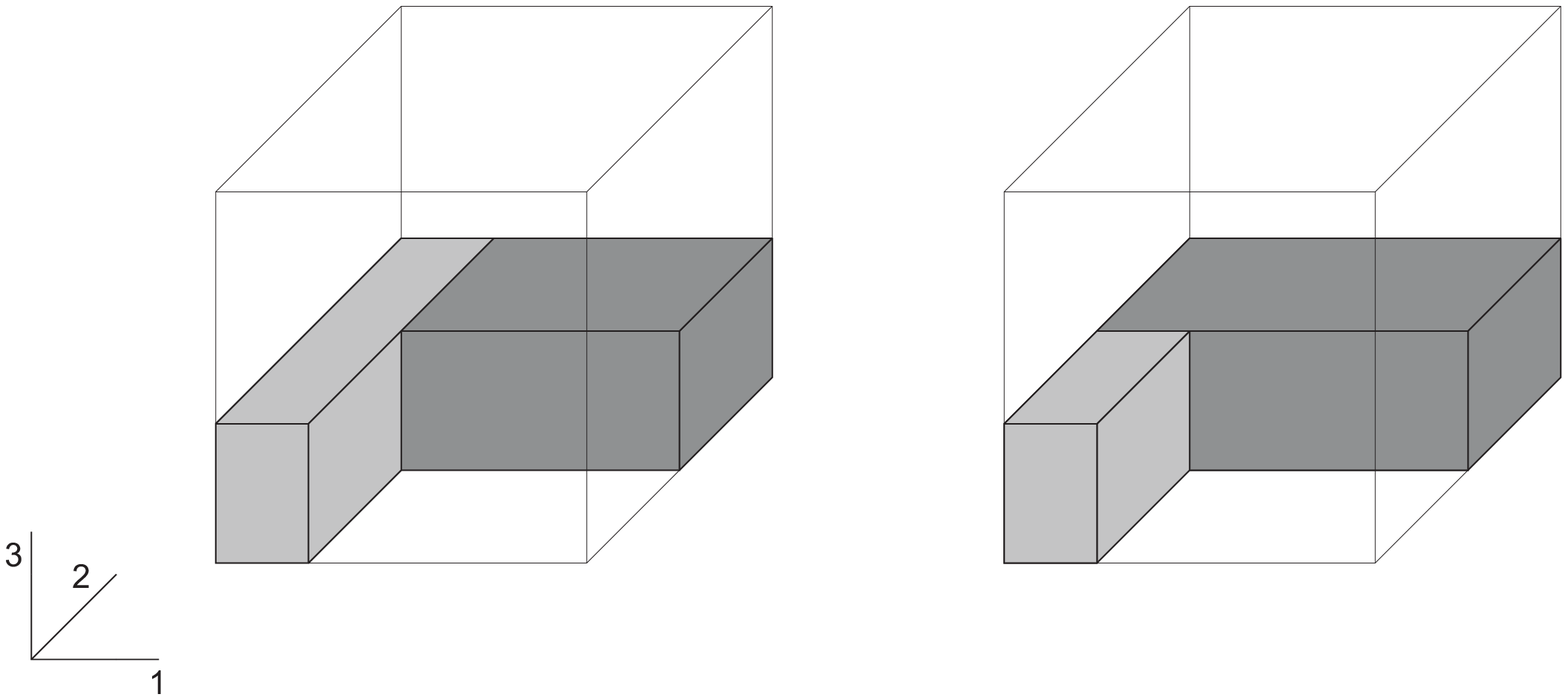}\\
}

\noindent{\footnotesize Fig. 11. Two realizations $f(U)$ and $f(P)$ in $X=[0,1]^3$ of polybox codes: $U=\{l_1*l_3,\; l_1'l_2l_3\}$, $l_1,l_2,l_3\neq *$, (on the left) and $P=\{l_1l_2'l_3, \; *l_2l_3\}$ (on the right), where $i_1=1,i_2=2$. Moreover, $f(l_1*l_3)=[0,1/4)\times [0,1]\times [0,1/3)$, $f(l_1'l_2l_3)=[3/4,1]\times [1/2,1]\times [0,1/3)$, $f(l_1l_2'l_3)=[0,1/4)\times [0,1/2)\times [0,1/3)$ and $f(*l_2l_3)=[0,1]\times [1/2,1]\times [0,1/3)$.   
  } 

\medskip 
Figure 11 presents realizations of polybox codes which are described in the following lemma. Recall that, if $v\in (*S)^d$ and $A\subset [d]$, then the word $v_A\in (*S)^{|A|}$ arises from $v$ by skipping  all the letters $v_i$ in the word $v$ which stand at the positions $i\in A^c=[d]\setminus A$. Moreover, if $V\subset (*S)^d$, then $V_A=\{v_A:v\in V\}$. 
  
\begin{lemat}
\label{l1}

Let $U,P\subset (*S)^d$ be disjoint sets which are equivalent polybox codes without twin pairs, and let $U=\{v,u\}$ and $P=\{w,q\}$. Then there are a set $A=\{i_1<i_2\}\subset [d]$ and letters $l_1,l_2\in S$ such that 
$$
U_A=\{*l_2, l_1l_2'\},\quad  P_A=\{l_1*, l_1'l_2\},
$$
and $v_i=u_i=w_i=q_i$ for every $i\in A^c$.

In particular, if $s,v,u,w,q\in S^d$ are such that $\breve{s}\cap \breve{p}\neq\emptyset$ for $p\in \{v,u,w,q\}$, the suits $\{\breve{s}\cap \breve{v},\breve{s}\cap \breve{u}\}$ and $\{\breve{s}\cap \breve{w}, \breve{s}\cap \breve{q}\}$ are disjoint, do not contain a twin pair and  $\breve{s}\cap \breve{v}\cup \breve{s}\cap \breve{u}=\breve{s}\cap \breve{w}\cup \breve{s}\cap \breve{q}$, then the codes $V=\{v,u\}$ and $W=\{w,q\}$ are of the forms
%then there are $i_1,i_2\in [d],i_1<i_2$ and $l_1,l_2\in S$, $l_j\not\in \{s_{i_j},s_{i_j}'\}$ for $j=1,2$, such that 
$$
V_A=\{s_{i_1}l_2, l_1l_2'\},\quad W_A=\{l_1s_{i_2}, l_1'l_2\}, 
$$
where $A=\{i_1<i_2\}\subseteq [d]$, $l_j\not\in \{s_{i_j},s_{i_j}'\}$ for $j=1,2$ and $v_i=u_i=w_i=q_i$,  $q_i\in S\setminus \{ s_i'\}$ for every $i\in A^c$.

\end{lemat}

\proof %Let $V=\{v,u\}$ and $W=\{w,q\}$.
Assume without loss of generality that $\breve{w}\cap \breve{v}\neq\emptyset$ and $\breve{q}\cap \breve{v}\neq\emptyset$. Thus, there is $k\in [d]$ such that $w_k=q_k'$, $q_k\neq *$ and $\breve{v}_{k^c}\subseteq \breve{w}_{k^c}\cap \breve{q}_{k^c}$. Since $w$ and $q$ are not a twin pair, we have $\breve{w}_{k^c}\setminus \breve{q}_{k^c}\cup \breve{q}_{k^c}\setminus \breve{w}_{k^c} \neq\emptyset$. Assume that $\breve{w}_{k^c}\setminus \breve{q}_{k^c}  \neq\emptyset$. Observe that a point $x\in \breve{w}$ such that $x_{k^c}\in \breve{w}_{k^c}\setminus \breve{q}_{k^c}$ cannot be contained in $\breve{v}$. Therefore, $x\in \breve{u}$. This means that $\breve{u}\cap \breve{q}=\emptyset$, for otherwise $\breve{u}\cap ((ES)^d\setminus (\breve{w}\cup \breve{q}))\neq\emptyset$, which is impossible. Then $\breve{v}=\breve{w}\cap \breve{v}\cup \breve{q}$ and $\breve{w}=\breve{u}\cup \breve{w}\cap \breve{v}$. Consequently,
%hence the boxes  $\breve{w}\cap \breve{v}$ and $\breve{q}$ are a twin pair, and similarly $\breve{u}$ and $\breve{w}\cap \breve{v}$ are a twin pair. Then 
$\breve{v}_{k^c}=\breve{q}_{k^c}$ and $\breve{u}_{j^c}=\breve{w}_{j^c}$, where $j\neq k$ is such that $u_j=v'_j$, $v_j\neq*$. Moreover, $\breve{v}_{\{k,j\}^c}=\breve{q}_{\{k,j\}^c}=\breve{u}_{\{k,j\}^c}=\breve{w}_{\{k,j\}^c}$, and then  $v_i=u_i=w_i=q_i$,  $q_i\in *S$ for every $i\in [d]\setminus \{i_1,i_2\}$, where $i_1=k$ and $i_2=j$.

Since $Eq_k\subsetneq Ev_k$, it follows that, by (\ref{dkostki}),  $Ev_k=ES$ and similarly $Ew_j=ES$, which means that $v_k=*$ and $w_j=*$. Let $v_j=l_2$ and $w_k=l_1$. Then $w_A=l_1*$ and $v_A=*l_2$. Since $\breve{v}_{k^c}=\breve{q}_{k^c}$, $w_k=q_k'$ and $\breve{u}_{j^c}=\breve{w}_{j^c}$, $u_j=v_j'$, we have $q_A=l_1'l_2$ and $u_A=l_1l_2'$. 

To prove the second part of the lemma let $U$ and $P$ be the polybox codes such that the suits $\ka F=\{\breve{s}\cap \breve{v},\breve{s}\cap \breve{u}\}$ and $\ka G=\{\breve{s}\cap \breve{w}, \breve{s}\cap \breve{q}\}$ are their exact realization in the $d$-box $X=\breve{s}$, respectively, where $U$ and $P$ are obtained in the manner described at the end of Section 2.6: $U=\{\bar{K}_1\ldots \bar{K}_d\colon K\in \ka F\}$,  $P=\{\bar{G}_1\ldots \bar{G}_d\colon G\in \ka G\}$ (recall that $\bar{K}_i=(K_i,i)$ if $K_i\neq X_i$ and $\bar{K}_i=*$ if $K_i=X_i$). 

Since $\bigcup \ka F=\bigcup \ka G$, by Lemma \ref{kodd}, the codes $U$ and $P$ are equivalent. As  $|U|=|P|=2$, by the first part of the lemma, we obtain  
$$
U_A=\{*(K_{i_2},i_2), (K_{i_1},i_1)(K_{i_2},i_2)'\}\quad {\rm and}\quad P_A=\{(K_{i_1},i_1)*, (K_{i_1},i_1)'(K_{i_2},i_2)\},
$$
where $A=\{i_1,i_2\}\subseteq [d]$, $K_{i_1}\neq X_{i_1}$, $K_{i_2}\neq X_{i_2}$. Moreover, $U_{A^c}=P_{A^c}=\{p_{A^c}\}$, where $p_i\in \{(Es\cap El,i):l\in S\setminus \{s,s'\}\}\cup \{*\}$ for $i\in A^c$.

Since $K_{i_1}=Es_{i_1}\cap El_1$, $K_{i_2}=Es_{i_2}\cap El_2$, where $l_j\not\in \{s_{i_j},s_{i_j}'\}$ for $j=1,2$ and $X_{i_1}=s_{i_1}, X_{i_2}=s_{i_2}$, we have $V_A=\{s_{i_1}l_2,l_1l_2'\}$, $W_A=\{l_1s_{i_2},l_1'l_2\}$ and $v_i=u_i=w_i=q_i$ for every $i\in A^c$.
%Since $\ka F$ and $\ka G$ do not contain a twin pair, the codes $V$ and $W$ do not contain such pairs and, by Lemma \ref{kodd}, they are equivalent.  Let  $\bar{V}=\{\bar{v}:v\in V\}$ be defined as in  Corollary \ref{uV}: For $v\in V$ if  $v_i\neq s_i$, then $\bar{v}_i=v_i$, and if  $v_i=s_i$, then $\bar{v}_i=*$. In the same way we define $\bar{W}$. Clearly, the polybox codes $\bar{V}, \bar{W}$ are equivalent and they do not contain a twin pair. Since $|\bar{V}|=|\bar{W}|=2$, by the first part of the lemma, $\bar{V}_A=\{*l_2, l_1l_2'\}$ and $\bar{W}_A=\{l_1*, l_1'l_2\}$, where $A=\{i_1,i_2\}$, $l_1,l_2\neq *$ and $\bar{v}_i=\bar{u}_i=\bar{w}_i=\bar{q}_i$ for every $i\in A^c$. This, by the definitions of the codes $\bar{V}$ and $\bar{W}$, completes the proof.
\hfill{$\square$}

\medskip
We now describe the structure of partition codes with five (Figure 12) and six words which do not contain twin pairs. 

\begin{lemat}
\label{5}
Let $V\subset (*S)^d$ be a partition code without twin pairs. 

\smallskip
\noindent
$(a)$ If $|V|>1$, then $d\geq 3$ and $|V|\geq 5$. The equality $|V|=5$ holds if and only if there are a set $A=\{i_1<i_2<i_3\}\subseteq [d]$ and letters  $l_1,l_2,l_3\in S$ such that 
$$
V_{A}=\{l_1l_2l_3, \; l_1'l_2'l_3', \; *l_2l_3', \; l_1'*l_3, \; l_1l_2'*\}
$$
 and $v=*\ldots *$ for every $v\in V_{A^c}$.

\smallskip
\noindent
$(b)$ If  $|V|=6$, then $d\geq 4$ and there are $i\in [d]$ and $l\in S$ such that 
$$
V=V^{i,l}\cup V^{i,l'}
$$
and $|V^{i,l}|=1$ and $|V^{i,l'}|=5$. Moreover, $v_{i^c}=*\ldots *$, where $V^{i,l}=\{v\}$, and 
$$
V^{i,l'}_{A}=\{l_1l_2l_3, \; l_1'l_2'l_3', \; *l_2l_3', \; l_1'*l_3, \; l_1l_2'*\}, 
$$
where $A=\{i_1<i_2<i_3\}\subseteq [d]\setminus \{i\}$, $l_1,l_2,l_3\in S$ and $u=*\ldots *$ for every $u\in V_{\{i_1,i_2,i_3,i\}^c}$.
\end{lemat}
{\it Proof of $(a)$} By (\ref{0indd}) there is a simple partition code $C\subset (*S)^d$ and there are words $v,u\in V\cap C$ such that the number $h=|\{i\in [d]: u_i=v'_i,\; u_i\neq *\}|$ is odd. Since $V$ does not contain a twin pair,
we have $h\geq 3$.
Thus, $d\geq 3$. Let $h=3$ and $\{i\in [d]: u_i=v'_i\}=\{i_1<i_2<i_3\}$. For every $j\in \{1,2,3\}$ let
$$
x^j\in Eu_1\times \cdots \times Eu_{i_j-1}\times Eu'_{i_j}\times Eu_{i_j+1}\times \cdots \times Eu_d.
$$
The points $x^1, x^2, x^3$ are pairwise different. Let us observe that for every $k,m\in \{1,2,3\}, k\neq m$, if $x^k,x^m\in \breve{w}$ for some $w\in (*S)^d$, then, as $\breve{w},\breve{u}$ are boxes, $\breve{u}\cap \breve{w}\neq\emptyset$ and consequently $w\notin V$.
Moreover, $x^1,x^2,x^3\not\in \breve{v}$ and $x^1,x^2,x^3\not\in \breve{u}$. Therefore, $|V|\geq 5$. In the same manner we  show that, if $h\geq 5$, then at least five words is needed to complete the set $\{u,v\}$ to a partition code. 

Let $|V|=5$, and let $v_{A}=l_1l_2l_3$ and $u_{A}=l'_1l'_2l'_3$. 
By Lemma \ref{2dd}, $\sum_{w\in V}\bar{g}(w,s_\ast)=2^{d}.$ Suppose that  $\bar{g}(w,s_\ast) =2^{d-1}$ for some $w\in V\setminus\{v,u\}$. (Recall that, $\bar{g}(v,s_*)=2^{d-i}$, $i\in \{0,\ldots,d\}$, if and only if the word $v$ contains $d-i$ stars.) Then there is exactly one $i\in [d]$ such that $w_i\neq *$ and $w_{i^c}=*\ldots *$. Since $V$ is a partition code, it follows that $(V\setminus\{w\})_{i^c}\subset (*S)^{d-1}$ is a partition code. This code does not contain a twin pair and consists of four words, which is impossible. Therefore, by Lemma \ref{2dd}, $\bar{g}(u,s_\ast)=\bar{g}(v,s_\ast)=2^{d-3}$ and $\bar{g}(w,s_\ast) =2^{d-2}$ for every $w\in V\setminus\{v,u\}$.  Since  every two words in $V$ are dichotomous, we have $(V\setminus\{v,u\})_{A}= \{*l_2l_3', \; l_1'*l_3, \; l_1l_2'*\}$ %or  $(V\setminus\{v,u\})_{A}=\{*l_2'l_3, \; l_1*l_3', \; l_1'l_2*\}$ 
and $w_{A^c}=*\ldots *$ for every $w\in V$.

\medskip
\noindent
{\it Proof of $(b)$}. Let $C$, $v,u$ and $h$ be such as in the proof of the part $(a)$.  
%By (\ref{0indd}), there is a simple partition code $C$ and two words $v=v_1\cdots v_d,u=u_1\cdots u_d\in C\cap V$ such that the number $h=|\{i\in [d]: u_i=v'_i\}|$ is odd. Since $V$ does not contain twin pairs, $n\geq 3$. 
If $h\geq 5$, then, in the similar way as in the proof of $(a)$, we show that at least five words is needed to complete the set $\{u,v\}$ to a partition code. Thus, $h=3$ and therefore $\bar{g}(u,s_*)=\bar{g}(v,s_*)=2^{d-3}$.  Let $\{i\in [d]: u_i=v'_i,\; u_i\neq *\}=\{i_1,i_2,i_3\}$, and $V=\{v^1,\ldots, v^6\}$, where  $v^3=v,v^4=u$. Assume on the contrary that for every $i\in [d]$ and $l\in S$ we have $V\neq V^{i,l}\cup V^{i,l'}$. Then, by (a), $\bar{g}(v^i,s_*)\leq 2^{d-2}$ for every $v^i\in V$. We will use again Lemma \ref{2dd}.  We consider the system of the equations   
\begin{equation}
\label{eqqu}
\sum_{i\geq 2}x_i2^{d-i}=2^d,\;\;\;\; \sum_{i\geq 2}x_i=6.
\end{equation}
where  $x_i$ are non-negative integers. Since $x_3\geq 2$, because $\bar{g}(u,s_*)=\bar{g}(v,s_*)=2^{d-3}$, we consider only one solution of this system: $x_2=2$, $x_3=4$ (the second solution is $x_4=2,x_3=1,x_2=3$). Therefore,
\begin{equation}
\label{ee}
\bar{g}(v^1,s_*)=\bar{g}(v^2,s_*)=2^{d-2}\;\;{\rm and}\;\; \bar{g}(v^i,s_*)=2^{d-3}
\end{equation}
for $i=3,\ldots ,6$. Let us consider the realization $E(V)$. Clearly, we can assume that for every $i\in [d]$ there is $v\in V$ such that $v_i\neq *$, for otherwise we consider the code $V_{i^c}\subset (*S)^{d-1}$. 

Observe now that, every non-empty $i$-cylinder $C_i=\bigcup \{\breve{v}^j: v^j_i\in \{l,l'\}\}$, $l\in S$, $i\in [d]$, has to contain at least four boxes. Indeed, if, on the contrary,  $C_i=\breve{v}^j\cup \breve{v}^k$, then $v^j$ and $v^k$ form a twin pair, and if $C_i=\breve{v}^j\cup \breve{v}^k\cup \breve{v}^n$, $v^j_i=v^k_i, v^n_i=(v_i^k)'$, then $v^j$ and $v^k$ are twins, as $\breve{v}^j_{i^c}\cup \breve{v}^k_{i^c}=\breve{v}^n_{i^c}$, and thus 
$\breve{v}^j\cup \breve{v}^k$ is a box, which is equivalent to say that $v^j$ and $v^k$ form a twin pair. A contradiction.
     
It follows from the above and (\ref{ee}) that we can always choose $i\in [d]$ such that $v^1_i=v^2_i=*$ or $v^k_i=v^n_i=*$ for some $k\in \{1,2\}$ and $n\in \{3,...,6\}$. In the first case we have  $C_i=\breve{v}^3\cup \breve{v}^4\cup \breve{v}^5\cup \breve{v}^6$. 
Since $C_i$ is an $i$-cylinder and (\ref{ee}), we can assume that 
$$
\breve{v}^3_{i^c}\cup \breve{v}^4_{i^c}=\breve{v}^5_{i^c}\cup \breve{v}^6_{i^c},
$$
where $v^3_i=v^4_i,v^5_i=v^6_i$ and $v^3_i=(v^5_i)'$. As  $V$ does not contain a twin pair, by Lemma \ref{l1}, the polybox codes $\{v^3_{i^c},v^4_{i^c}\}$ and $\{v^5_{i^c},v^6_{i^c}\}$ are of the forms given in the first part of this lemma. Thus, $\bar{g}(v^3,s_*)\neq \bar{g}(v^4,s_*)$, which contradicts (\ref{ee}). 

In the second case the set $\bigcup\{\breve{s}:s\in V\setminus\{v^k,v^n\},\;s_i\neq *\}$ cannot be an $i$-cylinder. Indeed, we have $\bar{g}(v^i,s_*)=2^{d-2}$ for exactly one $v^i\in V\setminus\{v^k,v^n\}$ and $\bar{g}(w,s_*)=2^{d-3}$ for the rest three words $w\in V\setminus\{v^k,v^m,v^i\}$. Thus, $|\breve{v}^i|=(1/4)|ES|^d$ and $|\breve{w}|=(1/8)|ES|^d$ for $w\in V\setminus\{v^k,v^m,v^i\}$. Therefore, the boxes from the set  $\{\breve{s}:s\in V\setminus\{v^k,v^n\},\;s_i\neq *\}$ cannot be divided into two parts (one with all words having $s_i$ at the $i$th position and the second with words having the letter $s_i'$ at the $i$th position) with the same sizes.   

Thus, there is a word, say $v^1$, in the code $V$ such that $\bar{g}(v^1,s_*)=2^{d-1}$.  If $i\in [d]$ is such that $v^1_i=l$, where $l\neq *$, then $v^2_i=\cdots =v^6_i=l'$, and hence $V^{i,l}=\{v^1\}$ and  $V^{i,l'}=\{v^2,\ldots v^6\}$. Clearly, $v^1_{i^c}=*\ldots *$. Since $\{v^2,\ldots ,v^6\}$ does not contain a twin pair and $v^1_{i^c}\sqsubseteq  \{v^2_{i^c},\ldots ,v^6_{i^c}\}$, by $(a)$,  $d\geq 4$. The form of $V^{i,l'}$ is guaranteed by the part $(a)$. This completes the proof of the part $(b)$.
\hfill{$\square$}   

{\center
\includegraphics[width=7cm]{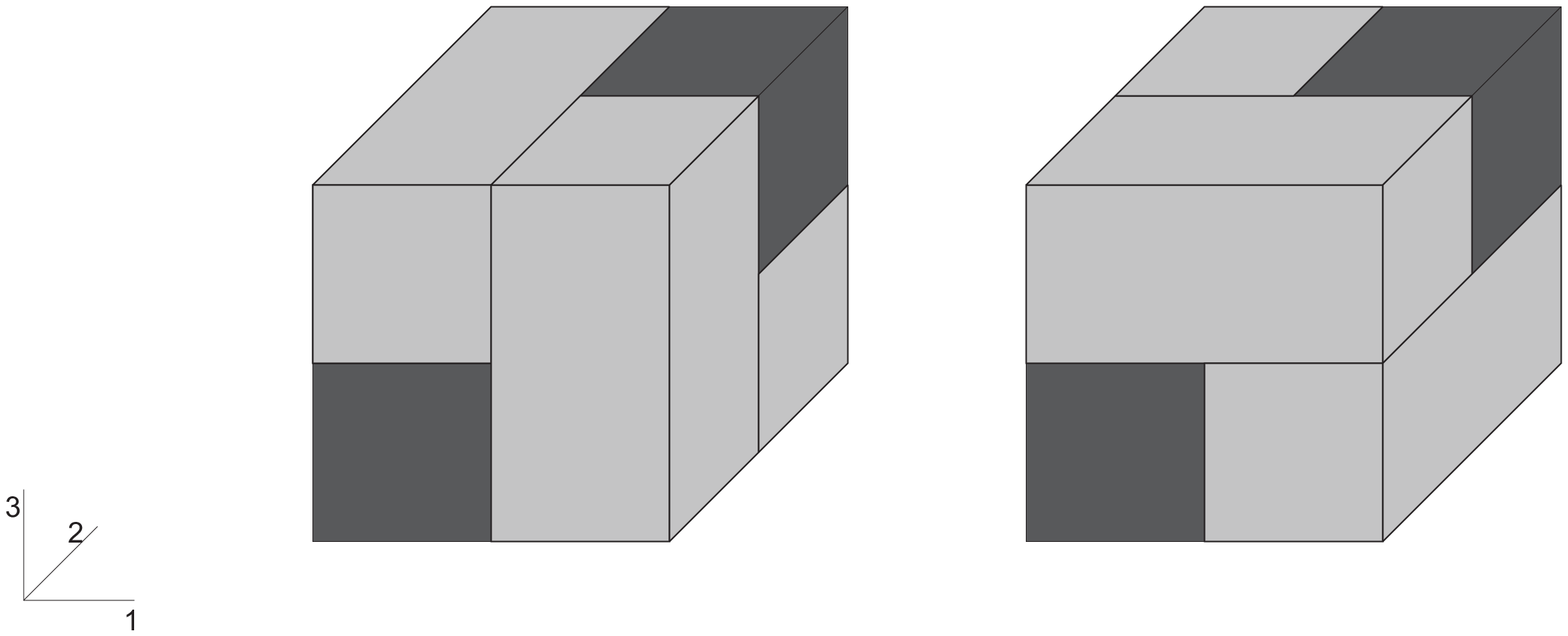}\\
}
\medskip
\noindent{\footnotesize Fig. 12. Let $*S=\{a,a',*\}$. On the right we see a realization $f(V)$ (in $X=[0,1]^3$) of the code $V= \{l_1l_2l_3, \; l_1'l_2'l_3', \; *l_2l'_3, \\ \; l'_1*l_3, \; l_1l'_2*\}$  for $l_1=l_2=l_3=a$ and on the left for $l_1=l_2=l_3=a'$, where $f_i(a)=[0,1/2)$ and $f_i(*)=[0,1]$ for $i=1,2,3$. The words $v$ and $w$, where $v=l_1l_2l_3$ and $w=l_1'l_2'l_3'$ are such as in (\ref{0indd}).%and $(V)_{A}=\{l_1l_2l_3, \; l_1'l_2'l_3', \; *l_2l_3', \; l_1'*l_3, \; l_1l_2'*\} $ (on the right). 
 
  }

\begin{wn}
\label{uww}
Let $V\subset S^d$ be a polybox code without twin pairs, $u\sqsubseteq V$ for some $u\in S^d\setminus V$ and $\breve{u}\cap \breve{v}\neq\emptyset$ for every $v\in V$. Then $|V|\geq 5$. If $|V|=5$, then there are a set  $A=\{i_1<i_2<i_3\}\subseteq [d]$ and letters $l_1,l_2,l_3\in S\setminus \{u_{i_j},u'_{i_j}\}$ for $j=1,2,3$ such that
%obtained from one of the partition codes described in the part $(a)$ of Lemma \ref{5}, by replacing each star standing at the $i$th position by a letter $u_i$ for every $i\in [d]$. Observe that such obtained code $V$, where  
$$
V_{A}=\{l_1l_2l_3, \; l_1'l_2'l_3', \; u_{i_1}l_2l_3', \; l_1'u_{i_2}l_3, \; l_1l_2'u_{i_3}\},
$$
%or
%$$
%(V)_{A}=\{l_1l_2l_3, \; l_1'l_2'l_3', \; u_{i_1}l'_2l_3, \; l_1u_{i_2}l'_3, \; l'_1l_2u_{i_3}\},
%$$
and $v_i=u_i$ for every $i\in A^c$ and $v\in V$. In particular, $V$ is rigid. 

If $|V|=6$ and $\breve{u}\cap \breve{v}\neq\emptyset$ for every $v\in V$, then there are $i\in [d]$ and a letter $l\in S$ such that $u_{i^c}=v_{i^c}$, where 
$V^{i,l}=\{v\}$ and
%a set  $A=\{i_1<i_2<i_3\}\subseteq [d]\setminus \{i\}$ and letters $l_1,l_2,l_3\in S$, $l_j\in S\setminus \{u_{i_j},u'_{i_j}\}$ for $j=1,2,3$ such that
$$
V^{i,l'}_{A}=\{l_1l_2l_3, \; l_1'l_2'l_3', \; u_{i_1}l_2l_3', \; l_1'u_{i_2}l_3, \; l_1l_2'u_{i_3}\},
$$
where $A=\{i_1<i_2<i_3\}\subseteq [d]\setminus \{i\}$, $l_1,l_2,l_3\in S\setminus \{u_{i_j},u'_{i_j}\}$ for $j=1,2,3$ and $v_j=u_j$ for every $j\in [d]\setminus \{i_1,i_2,i_3,i\}$ and $v\in V$. 
\end{wn}
\proof %Let  $\bar{V}=\{\bar{v}:v\in V\}$ be defined as in  Corollary \ref{uV}: for $v\in V$ if  $v_i\neq u_i$, then $\bar{v}_i=v_i$, and if  $v_i=u_i$, then $\bar{v}_i=*$.
We consider the code $\bar{V}=\{\bar{v}\colon v\in V\}$, where $\bar{v}_i=v_i$ if $v_i\neq u_i$ and $\bar{v}_i=*$ if $v_i=u_i$.
By  Corollary \ref{uV}, the set $\bar{V}$ is a partition code. Since $V$ does not contain a twin pair,  $\bar{V}$ does not contain a twin pair. By Lemma \ref{5} $(a)$, $|\bar{V}|\geq 5$, and thus $|V|\geq 5$, because $|\bar{V}|=|V|$. Similarly, by Lemma \ref{5}, if $|\bar{V}|=5$, the code $\bar{V}$  is of the form described in the part (a) of this lemma, and thus, by the definition of $\bar{V}$, the form of the polybox code $V$ has to be as given in the corollary. 

To show that $V$ is rigid, it is enough to observe that the only word which can be covered by $V$ is $u$.
Indeed, if $w\sqsubseteq V$ for some $w\in S^d$, then  $w_i=u_i$ for $i\in A^c$ (if not, by (\ref{dkostki}), there is $x\in \breve{w}$ with $x_i\in Ew_i\setminus Eu_i$ and then $x\not\in \breve{v}$ for every $v\in V$). 

If $w_{i_j}\in \{l_1,l'_1, l_2,l'_2,l_3,l'_3\}$ for some $j\in \{1,2,3\}$, then $\breve{w}\cap \breve{v}=\emptyset$ for some $v\in V$. Then $w\sqsubseteq V\setminus \{v\}$, which contradicts the first part of the lemma.
Thus $w_{i_j}\not\in \{l_1,l'_1, l_2,l'_2,l_3,l'_3\}$ for every $i\in \{1,2,3\}$. 

If $w_{i_j}\neq u_{i_j}$ for some $j\in \{1,2,3\}$, then $g(v,w)=2^{d-2}$ for at most two words $v\in V$ and $g(v,w)=2^{d-3}$ for the rest $v\in V$. Then $\sum_{v\in V}g(v,w)<2^d$, and by (\ref{2d}), $w\not\sqsubseteq V$. Thus, $w=u$, and consequently $V$ is rigid.    

The proof of the corollary in the case $|V|=6$ is the same as the proof of the case $|V|=5$, but instead of Lemma \ref{5} $(a)$ we use  Lemma \ref{5} $(b)$.
%, $d\geq 4$. Moreover, there is $i\in [d]$ and $\bar{v}\in \bar{V}$ such that $(\bar{v}_{i^c}=*\cdots *$. By the definition of $\bar{v}$, $(v_{i^c}=(u_{i^c}$.
\hfill{$\square$} 

\begin{wn}
\label{<3}
Let $V\subset S^d$ be a polybox code which does not contain a twin pair. If $d\leq 3$, then $V$ is rigid and $|V|\leq 5$.    
\end{wn}
\proof It is enough to prove the lemma for $d=3$. If $|V|\leq 5$, then the rigidity of $V$ is guaranteed by Corollary \ref{uww}. It is easy to show (\cite{LP2}) that for every partition code $U\subseteq S^3$ we have $U=U^{i,l}\cup U^{i,l'}$ for some $i\in [3]$ and some $l\in S$ or 
$$
U=\{l_1l_2l_3, l_1'l_2'l_3',  u_{1}l_2l_3', u'_{1}l_2l_3',  l_1'u_{2}l_3,  l_1'u'_{2}l_3,  l_1l_2'u_{3}, l_1l_2'u'_{3}\},
$$
where $u_i\not \in \{l_i,l'_i\}$ for $i\in [3]$. It was shown in \cite{DIP} that every polybox code $V\subseteq S^3$ with at least five words is extensible to a partition code, that is, there is a code $V'\subset S^3$ such that $V\cup V'$ is a partition code. Note that every three words in $U^{i,l}$ and $U^{i,l'}$ contain a twin pair because  $U^{i,l}_{i^c}$ and  $U^{i,l'}_{i^c}$ are partition codes in dimension two. Therefore, any polybox code $V\subset S^3$ with $|V|\geq 6$ contains a twin pair.
%In \cite{CS2} it was computed that the maximal clique in the $3$-dimensional Keller graph has  five vertices. This means that if $S=\{a,a',b,b'\}$ and $|V|\geq 6$, then there is a twin pair in $V$. The same is true for arbitrary $S$. 
\hfill{$\square$}

\section{Equivalent polybox codes without twin pairs}

In this section we prove Theorem \ref{12}. We first prove three lemmas on the properties of polybox codes
%: the first lemma give a condition which implies the existence of twin pairs in a polybox code, and the second lemma says that a polybox code with at most seven words is rigid. The third lemma in this section and the followed by it 
and next we give  the special cases in which Theorem \ref{12} holds.  
 
%We start from the following fact about twin pairs in a polybox code $V$.  

\begin{lemat}
\label{3}
Let $V\subseteq S^d$ be a polybox code, $w\in S^d$, $w\not\in V$ and $w\sqsubseteq V$. Suppose that there are $i\in [d]$ and $l\in S$ such that $l\not\in \{w_i,w_i'\}$, $ |\{v\in V^{i,l}\colon \breve{w}\cap \breve{v}\neq\emptyset\}|=1$ and $1\leq |\{v\in V^{i,l'}\colon \breve{w}\cap \breve{v}\neq\emptyset\}|\leq 4$. Then there is a twin pair in $V$. 
%In particular, if  $\breve{w}\cap \breve{v}\neq\emptyset$ for every $v\in V$ and $(|V|-3)\leq |V^{i,w_i}|<|V|$, then there is a twin pair in $V$. 
\end{lemat}
\proof  %Suppose that there is a word $w\in W$ such that $w_i\not \in \{l,l'\}$ and $\breve{w}\cap \bigcup E(V^{i,l})\neq\emptyset$. Then, by (P), we have also  $\breve{w}\cap \bigcup E(V^{i,l'})\neq\emptyset$. 
Since $w\sqsubseteq V$, the set $\breve{w}\cap \bigcup E(V^{i,l}\cup V^{i,l'})$ is an $i$-cylinder in the $d$-box $\breve{w}$ (compare Example \ref{p}). Therefore, 
$$
(\breve{w}\cap \breve{u})_{i^c}=\bigcup\{(\breve{w}\cap \breve{v})_{i^c}:v\in V^{i,l'}\},
$$
where $u$ is the sole word in $V^{i,l}$ such that $\breve{u}\cap\breve{w}\neq\emptyset$. 

Let $|\{\breve{w}\cap \breve{v}\neq\emptyset :v\in V^{i,l'}\}|=1$.  Then, $(\breve{w}\cap \breve{u})_{i^c}=(\breve{w}\cap \breve{v})_{i^c}$, where $v\in V^{i,l'}$. By Lemma \ref{=c}, $u_{i^c}=v_{i^c}$, which means, as $u_i=l$ and $v_i=l'$, that $u$ and $v$ form a twin pair.

Let now $(\breve{w}\cap \breve{u})_{i^c}=(\breve{w}\cap \breve{v}^1)_{i^c}\cup \cdots \cup (\breve{w}\cap \breve{v}^k)_{i^c}$, where $v^1,\ldots ,v^k$ are the words in $V^{i,l'}$ such that $\breve{w}\cap \breve{v}^j\ne\emptyset$ for $j\in [k]$ and $k\in \{2,3,4\}$. Since the family of boxes $\{(\breve{w}\cap \breve{v}^1)_{i^c}, \ldots ,(\breve{w}\cap \breve{v}^k)_{i^c}\}$ is a suit for the $(d-1)$-box $(\breve{w}\cap \breve{u})_{i^c}$ 
%(that is, boxes from this family are pairwise dichotomous) 
and  $1<k<5$, by Lemma \ref{5}, there are $n,m\in [k]$, $n\neq m$, such that the boxes $(\breve{w}\cap \breve{v}^n)_{i^c}$ and $(\breve{w}\cap \breve{v}^m)_{i^c}$ are a twin pair (indeed, if we write the polybox code of this suit in the manner described in Section 2.6, where $X=(\breve{w}\cap \breve{u})_{i^c}$, we obtain a partition code with less than five words). Then, by Lemma \ref{=c}, $v^n_{i^c}$ and $v^m_{i^c}$ are a twin pair. Since $v_i^n=v^m_i=l'$, the words $v^n$ and $v^m$ form a twin pair. 
\hfill{$\square$}

\begin{lemat} 
\label{komp}
Let $S=\{a,a',b,b'\}$, and let $V\subset S^d$  be a polybox code without twin pairs. If $|V|\leq 7$, then $V$ is rigid.%, and for $d\geq 5$ if $|V|\leq 8$, then $V$ is rigid.
\end{lemat}
\proof %The observation made at the beginning of the proof of Statement \ref{pusty} is also helpful in the present proof.  
Suppose on the contrary that there is a polybox code  $W$ which is equivalent to $V$ and $V\cap W=\emptyset$ (observe that the code $W$ can contain a twin pair). % for $d=4$, and $|V|>8$ for $d\geq 5$. 
We proceed by induction on $d$. By Corollary \ref{<3}, the lemma is true for $d \leq 3$. Let $d\geq 4$. We will show first that $V^{i,l}\neq\emptyset$ for every $i\in [d]$ and $l\in S$. 

Suppose that this is not true. We may assume that $V^{i,a}\neq\emptyset$ and $V^{i,b}=\emptyset$ for some $i\in [d]$. Let $x\in Ea\cap Eb$. Observe that, since 
$\pi^i_x\cap \bigcup E(V^{i,a}\cup V^{i,b})= \pi^i_x\cap \bigcup E(W^{i,a}\cup W^{i,b})$, we have 
$\bigcup \{\breve{v}_{i^c}:v\in V^{i,a}\cup V^{i,b}\}=\bigcup \{\breve{w}_{i^c}:w\in W^{i,a}\cup W^{i,b}\}$.
It follows from  (\ref{2d}) (compare the comments below (\ref{2d})) that the polybox codes $V^{i,a}_{i^c}\cup V_{i^c}^{i,b}$ and $W^{i,a}_{i^c}\cup W^{i,b}_{i^c}$ are equivalent, and since $V^{i,b}=\emptyset$, the codes  $V^{i,a}_{i^c}$ and $W^{i,a}_{i^c}\cup W^{i,b}_{i^c}$ are equivalent . By the inductive hypothesis the code $V^{i,a}_{i^c}\subset S^{d-1}$ is rigid. Thus,  $V^{i,a}_{i^c}=W^{i,a}_{i^c}\cup W^{i,b}_{i^c}$. Note that $W^{i,a}_{i^c}=\emptyset$, for otherwise $V^{i,a}\cap W^{i,a}\neq\emptyset$, which is impossible. Hence, $V^{i,a}_{i^c}=W^{i,b}_{i^c}$. Then, by ($\textbf{Co}$) in Section 2.3,  $V^{i,a}_{i^c}\sqsubseteq V^{i,a'}_{i^c}$, and thus $V^{i,a'}\neq\emptyset$. In the same way we show that $V^{i,a'}_{i^c}\sqsubseteq V^{i,a}_{i^c}$. Consequently the codes $V^{i,a'}_{i^c},V^{i,a}_{i^c}$ are equivalent. Since  $V^{i,a}_{i^c}$ is rigid, these codes are equal. Then the set $V^{i,a}\cup V^{i,a'}$ contains a twin pair, which is impossible. 

Hence, $V^{i,l}\neq\emptyset$ for every $i\in [d]$ and $l\in S$. 

\smallskip
We now show that for every $l,s\in S$, $l\not \in \{s,s'\}$, the code $V^{i,l}\cup V^{i,s}$ contains $i$-siblings.

To do this, let us suppose on the contrary that there are  $i\in [d]$ and two letters in $S$, say $a$ and $b$, such that there are no $i$-siblings in $V^{i,a}\cup V^{i,b}$, that is, the polybox code $V^{i,a}_{i^c}\cup V^{i,b}_{i^c}$ does not contain a twin pair. This means, by the inductive hypothesis,  that the polybox code $V^{i,a}_{i^c}\cup V^{i,b}_{i^c}$ is rigid. As it was shown, the polybox codes  $V^{i,a}_{i^c}\cup V^{i,b}$ and $W^{i,a}_{i^c}\cup W^{i,b}_{i^c}$ are equivalent and therefore, by the rigidity of  $V^{i,a}_{i^c}\cup V^{i,b}_{i^c}$, we have  $V^{i,a}_{i^c}\cup V^{i,b}_{i^c}=W_{i^c}^{i,a}\cup W_{i^c}^{i,b}$. Since $V\cap W=\emptyset$, it follows that $V^{i,a}_{i^c}=W_{i^c}^{i,b}$ and $V^{i,b}_{i^c}=W_{i^c}^{i,a}$. By  ($\textbf{Co}$) in Section 2.3,  $V^{i,a}_{i^c}\sqsubseteq V_{i^c}^{i,a'}$ and $V^{i,b}_{i^c}\sqsubseteq V_{i^c}^{i,b'}$. Note that $V^{i,a}_{i^c}\cap  V_{i^c}^{i,a'}=\emptyset$ and $V^{i,b}_{i^c}\cap V_{i^c}^{i,b'}=\emptyset$ because $V$ does not contain twin pairs. Therefore, by Lemma \ref{uww},  $|V|\geq |V^{i,a'}|+|V^{i,b'}|\geq 10$, a contradiction.

\smallskip
Thus we may assume that for every $i\in [d]$ and every two letters $l,s\in S$, $l\not\in \{s,s'\}$, there are $i$-siblings in the set $ V^{i,l}\cup V^{i,s}$. 
Then for every $i\in[d]$ there are at least $4$ edges with the colour $i$ in the graph of siblings $G=(V,\ka E)$ on $V$, and consequently there are at least $4d$ edges in the set $\ka E$.
To finish the proof we will show that $V$ has to contain more than $7$ vertices. 

Let $u^0,v^0\in V$ be such that 
$$
d(v^0)+d(u^0)=\max\{d(v)+d(u)\colon v,u\in V\;{\rm and}\; v,u\;{\rm are}\;{\rm adjacent}\}.     
$$
Since $|V|\leq 7$, we have $d(v^0)+d(u^0)\leq 7$, and then it follows from Lemma \ref{graph}  that
$$
d(G)\leq \frac{7}{2}.
$$ 
But $d(G)|V|=2|\ka E|$ and $2|\ka E|\geq 32$. Therefore, $|V|>7$, a contradiction. 
\hfill{$\square$}

\medskip
\begin{lemat}
\label{dis}
Let $S=\{a,a',b,b'\}$, and let $V,W\subset S^4$ be disjoint sets which are equivalent polybox codes without twin pairs.
% and assume that $V$ has the smallest possible cardinality in the sense that if $\tilde{V},\tilde{W}\subset S^d$ are disjoint and equivalent polybox codes and $|\tilde{V}|<|V|$, then there is a twin pair in $\tilde{V}$ or in $\tilde{W}$. 
Then $V^{i,l}\neq\emptyset$ and $W^{i,l}\neq\emptyset$ for every $l\in S$. 
If $d\geq 5$, $U,P\subset S^{d}$ are disjoint sets which are equivalent polybox codes without twin pairs and $U^{i,l}=\emptyset$ for some $i\in [d]$ and some $l\in S$, then $|U|\geq {\rm min}\{M,12\}$, where $M={\rm min}\{|V|:  V,W\subset S^4 \; are \;disjoint \;and \;equivalent \;polybox \;codes  \;with-\\out \;twin \;pairs\}$.
\end{lemat}

\proof To prove the first part of the lemma let us suppose that it is not true. We may assume that   $V^{i,a}\neq\emptyset$ and $V^{i,b}=\emptyset$ for some $i\in [d]$. Observe that, by Corollary \ref{<3}, the polybox code $V^{i,l}_{i^c}$ is rigid for every $i\in [4]$ and $l\in S$. Therefore, exactly in the same way as in the first part of the proof of Lemma \ref{komp} we show that the codes $V^{i,a}_{i^c}$ and $V^{i,a'}_{i^c}$ are equal, and consequently there is a twin pair in the set  $V^{i,a}\cup V^{i,a'}$, which is a contradiction.  

Thus, $V^{i,l}\neq\emptyset$ for every $l\in S$ and $i\in [4]$. In the same way we show that $W^{i,l}\neq\emptyset$ for every $l\in S$ and $i\in [4]$.  

\smallskip

To prove the second part of the lemma let $d=5$. We will consider four cases. 

In the first case we assume that $U=U^{i,l}$ for some $i\in [d]$ and some $l\in S$ (thus, all the sets  $U^{i,l'},U^{i,s},U^{i,s'},l\not\in \{s,s'\}$, are empty). By (\textbf{P}) in Section 2.3, the codes $U^{i,l}$ and $P^{i,l}$ are equivalent, and then $U_{i^c}^{i,l},P_{i^c}^{i,l}\subset S^{4}$ are equivalent. Therefore, $|U^{i,l}|\geq M$. Thus, $|U|\geq M$. 

\smallskip
The second case is $U^{i,l}\neq\emptyset,U^{i,s}\neq\emptyset$ and  $U=U^{i,l}\cup U^{i,s}$ for some $i\in [d]$ and some $l,s\in S$, $s\not\in \{l,l'\}$. In the same manner as in the first case we show that  $|U^{i,l}|\geq M$ and  $|U^{i,s}|\geq M$. Thus, $|U|> M$.

\smallskip
The third case is $U^{i,l}\neq\emptyset,U^{i,l'}\neq\emptyset$ and  $U=U^{i,l}\cup U^{i,l'}$. 

If $P^{i,s}\neq \emptyset$, $s\not\in \{l,l'\}$, then, by (\textbf{C}) in Section 2.3, the codes $P^{i,s}_{i^c}$, $P^{i,s'}_{i^c}\subset S^4$ are equivalent. Thus, $|P^{i,s}|\geq M$, and then $|U|\geq M$.

If $P^{i,s}\cup P^{i,s'}=\emptyset$, then $P=P^{i,l}\cup P^{i,l'}$, and therefore the codes $U^{i,l}$ and $P^{i,l}$ are equivalent and similarly,  $U^{i,l'}$ and $P^{i,l'}$ are equivalent. Hence, $|U^{i,l}|\geq M$, $|U^{i,l'}|\geq M$ and then $|U|>M$ 

\smallskip
Finally, the fourth case is $U^{i,l}\neq\emptyset,U^{i,l'}\neq\emptyset$, $U^{i,s}\neq\emptyset$, $s\not\in \{l,l'\}$, and  $U=U^{i,l}\cup U^{i,l'}\cup U^{i,s}$. 

If $|U^{i,l}|\leq 7$ and $|U^{i,l'}|\leq 7$, then, by  Lemma \ref{komp}, the codes $U^{i,l},U^{i,l'}$ are rigid, and thus,  the codes $U_{i^c}^{i,l},U_{i^c}^{i,l'}$ are rigid.  Since the codes $U^{i,l}_{i^c}$ and $P^{i,l}_{i^c}\cup P^{i,s'}_{i^c}$ are equivalent (compare the first part of the proof of Lemma \ref{komp}) and $U^{i,l}_{i^c}$ is rigid, it follows that   $U^{i,l}_{i^c}=P^{i,l}_{i^c}\cup P^{i,s'}_{i^c}$. Then $P^{i,l}=\emptyset$, for otherwise $U^{i,l}\cap P^{i,l}\neq\emptyset$, which is a contradiction. Thus, $U^{i,l}_{i^c}=P^{i,s'}_{i^c}$. In the same way we show that $U^{i,l'}_{i^c}=P^{i,s'}_{i^c}$. Hence, $U^{i,l}_{i^c}=U^{i,l'}_{i^c}$, and then there is a twin pair in the set $U^{i,l}\cup U^{i,l'}$. A contradiction.

%In the same way as in the case of the code $V$ we show that $(U^{i,l}_{i^c}$ and $(U^{i,l'}_{i^c}$ are equivalent.  

Therefore, we may assume that $|U^{i,l}|\geq 8$ and $|U^{i,l'}|\geq 1$. Since $U^{i,s}\neq\emptyset$ and $U^{i,s'}=\emptyset$, it follows that, by (\textbf{P}), $U^{i,s}\sqsubseteq P^{i,s}$, and then, by Corollary \ref{uww}, $|P^{i,s}|\geq 5$.  Consequently, $|(\pi^i_x\cap \bigcup E(U))_{i^c}|_0\geq 8$ for every $x\in El\cap Es'$ and $|(\pi^i_y\cap \bigcup E(U))_{i^c}|_0\geq 5$ for every $y\in El'\cap Es$. By Lemma \ref{tech}, we obtain $|U|>12$. 
Thus, $|U|\geq {\rm min}\{M,12\}$ for $d=5$.

Now assume that the second part of the lemma is true for $d-1\geq 5$. To show that it is also true for $d\geq 6$ we consider the same four cases as for $d=5$. They are examined in exactly the same way as previously with only one change: In the first three cases instead of the inequalities $|U|\geq M$ obtained for $d=5$ we get, by the inductive hypothesis, $|U|\geq  {\rm min}\{M,12\}$. In the four case we obtain, as for $d=5$, $|U|>12$.

Thus, by induction on the dimension $d$, $|U|\geq  {\rm min}\{M,12\}$ for every $d\geq 5$.
\hfill{$\square$}

\medskip
Now we list the special cases in which Theorem \ref{12} holds. 
 
\begin{st}
\label{pusty}
 Let $V,W\subset S^d$ be disjoint sets which are equivalent polybox codes without twin pairs.
 
 (a) If there are words $v\in V$ and $w\in W$ such that $v_{i^c}=w_{i^c}$ for some $i\in [d]$, then $|V|\geq 12$.
 
 (b) If there are $i\in [d]$ and $l\in S$ such that $|V^{i,l}|\geq 5$ and  $|V^{i,l'}|\geq 5$, then $|V|\geq 12$.
 
 (c) If there are $i\in [d]$ and $l\in S$ such that $|V^{i,l}|=1$ and  $1\leq |V^{i,l'}|\le 4$, then $|V|\geq 12$.
 
 (d) If there are $i\in [d]$ and letters $l,s,p\in S$ such that $l\not\in \{s,s',p,p'\}$, $s\not\in \{p,p'\}$ and $V^{i,r}\cup V^{i,r'}\neq \emptyset$ for every $r\in\{l,s,p\}$, then $|V|\geq 12$. 
 
 %(e) If $d=4$ and there are $i\in [d]$ and $l\in S$ such that $V^{i,l}\neq\emptyset$ and $V^{i,l'}=\emptyset$, then $|V|\geq 12$.
   
(e) Assume that there are $i\in [d]$ and $l\in S$ such that $|V^{i,l}|=1$. If $d=4$, then $|V|\geq 12$, and if $d\geq 5$, then $|V|\geq {\rm min}\{M,12\}$, where $M$ such as in Lemma \ref{dis}.
\end{st} 
 
\noindent
{\it Proof of (a)} Since $V\cap W=\emptyset$, we have $v\in V^{i,l}$ and $w\in W^{i,s}$, where $l\neq s$. By (\textbf{Co}) in Section 2.3,  $v_{i^c}\sqsubseteq V^{i,l'}_{i^c}$ and $w_{i^c}\sqsubseteq W^{i,s'}_{i^c}$. Clearly, $v_{i^c}\not\in  V^{i,l'}_{i^c}$ and $w_{i^c}\not\in  W^{i,s'}_{i^c}$, for otherwise there are  twin pairs in $V^{i,l}\cup V^{i,l'}$ or $W^{i,s}\cup W^{i,s'}$, which is impossible.
By Corollary \ref{uww}, $|V^{i,l'}|\geq 5$ and $|W^{i,s'}|\geq 5$.  Let us assume that $|V^{i,l'}|=5$. 
Then, again by Corollary \ref{uww}, the code $V^{i,l'}$ is rigid, and thus the code $V^{i,l'}_{i^c}\subset S^{d-1}$ is rigid. Let
%Let us suppouse that there is a point $x\in Ea'$ such that
$$
\pi_x^i\cap \bigcup E(V)=\pi_x^i\cap \bigcup E(V^{i,l'})
$$
for some $x\in El'\cap Es$. Since $\pi_x^i\cap \bigcup E(W)=\pi_x^i\cap \bigcup E(V)$, we have $(\pi_x^i\cap \bigcup E(V^{i,l'}))_{i^c}=(\pi_x^i\cap \bigcup E(W))_{i^c}$, and hence, by (\ref{2d}), the polybox codes $V^{i,l'}_{i^c}$ and $\{u_{i^c}\colon u\in W\; {\rm and}\; \pi^i_x\cap \breve{u}\neq\emptyset\}$ are equivalent. As $V^{i,l'}_{i^c}$ is rigid,  $V^{i,l'}_{i^c}= \{u_{i^c}\colon u\in W\; {\rm and}\; \pi^i_x\cap \breve{u}\neq\emptyset\}$, and then $w_{i^c}\in V^{i,l'}_{i^c}$ because $\pi^i_x\cap \breve{w}\neq\emptyset$. It follows that there is a twin pair in $V^{i,l}\cup V^{i,l'}$ because $w_{i^c}=v_{i^c}$ and $v\in V^{i,l}$, which is a contradiction. Therefore, $|V^{i,l'}|\geq 6$ or $\pi_x^i\cap \bigcup E(V)\neq \pi_x^i\cap \bigcup E(V^{i,l'})$ for every $x\in El'\cap Es$. Thus
$$
|(\pi_x^i\cap \bigcup E(V))_{i^c}|_0\geq 6
$$
for every  $x\in El'\cap Es$. In the same manner we show that
$$
|(\pi_y^i\cap \bigcup E(W))_{i^c}|_0\geq 6
$$
for every  $y\in El\cap Es'$. Clearly, $(\pi_y^i\cap \bigcup E(W))_{i^c}=(\pi_y^i\cap \bigcup E(V))_{i^c}$. By Lemma \ref{tech}, $|V|\geq 12.$ 

\medskip

\noindent
{\it Proof of (b)}. %The proof is quite similar to that of the part $(a)$. 
Let $|V^{i,l}|=5$ and suppose that there is $x\in El$ such that
$$
\pi_x^i\cap \bigcup E(V)=\pi_x^i\cap \bigcup E(V^{i,l}).
$$
Since $(\pi_x^i\cap \bigcup E(V^{i,l}))_{i^c}=(\pi_x^i\cap \bigcup E(W))_{i^c}$, the polybox codes $V^{i,l}_{i^c}$ and $\{w_{i^c}\colon w\in W\; {\rm and}\; \pi_x^i\cap \breve{w}\neq\emptyset\}$ are equivalent and thus equal because, by Corollary \ref{uww}, the polybox code  $V^{i,l}_{i^c}$ is rigid. Hence,  $v_{i^c}=w_{i^c}$  for some $v\in V^{i,l}$ and some $w\in W$.
Then, by $(a)$, $|V|\geq 12$. 

Let 
$$
|V^{i,l}|\geq 6\;\; {\rm or}\;\; \pi_x^i\cap \bigcup E(V)\neq \pi_x^i\cap \bigcup E(V^{i,l})
$$
for every $x\in El$ and 
$$
|V^{i,l'}|\geq 6\;\; {\rm or} \;\; \pi_y^i\cap \bigcup E(V)\neq \pi_y^i\cap \bigcup E(V^{i,l'})
$$
for every $y\in El'$. Then, for every $x\in El$ and $y\in El'$, we have 

$$
|(\pi_x^i\cap \bigcup E(V))_{i^c}|_0\geq 6\;\;{\rm and}\;\; |(\pi_y^i\cap \bigcup E(V))_{i^c}|_0\geq 6.
$$
It follows from Lemma \ref{tech}, in which we take $l_1=l_2=l$,  that $|V|\geq 12$.

\medskip 
 
\noindent
{\it Proof of (c).} By Lemma \ref{3}, for every $w\in W$ such that $w_i\not \in \{l,l'\}$ we have $\breve{w}\cap \bigcup E(V^{i,l}\cup V^{i,l'})=\emptyset$. 
Therefore, $V^{i,l}\sqsubseteq W^{i,l}$ and $V^{i,l'}\sqsubseteq W^{i,l'}$. By Corollary \ref{uww}, $|W^{i,l}|\geq 5$ and $|W^{i,l'}|\geq 5$, and from (b) we get $|W|\geq 12$.   

\medskip
\noindent
{\it Proof of (d).} If $V^{i,r}\neq \emptyset$  and $V^{i,r'}\neq \emptyset$ for some $r\in \{l,s,p\}$, then we may assume, by $(c)$, that $|V^{i,r}\cup V^{i,r'}|\geq 4$.  Thus, if  $V^{i,r}\neq \emptyset$  and $V^{i,r'}\neq \emptyset$ for every $r\in \{l,s,p\}$, then $|V|\geq 12$.

Observe that we may assume that  
$$
W^{i,r}\cup W^{i,r'}\neq \emptyset
$$
for every $r\in \{l,s,p\}$. Indeed, if  $W^{i,r}\cup W^{i,r'}= \emptyset$ for some $r\in \{l,s,p\}$, then, by (\textbf{C}) in Section 2.3,  the polybox codes $V^{i,r}_{i^c}$ and $V^{i,r'}_{i^c}$ are equivalent. If $|V^{i,r}|\leq 5$, then, by Corollary \ref{uww} , these two codes are rigid. Thus,  $V^{i,r}_{i^c}=V^{i,r'}_{i^c}$ and consequently there is a  twin pair in $V^{i,r}\cup V^{i,r'}$, a contradiction. Therefore, $|V^{i,r}|\geq 6$, $|V^{i,r'}|\geq 6$ and then $|V|\geq 12$.

\smallskip
Suppose that  $V^{i,r}=\emptyset$ for some $r\in \{l,l',s,s',p,p'\}$. We will consider three cases.

\smallskip
In the first case we assume that $V^{i,r}\neq \emptyset$, $V^{i,r'}\neq \emptyset$ for every $r\in \{l,s\}$ and $V^{i,p}\neq\emptyset$, $V^{i,p'}=\emptyset$. Then $V^{i,p}\sqsubseteq W^{i,p}$, by (\textbf{P}) in Section 2.3 and then $|W^{i,p}|\geq 5$, by Corollary \ref{uww}. If $W^{i,r}\neq \emptyset$  and $W^{i,r'}\neq \emptyset$ for every $r\in \{l,s\}$, then, by $(c)$, we may assume that $|W^{i,r}\cup W^{i,r'}|\geq 4$ for $r\in \{l,s\}$ which gives $|W|>12$. Thus, assume first that  $W^{i,l}\neq \emptyset$, $W^{i,l'}\neq \emptyset$ and $W^{i,s}\neq \emptyset$, $W^{i,s'}=\emptyset$. If $W^{i,p'}=\emptyset$, then, by (\textbf{P}), we have  $W^{i,p}\sqsubseteq V^{i,p}$ which means, by Corollary \ref{uww}, that $|V^{i,p}|\geq 5$. Consequently, $|V|> 12$ because $|V^{i,r}\cup V^{i,r'}|\geq 4$ for $r\in \{l,s\}$. If  $W^{i,p'}\neq\emptyset$, then only one distributions of words in $W$ has to be considered: $|W^{i,l}\cup W^{i,l'}|=4$,  $|W^{i,s}|=1$, and $|W^{i,p}|=5$, $|W^{i,p'}|=1$. It follows from  (\textbf{V}) in Section 2.3 that $|V^{i,p}|\geq 4$, and from (\textbf{P}) and Corollary \ref{uww} we obtain $|V^{i,s}|\geq 5$. Since  $|V^{i,l}\cup V^{i,l'}|\geq 4$, we have $|V|>12$. Assume now $W^{i,l}\neq \emptyset$, $W^{i,l'}=\emptyset$ and $W^{i,s}\neq \emptyset$, $W^{i,s'}=\emptyset$. Then, by (\textbf{P}) and Corollary \ref{uww}, $|V^{i,l}|\geq 5$ and $|V^{i,s}|\geq 5$. Since $V^{i,l'}\neq\emptyset$ and $V^{i,s'}\neq\emptyset$, we have $|V|\geq 12$.

\smallskip
In the second case we assume that $V^{i,l}\neq \emptyset$, $V^{i,l'}\neq \emptyset$, $V^{i,s}\neq\emptyset$, $V^{i,s'}=\emptyset$ and $V^{i,p}\neq\emptyset$, $V^{i,p'}=\emptyset$. Then $|W^{i,s}|\geq 5$ and  $|W^{i,p}|\geq 5$,  by (\textbf{P}) and Corollary \ref{uww}. If $W^{i,s'}=\emptyset$ and $W^{i,p'}=\emptyset$, then, by (\textbf{P}) and Corollary \ref{uww}, $|V^{i,s}|\geq 5$ and  $|V^{i,p}|\geq 5$. Then $|V|\geq 12$ because  $V^{i,l}\neq \emptyset$, $V^{i,l'}\neq \emptyset$. Therefore, we may assume that   $W^{i,s'}\neq\emptyset$ or $W^{i,p'}\neq\emptyset$. Since  $|W^{i,s}|\geq 5$, $|W^{i,p}|\geq 5$ and $W^{i,l}\cup W^{i,l'}\neq \emptyset$, we have $|W|\geq 12$.

\smallskip
Finally, in the third case we assume that $V^{i,r}\neq \emptyset$  and $V^{i,r'}=\emptyset$ for every $r\in \{l,s,p\}$. Then, by (\textbf{P}) and Corollary \ref{uww}, $|W^{i,r}|\geq 5$ for  every $r\in \{l,s,p\}$, and consequently $|W|>12$.
%\hfill{$\square$}

%\medskip
%\noindent
%{\it Proof of (e).} By Lemma \ref{<3} every polybox code $V\subset S^d$ without twin pairs, where $d\leq 3$, is rigid. Thus, by Lemma \ref{dis}, we may assume that there are   at least two letters $p,s\in S$ such that $V^{i,r}\neq\emptyset$ for every $r\in \{p,p',s,s'\}$. Since $V^{i,l}\neq\emptyset$ and  $V^{i,l'}=\emptyset$, we have $l\not\in \{p,p'.s,s'\}$. Thus, by $(d)$, $|V|\geq 12$. 

\medskip
\noindent
{\it Proof of (e).} Suppose that $V^{i,l}\neq\emptyset$ and $W^{i,l}\neq\emptyset$ for every $i\in [d]$ and every $l\in S$. 

Let $V^{i,l}=\{u\}$. Assume first that $V^{i,l}\not \sqsubseteq W^{i,l}$. Then there is $w\in W\setminus ( W^{i,l}\cup W^{i,l'})$ such that $\breve{w}\cap \breve{u}\neq\emptyset$, and thus, by (\textbf{P}), $\breve{w}\cap \bigcup E(V^{i,l'})\neq\emptyset$.  %Thus $V^{i,l'}\not \sqsubseteq W^{i,l'}$. 
Since $w\sqsubseteq V$, the set  $\breve{w}\cap \breve{u}\cup \breve{w}\cap  \bigcup E(V^{i,l'})$ is an $i$-cylinder in the $d$-box $\breve{w}$. %(as  $\{\breve{w}\cap \breve{v}\neq\emptyset: v\in V\}$ is a suit for $\breve{w}$). 
Therefore, 
$$
(\breve{w}\cap \breve{u})_{i^c} =\bigcup\{(\breve{w} \cap \breve{v})_{i^c}\colon v\in V^{i,l'}\},
$$
and thus 
\begin{equation}
\label{f}
Ew_j\cap Ev_j\subseteq Ew_j\cap Eu_j
\end{equation}
for every $j\in[d]\setminus \{i\}$ and every $v\in V^{i,l'}$ for which $(\breve{w} \cap \breve{v})_{i^c}\neq\emptyset$. %(compare the proof of Lemma \ref{dis}). 

If $w_j=u_j$ for every $j\in [d]\setminus\{i\}$, then $w_{i^c}=u_{i^c}$ and, by $(a)$, $|V|\geq 12$. 

Let $w_k\neq u_k$ for some $k\in [d]\setminus\{i\}$.  Then $v_k=u_k$ for every $v\in V^{i,l'}$ with $(\breve{w} \cap \breve{v})_{i^c}\neq\emptyset $. Indeed, if $v_k\neq u_k$ for some $v$, then, by (\ref{dkostki}), $Ew_k\cap Ev_k\not\subseteq Ew_k\cap Eu_k$, which contradicts (\ref{f}). Therefore, $v\in V^{k,u_k}$ for every $v\in V^{i,l'}$ such that $(\breve{w} \cap \breve{v})_{i^c}\neq\emptyset$. Since the set of boxes $\{(\breve{w} \cap \breve{v})_{i^c}\colon v\in V^{i,l'}\}$ is a suit for the $(d-1)$-box $(\breve{w}\cap \breve{u})_{i^c}$, which, by Lemma \ref{=c} and the fact that $V$ is a twin pair free, does not contain a twin pair, it has to contain, by Lemma \ref{5} $(a)$, at least five boxes. Therefore, $|V^{k,u_k}|\geq 5$. But $u\in V^{k,u_k}$, and thus $|V^{k,u_k}|\geq 6$. 

If $|V^{k,u'_k}|\geq 2$, then assuming,  by $(c)$, that $|V^{k,l}\cup V^{k,l'}|\geq 4$, where $l\not\in \{u_k,u_k'\}$, (recall that it was assumed that $V^{i,l}\neq\emptyset$ for every $i\in [d]$ and every $l\in S$) we get $|V|\geq 12$. 

Assume that  $V^{k,u'_k}=\{p\}$ for some $p\in S^d$. Clearly, $\breve{w}\cap \breve{p}\neq\emptyset$ because $w_k\not\in \{u_k,u_k'\}$ and $\breve{w}\cap \breve{v}\neq\emptyset$ for at least six $v\in V^{k,u_k}$. Since the set $\breve{w}\cap \bigcup E( V^{k,u_k}\cup  V^{k,u_k'})$ is a $k$-cylinder in $\breve{w}$, we have 
$$
(\breve{w}\cap \breve{p})_{k^c}=\bigcup \{(\breve{w} \cap \breve{v})_{k^c}\colon v\in V^{k,u_k}\}.
$$ 
%(if $V^{k,u'_k}=\emptyset$, then, by $(d)$, we have $|V|\geq 12$). 
In the same way as above, we show that $w_{k^c}=p_{k^c}$, and then, by $(a)$, $|V|\geq 12$ or $|V^{m,p_m}|\geq 7$ for some $m\in [d]$. This last inequality follows from that fact that the set $\{(\breve{w} \cap \breve{v})_{k^c}\colon v\in V^{k,u_k}\}$ contains at least six boxes. Using the same arguments as before we show that there is $m\in [d]$ such that $w_m\neq p_m$ and $v_m=p_m$ for every $v\in  V^{k,u_k}$ such that $(\breve{w} \cap \breve{v})_{k^c}\neq\emptyset$, which  gives $|V^{m,p_m}|\geq 6$. Since $p\in V^{m,p_m}$ and $p\notin  V^{k,u_k}$, we have $|V^{m,p_m}|\geq 7$.  
%there is $\bar{w}\in W$ such that $(\bar{w})_{k'}=(\bar{u})_{k'}$, where $\bar{u}\in V^{k,u'_k}$, which gives us $|V|\geq 12$ or $|V^{m,\bar{u}_m}|\geq 7$ for some $m\in [d]$. 
By  $(c)$, we may assume that $|V^{m,l}\cup V^{m,l'}|\geq 4$ for some $l\in S\setminus\{p_m,p'_m\}$. Since $|V^{m,p'_m}|\geq 1$, it follows that $|V|\geq 12$.

We have shown that if $V^{i,l}\not \sqsubseteq W^{i,l}$, then $|V|\geq 12$.

Let now $V^{i,l}\sqsubseteq W^{i,l}$ and $V^{i,l'}\sqsubseteq W^{i,l'}$ (note that if $V^{i,l'}\not \sqsubseteq W^{i,l'}$, then $V^{i,l}\not \sqsubseteq W^{i,l}$). Since $V^{i,l}\neq\emptyset$ and $V^{i,l'}\neq\emptyset$, it follows that, by Corollary \ref{uww}, $|W^{i,l}|\geq 5$ and $|W^{i,l'}|\geq 5$. Thus, by $(b)$, we obtain $|W|\geq 12$.  

Let now $V^{i,l}=\emptyset$ for some $i\in [d]$ and some $l\in S$. By $(d)$ we assume that $S=\{a,a',b,b'\}$. By Lemma \ref{dis}, the assumption   $V^{i,l}=\emptyset$ leads to a contradiction for $d=4$, and for $d\geq 5$ it implies that $|V|\geq {\rm min}\{M,12\}$.  
\hfill{$\square$}

\medskip
If $v\in S^d$, and $\sigma$ is a permutation of the set $[d]$, then $\sigma^*(v)=v_{\sigma(1)}\ldots v_{\sigma(d)}$. For every $i\in [d]$ let $h_i:S\rightarrow S$ be a bijection such that $h_i(l')=(h_i(l))'$ for every $l\in S$, and let $h:S^d\rightarrow S^d$ be defined by the formula $h(v)=h_1(v_1)\ldots h_d(v_d)$. We say that polybox codes $P,Q\subset S^d$ are {\it isomorphic} if there are $\sigma$ and $h$ such that $Q=\{h_1(v_{\sigma(1)})\ldots h_d(v_{\sigma(d)}): v\in P\}$. %The composition $h\circ \sigma$ is an {\it isomorphism} between $P$ and $Q$. It follows from the definition of the isomorphic polybox codes that if $V$ and $W$ are disjoint, equivalent and do not contain a twin pair, then  the isomorphic codes $h\circ \sigma(V)$ and  $h\circ \sigma(W)$ are also disjoint, equivalent and do not contain a twin pair.

\medskip
Now we can prove Theorem \ref{12}. 

\medskip
{\it The proof of Theorem \ref{12}}. By  Corollary \ref{<3}, we may assume that $d\geq 4$. 

\medskip
Let $d=4$. By Statement \ref{pusty} $(d),(e)$ and Lemma \ref{dis} we make the following assumptions:

\;\;\; {\bf (A)} \;  $S=\{a,a',b,b'\}$, $|V^{i,l}|\geq 2$ and  $|W^{i,l}|\geq 2$ for every $i\in [d]$ and $l\in S$.

\medskip

For $w\in W$ let $V_w\subseteq V$ be such that $w\sqsubseteq V_w$ and $\breve{w}\cap \breve{v}\neq\emptyset$ for all $v\in V_w$. 
%Since $V$ does not contain a twin pair, by Lemma \ref{=c}, for every $w\in W$ the code $V_w$ does not contain a twin pair.

If $|V_w|<5$ for some $w\in W$, then, by Corollary \ref{uww}, there is a twin pair in $V_w$, which is impossible. If $|V_w|=6$ for some $w\in W$, then, by Corollary \ref{uww}, there are $v\in V$ and $i\in [d]$ such that $v_{i^c}=w_{i^c}$. Thus, by Statement \ref{pusty} (a),  $|V|\geq 12$. Finally, if $|V_w|\geq 10$ for some $v\in V$, then $|V|\geq 12$ because $V_w\cap V^{i,w_i'}=\emptyset$ for every $i\in [d]$ and, by (A), $| V^{i,w_i'}|\geq 2$.  Therefore in what follows we assume that $|V_w|\in \{5,7,8,9\}$. Since $v\sqsubseteq W$ for $v\in V$ and $W$ does not contain twin pairs, it follows from (\ref{00indd}) that there are $w^1,w^2\in W$ such that $|\{i:w^1_i=(w^2_i)'\}|=3$ and $w^1_j=w^2_j$ for $j\not\in \{i:w^1_i=(w^2_i)'\}$. 
For simplicity we assume that $w^1=bbbb$ and $w^2=b'b'b'b$. By ${\rm Cov}_{w^1}$ we denote the family of all polybox codes $V_{w^1}$ such that $|V_{w^1}|\in \{5,7,8,9\}.$ The family ${\rm Cov}_{w^1}$ can by easily computed. We can simplify computations using (\ref{2d}). Recall that, by (\ref{2d}), we have $\sum_{v\in V_{w^1}}g(v,w^1)=16$. Note that we may assume that $g(v,w^i)\leq 4$ for every $v\in V_{w^i}$ and $i=1,2$ because if  $g(v,w^i)=8$ for some $v\in V_{w^i}$ and some $i\in \{1,2\}$, then $v_{j^c}=w^i_{j^c}$ for some $j\in [4]$. Then, by Statement \ref{pusty} $(a)$, $|V|\geq 12$.

Let us consider the system of the equations:

\vspace{-3mm}
\begin{equation}
\label{eqqq}
4x+2y+z=16, \; x+y+z=k,
\end{equation}
where $k\in \{7,8,9\}$ and $x,y,z\in \{0,1\ldots \}$. For every $k\in \{7,8,9\}$ this system has only three solutions: 
$$
k=7:\;  x =3,\;y = 0,\;z =4;\;\; x= 1,\; y =6,\; z =0;\;\; x =2,\; y = 3,\; z =2.
$$
$$
k=8:\;\; x =2,\;y = 2,\; z =4;\;\;x= 0,\;y =8,\;z =0;\;\;x =1,\;y = 5\;z =2.
$$
$$
k=9:\;\;x= 0,\;y =7,\;z =2;\;\; x =1,\;y = 4,\;z =4;\;\; x =2,\; y = 1,\; z =6.
$$
(The form of $V_{w^1}$ in the case $|V_{w^1}|=5$ was described in Corollary \ref{uww}.)
In the last comments in Section 2.8 we showed how to use these solutions to decipher partially the structure of $V_{w^1}$. Let ${\rm Cov}_{w^1}(x,y,z)$ be the family of all codes $V_{w^1}$ such that every code $V_{w^1}$ (which does not contain a twin pair) contains exactly $z$ words without the letter $b$, $y$ words with one letter $b$ and $x$ words with two letters $b$. Simple computations show that

\vspace{-6mm} 
\begin{center}
\begin{displaymath}
\begin{tabular}{|l|l|}
%\hline
%\multicolumn{2}{|c|}{If $|V|=|W|=3$ and $d=3$, then} \\ 
%Total & Non-isomorphic codes \\ \cline{1-2}
\hline
 %\multicolumn{2}{|c|} 
 $|{\rm Cov}_{w^1}(3,2,0)|=32$ & $B_{w^1}=\{aaab,a'a'a'b,baa'b,a'bab,aa'bb\}$ \\ \cline{1-2}
 $$ & $C_{w^1}=\{aaaa, aaa'b, aa'a'a', a'aba,a'ba'a',ba'ba,bbaa'\}$  \\ 
 $|{\rm Cov}_{w^1}(2,3,2)|=576$                               &  $D_{w^1}=\{aaab,aa'ba,a'baa',aaa'a,aa'aa',a'bba,bba'a'\}$ \\ \cline{1-2}
 $|{\rm Cov}_{w^1}(1,5,2)|=192$ & $E_{w^1}=\{aaab,aba'a', a'aaa,a'aba',a'a'a'a',baa'a,ba'aa',ba'ba\}$ \\ \cline{1-2}
  $|{\rm Cov}_{w^1}(0,8,0)|=8$ & $F_{w^1}=\{aaab, aa'ba', aba'a, a'aba, a'a'a'b, a'baa', baa'a', ba'aa\}$ \\ \cline{1-2}
  $|{\rm Cov}_{w^1}(1,4,4)|=48$ & $G_{w^1}=\{aaa'a,aaba',aa'a'a',aa'ba,a'aaa',a'a'aa,a'ba'b,baaa,ba'aa'\}$ \\ \cline{1-2}
  \hline
  %\multicolumn{2}{|c|}{$A=\{i_1,i_2\}$, $(V)_{A}=(W)_{A}=(\{s\})_{A}$, $s\in (*S)^d$ is fixed}\\
\end{tabular}
\end{displaymath}
\end{center}

\smallskip
\noindent{\footnotesize 
Tab. 1. The structures of the polybox codes $V_{w^1}$ without twin pairs  for $|V_{w^1}|\in \{5,7,8,9\}$, where $g(v,w^1)\leq 4$ for every $v\in V_{w^1}$. 
  }

\medskip
where ${\rm Cov}_{w^1}(x,y,z)=\emptyset$ for every solution $(x,y,z)$ of (\ref{eqqq}) which is different from $(3,2,0),(2,3,2),\\(1,5,2),(0,8,0)$ and $(1,4,4)$. The codes $B_{w^1},\ldots ,G_{w^1}$ are non-isomorphic codes in the corresponding ${\rm Cov}_{w^1}(x,y,z)$, but in our case for every $i\in [d]$ we consider only two bijection $h_i$: The first is the identity and the second is $h_i(a)=a',h_i(a')=a,h_i(b)=b$ and $h_i(b')=b'$ (compare the definition of the isomorphic codes before the proof). Making the substitution $b\rightarrow b'$ at the first three position in every $v\in P$ for every  $P\in  {\rm Cov}_{w^1}(x,y,z)$ we obtain  ${\rm Cov}_{w^2}(x,y,z)$. 
%Having the above solutions and knowing the structure of $V_{w^1}$ when $|V_{w^1}|=5$, for every $k\in \{5,7,8,9\}$ 

For every $k\in \{5,7,8,9\}$ we will consider the system  $|V_{w^1}|=k$ and  $|V_{w^2}|\geq k$. Recall that if a word $v\in S^4$ contains a letter $b$ at the the position $1,2$ or $3$, then $v\not \in V_{w^2}$ and similarly,  if $v$ contains a letter $b'$ at the the position $1,2$ or $3$, then $v\not \in V_{w^1}$.

\smallskip
Let $|V_{w^1}|=5$ and $|V_{w^2}|=5$.  Then, there is $i\in [4]$ such that $v_i=b$ for every $v\in V_{w^1}$ (see Table 1). 

If $i\in \{1,2,3\}$, then $V_{w^1}\cap V_{w^2}=\emptyset$, and consequently  $|V_{w^1}\cup V_{w^2}|=10$. Since $V^{4,b'}\cap (V_{w^1}\cup V_{w^2})=\emptyset$ and, by (A), $|V^{4,b'}|\geq 2$, we have $|V|\geq 12$. 

If $i=4$, then $V_{w^1}\subset V^{4,b}$. Since $|V_{w^2}|=5$, we can assume, by just considered case, that $v_4=b$ for all $v\in V_{w^2}$. Then $|V^{4,b}|\geq 8$ because for every $v\in V_{w^2}$ with $v_i=b'$ for some $i\in \{1,2,3\}$ we have $v\in V^{4,b}$ and $v\not\in V_{w^1}$. By (A), we have $|V^{4,l}|\geq 2$ for $l\in \{a,a',b'\}$, and thus $|V|> 12$.

\smallskip
Let $|V_{w^1}|=5$ and $|V_{w^2}|\geq 7$. Observe that  $V_{w^2}\not \in {\rm Cov}_{w^2}(2,3,2)\cup {\rm Cov}_{w^2}(1,5,2)\cup {\rm Cov}_{w^2}(1,4,4)$ because every code in this union contains a word $u$ without the letter $b$, and every such word $u$ has to belong to $V_{w^1}$, as, by (\ref{dkostki}), $\breve{u}\cap \breve{w}^1\neq\emptyset$. On the other hand, $u\not\in V_{w^1}$ because every word in the set $V_{w^1}$ contains at least one letter $b$. Thus, $V_{w^2} \in {\rm Cov}_{w^2}(0,8,0)$. Since every  code in $ {\rm Cov}_{w^2}(0,8,0)$ contains exactly six words with the letter $b'$ at the positions $1,2,3$, we have $|V_{w^1}\cup V_{w^2}|\geq 11$. Then $|V|>12$ because $|V^{4,b'}|\geq 2$ and $(V_{w^1}\cup V_{w^2})\cap V^{4,b'}=\emptyset$.

\smallskip
Let $|V_{w^1}|\geq 7$ and $|V_{w^2}|\geq 7$. An inspection of the codes $C_{w^1},\ldots ,G_{w^1}$ shows that every polybox code $V_{w^2}$ from the family  ${\rm Cov}_{w^2}(2,3,2)\cup {\rm Cov}_{w^2}(1,5,2)\cup  {\rm Cov}_{w^2}(0,8,0)\cup {\rm Cov}_{w^2}(1,4,4)$ contains at least three words  with the letter $b'$ at the positions $1,2,3$. Thus,  in every such $V_{w^2}$ there are at least three words which do not belong to $V_{w^1}$. Since  $|V_{w^1}|\geq 7$, we have $|V_{w^1}\cup V_{w^2}|\geq 10$. In the same way as above we obtain $|V|\geq 12$. This completes the proof of the theorem for $d=4$.

\medskip 
Let $d\geq 5$. By Statement \ref{pusty} $(c)$, we may assume that $S=\{a,a',b,b'\}$. Since, by the first part of the proof,  $M\geq 12$, where $M$ is defined in Lemma \ref{dis}, we assume, by Lemma \ref{dis}, that $V^{i,l}\neq\emptyset$ for every $i\in [d]$ and every $l\in S$. By Statement \ref{pusty} $(e)$ we may assume that $|V^{i,l}|\geq 2$ for every $i\in [d]$ and every $l\in S$. Thus we assume that, 
%By Statement \ref{4} we need to prove the theorem for $d\geq 5$ and by Statement \ref{pusty} we can assume that $V\subset \{a,a'b,b'\}^d$, $|V^{i,l}\cup V^{i,l'}|\geq 4$ and $V^{i,l}\neq\emptyset$ for every $i\in [d]$ and $l\in \{a,b,a',b'\}$.
\begin{equation}
\label{7}
|V^{i,l}\cup V^{i,s}|\leq  7,
\end{equation}
for every $l,s\in \{a,a',b,b'\}$ and every $i\in [d]$ because if $|V^{i,l}\cup V^{i,s}|\geq  8$ for some $l,s\in \{a,a',b,b'\}$ and some $i\in [d]$ then $|V|\geq 12$ as  $|V^{i,l'}|\geq  2$  and $|V^{i,s'}|\geq  2$ .

Assume that there are $i\in [d]$ and two letters, say $a$ and $b$, such that there are no $i$-siblings $u$, $v$ in $V$ such that $u_i=a$ and $v_i=b$. Then, by (\ref{7}) and Lemma \ref{komp}, the polybox code  $(V^{i,a}\cup V^{i,b})_{i^c}$ is rigid. Since the codes $(V^{i,a}\cup V^{i,b})_{i^c}$ and $(W^{i,a}\cup W^{i,b})_{i^c}$ are equivalent (as we have $\pi^i_x\cap \bigcup E(V^{i,a}\cup V^{i,b})=\pi^i_x\cap \bigcup E(W^{i,a}\cup W^{i,b})$, they are equal. Thus,  there are  $v\in V^{i,a}$ and $w\in W^{i,b}$ such that $v_{i^c}=w_{i^c}$. By Statement \ref{pusty} $(a)$, $|V|\geq 12$.

Therefore, in what follows we assume that for every $i\in [d]$ and every two letters $l,s\in \{a,a'b,b'\}$ such that $l\not\in \{s,s'\}$ there are $i$-siblings $u$ and $v$ in $V$ such that $u_i=l$ and $v_i=s$. In particular, there are at least $4d$ edges in the set $\ka E$, where $G=(V,\ka E)$ is a graph of siblings on $V$ (see Section 2.9).   

Let $u^0,v^0\in V$ be such that 
$$
d(v^0)+d(u^0)=\max\{d(v)+d(u):v,u\in V\;{\rm and}\; v,u\;{\rm are}\;{\rm adjacent}\}.     
$$ 
Let $d=5$. Observe that we may assume that $d(v^0)+d(u^0)\leq 8$ because if $d(v^0)+d(u^0)\geq 9$, then it follows from (\ref{nmm3}) or  (\ref{nm3})  that there are $i\in [d]$ and $l\in S$ such that $|V^{i,l}\cup V^{i,l'}|\geq 8$ which contradicts the assumption (\ref{7}).

Let $d(v^0)+d(u^0)= 8$, and let $N(u^0)$ and $N(v^0)$ be the sets of all neighbors of $u^0$ and $v^0$, respectively. Taking into account (\ref{7}),
%Since  $|V^{j,l}\cup V^{j,l'}|\geq 4$ for $l\in \{a,b\}$ we can assume that  $|V^{j,l}\cup V^{j,l'}|\leq 7$ for $l\in \{a,b\}$. 
it can be easily shown (in the similar manner as in the proof of Lemma \ref{sybb}) that there are $j,k\in [5], k\neq j$, and $l,s\in \{a,b\}$
%the sets $(V^{j,a}\cup V^{j,a'})\cap (N(v^0)\cup N(u^0))$ and $(V^{k,a}\cup V^{k,a'})\cap (N(v^0)\cup N(u^0))$, where $k\neq j$ 
such that        
$$
|(V^{j,l}\cup V^{j,l'})\cap (N(v^0)\cup N(u^0))|= 7
$$
and
$$
|(V^{k,s}\cup V^{k,s'})\cap (N(v^0)\cup N(u^0))|= 7.
$$
Without loss of generality we can take $l=s=a$ because $k\neq j$.
 
Since $|V^{j,b}\cup V^{j,b'}|\geq 4$ and $|V^{k,b}\cup V^{k,b'}|\geq 4$, there are at least three words $x,y,z$ in the set $(V^{j,b}\cup V^{j,b'})\cap (V\setminus(N(u^0)\cup N(v^0)))$ and at least three words $\bar{x},\bar{y},\bar{z}$ in the set $(V^{k,b}\cup V^{k,b'})\cap (V\setminus(N(u^0)\cup N(v^0)))$. 
%Then $x_j,y_j,z_j\in \{b,b'\}$ and  $\bar{x}_k,\bar{y}_k,\bar{z}_k\in \{b,b'\}$. 
If $\{x,y,z\}\neq \{\bar{x},\bar{y},\bar{z}\}$, then $|V|\geq 12$. Let us assume that $\{x,y,z\}= \{\bar{x},\bar{y},\bar{z}\}$ and $|V|=11$. Then $x_j,y_j,z_j\in \{b,b'\}$ and $x_k,y_k,z_k\in \{b,b'\}$. On the other hand, the vertices $u\in (N(v^0)\cup N(u^0))\setminus (V^{j,a}\cup V^{j,a'})$ and $v\in (N(v^0)\cup N(u^0))\setminus (V^{k,a}\cup V^{k,a'})$ are such that $u_j\in \{b,b'\}$ and $v_k\in \{b,b'\}$. Note that  $u_k\in \{a,a'\}$ and $v_j\in \{a,a'\}$, for otherwise $u_j,u_k\in \{b,b'\}$ or $v_j,v_k\in \{b,b'\}$  and then $u\not\in N(v^0)\cup N(u^0)$ or $v\not\in N(v^0)\cup N(u^0)$ which is not true. 
Since $w_j,w_k\in \{a,a'\}$ for every $w\in (N(v^0)\cup N(u^0))\setminus \{u,v\}$, it follows that if $x,y$ and $z$ %can be adjacent only to  $u$ or $v$.
are joined with some vertices from the set $N(u^0)\cup N(v^0)$, then these vertices must be  $u$ or $v$. 
Assume without loss of generality that $u_j=b$ and $v_j=a$. Since $x,y,z$ can be joined only with $u$ or $v$, there are no $j$-siblings $p$, $q$ in $V$ with $p_j=b'$ and $q_j=a'$, which contradicts the assumption on $i$-siblings in $V$. 

If $d(v^0)+d(u^0)\leq 7$, then, by Lemma \ref{graph}, we have $d(G)\leq 7/2$. As $d(G)|V|=2|\ka E|$ and $2|\ka E|\geq 40$, we have $|V|>11.$

Let $d\geq 6$. For the same reason as for $d=5$ we assume that $d(v^0)+d(u^0)\leq 8$. Then it follows from Lemma \ref{graph} that $d(G)\leq 4$. Since now $2|\ka E|\geq 48$, we have $|V|\geq 12.$    
\hfill{$\square$}

\begin{wn}
\label{w11} 
If $V,W\subset S^d$ are equivalent polybox codes without twin pairs and $|V|\leq 11$, then $W=V$.
\end{wn}
\proof By Theorem \ref{12}, $V\cap W\neq\emptyset$. Suppose that $V\setminus V\cap W\neq\emptyset$. Since $V$ and $W$ are equivalent, $V\setminus V\cap W$ and $W\setminus V\cap W$ are equivalent (as $\bigcup\{\breve{v}: v\in V\setminus V\cap W\}=\bigcup\{\breve{w}: w\in W\setminus V\cap W\}$), and thus, again by Theorem \ref{12},  $(V\setminus V\cap W)\cap  (W\setminus V\cap W)\neq\emptyset$, a contradiction. Hence,  $V\setminus V\cap W=\emptyset$, and then $V=W$.
\hfill{$\square$}

\begin{wn}
\label{Sf}
For every $d$-box $X$ and every polybox $F\subset X$, if there is a proper suit $\ka F$ for $F$ which does not contain a twin pair and $|F|_0\leq 11$, then every proper suit $\ka G\neq \ka F$ for $F$ contains a twin pair or $\ka G=\ka F$.
\end{wn}
\proof Let $V$ be a code for $\ka F$ such that $\ka F$ is an exact realization of $V$ (see Section 2.6). Since $|V|\leq 11$ and $V$ does not contain a twin pair, by Corollary \ref{w11}, if $G$ is a polybox code which is equivalent to $V$, then $G$ contains a twin pair or $G=V$. Thus, if $\ka G$ is a realization of $G$, then there is a twin pair in $\ka G$ or $\ka G=\ka F$. \hfill{$\square$} 

\begin{uw}
{\rm The estimation given in Theorem \ref{12} is optimal. There are two equivalent polybox codes $V,W\subset \{a,a',b,b'\}^4$ both without twin pairs and $V\cap W=\emptyset$ (\cite{LS2}). These codes were used by Lagarias and Shor\cite{LS1} and later on by Mackey\cite{M} to construct the counterexamples to Keller's cube tiling conjecture. In the context of this conjecture one of these codes was given first by Corr\'adi and Szab\'o  in \cite{CS2}, as an example of the maximum clique in a 4-dimensional Keller graph.  
}
\end{uw}

\section{Twin pairs in cube tilings of $\er^d$}

From Theorem \ref{12} we obtain the following result.

\begin{tw} 
\label{cylin}
Let $W=W^{i, l_1}\cup W^{i, l'_1}\cup \cdots \cup W^{i, l_k}\cup W^{i, l'_k}\subset S^d$ be a partition code, where $W^{i,l}=\{w\in W: w_i=l\}$ for $l\in S$, $i\in [d]$, $l_j\not\in \{l_n,l_n'\},j\neq n$ and $W^{i, l_j}\cup W^{i, l'_j}\neq\emptyset$ for $j\in [k]$.  
%and $W^{i, l'_j}=\{w=w_1\cdots w_d\in W: w_i=l'_j\}$. 
If $k> 2^{d-3}/3$, then there is a twin pair in $W$.
\end{tw}
\proof Since $W$ is a partition code, for every $j\in [k]$ the polybox codes $W^{i,l_j}_{i^c}, W^{i, l'_j}_{i^c}\subset S^{d-1}$ are equivalent. 
(For every $j\in [k]$ the set  $\bigcup E(W^{i,l_j}\cup W^{i, l'_j})$ is an $i$-cylinder in the $d$-box $(ES)^d$, and thus  $\bigcup E(W^{i,l_j}_{i^c})=\bigcup E(W^{i, l'_j}_{i^c})$. By (\ref{2d}), the codes $W^{i,l_j}_{i^c}$ and $W^{i, l'_j}_{i^c}$ are equivalent.) 
%(For any realization $f(W)$, the set $f(W^{i,l_j}\cup W^{i, l'_j})$ is an $i$-cylinder, and thus  $\bigcup f((W^{i,l_j}_{i^c})=\bigcup f((W^{i, l'_j}_{i^c})$, which, by the definition, means that $(W^{i,l_j}_{i^c}, (W^{i, l'_j}_{i^c}$ are equivalent.) 
Thus $|W^{i,l_j}|=|W^{i, l'_j}|$, and then $\sum_{j=1}^k|W^{i,l_j}|=2^{d-1}$,  as $|W|=2^d$.  By the assumption on the number $k$,  there is at least one $j\in [k]$ such that
$$
|W^{i, l_j}|\leq 11.
$$
If the polybox codes $W^{i, l_j}_{i^c}$ and $W^{i, l'_j}_{i^c}$ do not contain twin pairs, then, by Corollary \ref{w11}, they are  equal. Therefore the polybox code $W^{i, l_j}\cup W^{i, l'_j}$ consists of twin pairs. If  $W^{i, l_j}_{i^c}$ contains a twin pair $w_{i^c}$, $v_{i^c}$, then $W^{i, l_j}$ contains the twin pair $w=w_1\ldots w_{i-1}l_jw_{i+1}\ldots w_d$, $v=v_1\ldots v_{i-1}l_jv_{i+1}\ldots v_d.$
\hfill$\square$

\medskip
From Theorem \ref{cylin} we obtain a general theorem on twin pairs in cube tilings of $\er^d$.
\begin{tw}
\label{tll}
Let $[0,1)^d+T$ be a cube tiling of $\er^d$. If $r^+(T)> 2^{d-3}/3$, then there is a twin pair in the tiling  $[0,1)^d+T$.
\end{tw}  

\proof  Let $x\in \er^7$ and $i\in [7]$ be such that $|L(T,x,i)|=r^+(T)$.  In Section $2.2$ we showed that the  family of boxes $\ka F_x=\{([0,1]^d+x)\cap ([0,1)^d+t)\neq\emptyset:t\in T\}$ is a minimal partition of the $d$-box $[0,1]^d+x$. Let $W=W^{i, l_1}\cup W^{i, l'_1}\cup \cdots \cup W^{i, l_k}\cup W^{i, l'_k}$ be a partition code such that $\ka F_x$ is an exact realization of $W$ (this code can be obtained in the manner described at the end of Section 2.6). Note that $|L(T,x,i)|=k$, that is   $|L(T,x,i)|$ is the number of all $i$-cylinders (such as in the representation (\ref{rep})) in the partition $\ka F_x$. Indeed,   $t_i\in L(T,x,i)$ if and only if there is $t'\in T$ such that $([0,1)^d+t')\cap ([0,1]^d+x)\neq\emptyset$ and $t'_i-t_i=1$. Since $k> 2^{d-3}/3$, by Theorem \ref{cylin}, there is a twin pair in $W$, and thus there is a twin pair in $\ka F_x$. Consequently, there is a twin pair in the tiling  $[0,1)^d+T$.  

\hfill{$\square$}

\medskip
Now we can prove Theorem \ref{keli}

\medskip
{\it Proof of Theorem \ref{keli}} By Theorem \ref{tll}, there is a twin pair in every cube tiling  $[0,1)^7+T$ of $\er^7$ such that $r^+(T)>16/3$. \hfill{$\square$}    

%\begin{wn} 
%\label{kee}
%Keller's conjecture is true for a cube tiling $[0,1)^7+T$ of $\er^7$ for which there are $x\in \er^d$ and $i\in [d]$ such that the set $L(T,x,i)$ contains at least six elements. \hfill{$\square$}
%\end{wn}

\begin{wn} 
\label{ke2}
If $[0,1)^7+T$  is a counterexample to Keller's conjecture in dimension seven, then $r^-(T),r^+(T)\in \{3,4,5\}$.
\end{wn}
\proof By Theorem \ref{keli}, we assume that  $r^-(T)\leq 2$. We will use somewhat different arguments to show that there is a twin pair in  $[0,1)^7+T$, than that used in Section 1. Since $r^-(T)\leq 2$, there is $x\in \er^7$ such that $|L(T,x,i)|\leq 2$ for every $i\in [7]$. Thus, the minimal partition $\ka F_x=\{([0,1]^d+x)\cap ([0,1)^d+t)\neq\emptyset:t\in T\}$ contains at most two $i$-cylinders (such as in the representation (\ref{rep})) for every $i\in [7]$, and therefore a partition code whose exact realization is $\ka F_x$ can be written in the alphabet $S=\{a,a',b,b'\}$. More precisely, for every $i\in [7]$ there are at most two numbers $a_i\leq b_i$ in $L(T,x,i)$. For every box $K\in \ka F_x$ and $i\in [7]$ the set $K_i$ can be of the form: $[x_i,1+a_i),[x_i,1+b_i),[1+a_i,1+x_i]$ or $[1+b_i,1+x_i]$. Let $A$ consist of all $i\in [7]$ for which $a_i=b_i$. Let $h_i([x_i,1+a_i))=a$ and  $h_i([1+a_i,1+x_i])=a'$ for $i\in A$ and $h_i([x_i,1+a_i))=a$, $h_i([1+a_i,1+x_i])=a'$, $h_i([x_i,1+b_i))=b$ and  $h_i([1+b_i,1+x_i])=b'$ for $i\in [7]\setminus A$. Clearly, $h(K)=h_1(K_1)\ldots h_7(K_7)\in S^7$ and $V=\{h(K): K\in \ka F_x\}$ is a partition code. Moreover,
$h^{-1}=h^{-1}_1\times \cdots \times h^{-1}_7$ preserves dichotomies, and $h^{-1}(V)=\ka F_x$  

If we take $a=0$, $a'=2$, $b=1$ and $b'=3$, then, as it was shown in \cite{De}, the maximum clique in a Keller graph has 124 vertices. This means that any polybox code $W\subset S^7$, $S=\{a,a',b,b'\}$, without twin pairs has at most 124 words. Since $|V|=128$, there is a twin pair in $V$, and consequently there is a twin pair in  $[0,1)^7+T$.   
\hfill{$\square$}

\medskip
Keller's conjecture is true in all dimensions $d\leq 6$, which was proved by Perron in 1940. Using the result of David Applegate on the size of the maximum clique in the Keller graph in dimension six announced in Section 1 and Theorem \ref{tll} we can give a new proof of Keller's conjecture for $d\leq 6$. 

\begin{tw}
\label{perr}
If $d\leq 6$, then every cube tiling $[0,1)^d+T$ of $\er^d$ contains a twin pair.
\end{tw}  
\proof Since every cube tiling of $\er^d$ with $r^+(T)=1$ contains a twin pair, it follows from Theorem \ref{tll} that Keller's conjecture is true for $d<6$. Let $d=6$. The maximum clique in the Keller graph in dimension six contains 60 vertices which was computed by David Applegate (see Section 1 in \cite{De}). Thus, in the same way as in the proof of  Corollary \ref{ke2} (compare also the explanation made in Introduction) we show that Keller's conjecture is true for every tiling $[0,1)^6+T$ with $r^-(T)\leq 2$. %Obviously, this means that the size of the maximal clique in a $d$-dimensional Keller graph  is less than $2^d$ for all dimensions $d< 6$.  Thus, in the same way as in we show that there is a twin pair in every cube tiling $[0,1)^d+T$ of $\er^d$ such that $r^-(T)\leq 2$ for $d\leq 6$. 
By Theorem \ref{tll}, the conjecture is true for all tilings $[0,1)^6+T$ of $\er^6$ for which $r^+(T)\geq 3$. Since $r^-(T)\leq 2$ or $r^+(T)\geq 3$ for every cube tiling $[0,1)^d+T$, the proof is completed.
\hfill{$\square$}  

\medskip
Since cliques in a $d$-dimensional Keller graph are polybox codes without twin pair which are written down in the alphabet $S=\{0,1,2,3\}$ with the complementation given by $0'=2$ and $1'=3$ we can define equivalent cliques in the Keller graph: Two cliques in a $d$-dimensional Keller graph with the vertex sets $V$ an $W$  are {\it equivalent} if  $\sum_{v\in V} g(v,w)= 2^d$  for every $w\in W$ and $\sum_{w\in W} g(w,v)= 2^d$ for every $v\in V$ (compare the comments just after (\ref{2d})). Thus,  Corollary \ref{w11} for cliques in the Keller graph reads as follows:   

\begin{wn}
\label{clik11} 
Two equivalent cliques in a $d$-dimensional Keller graph which have at most 11 vertices are equal.
\end{wn}

\medskip
We extend the notion of a $d$-dimensional Keller graph.  If $S$ is an alphabet with a complementation, then {\it a d-dimensional Keller graph on the set $S^d$} is the graph in which two vertices $u,v\in S^d$ are adjacent if they are dichotomous but do not form a twin pair. 

The only difference between a $d$-dimensional Keller graph and a  $d$-dimensional Keller graph on $S^d$ is that in the later 
 the set of vertices is $S^d$, where $S$ an arbitrary alphabet with a complementation, while in a $d$-dimensional  Keller graph we have $S=\{0,2,1,3\}$, where  $0'=2$ and $1'=3$.

From Corollary \ref{cylin} we obtain the following 
\begin{wn}
Every clique in a $d$-dimensional Keller graph on $S^d$ which contains at least $k> 2^{d-3}/3 $ vertices $u^1,\ldots ,u^k$ such that $u^n_i\not\in \{u^m_i,(u^m_i)'\}$ for some $i\in [d]$  and every $n,m\in \{1,...,k\},n\neq m$, has less than $2^d$ elements. 
In particular, any clique in the 7-dimensional Keller graph on $S^7$ which contains at least six vertices $u^1,\ldots ,u^6$ such that $u^n_i\not\in \{u^m_i,(u^m_i)'\}$ for some $i\in [7]$  and every $n,m\in \{1,...,6\},n\neq m$, has less than 128 elements. 
\end{wn}
\proof Assume on the contrary that there is a clique $W$ containing vertices $u^1,\ldots ,u^k$ and $|W|=2^d$. Thus, $W$ is a partition code without twin pairs. Let  $W=W^{i, l_1}\cup W^{i, l'_1}\cup \cdots \cup W^{i, l_r}\cup W^{i, l'_r}$. Since  $u^n_i\not\in \{u^m_i,(u^m_i)'\}$ for every $n,m\in \{1,...,k\},n\neq m$ and $u^n\in  W^{i, u^n_i}\cup W^{i, (u^n_i)'}$,  $u^m\in  W^{i, u^m_i}\cup W^{i, (u^m_i)'}$, it follows that $r>2^{d-3}/3 $. By Corollary \ref{cylin}, there is a twin pair in $W$, a contradiction.  
\hfill{$\square$} 

\medskip
\begin{uw}
{\rm 
Let us mention the three missing cases $r^+(T)\in \{3,4,5\}$. It seems that a proof of the case $r^-(T)\leq 2$ computed by Debroni et al. by hand computation is as yet out of reach. Similarly, the case $r^+(T)\geq 6$ are out of reach for computers (see e.g. Section 3 in  \cite{De}). Thus, these three missing cases might be an "intermediate" from the perspective of these two methods of proof. The easiest is the case $r^+(T)=5$, as now we can use Theorem \ref{12}. This opinion is based on some successful experiments which has been made in cooperation with Magdalena \L ysakowska. We believe that the case $r^+(T)=4$ can be attack along the same lines like $r^+(T)=5$ (using Theorem \ref{12} and the methods presented in the paper), but a computer support for this case will be much wider than that for $r^+(T)=5$.       

%we can only speculate: the case $|L(T,x,i)|=3$  will rather depend on a large-scale computations.

%, while the case $|L(T,x,i)|=4$ will evenly depend on the hand and computer computations.  
}
\end{uw}

%Let us mentioned that these three missing cases might be an "intermediate" from the perspective of a possible methods of proofs. A hand verification of the case computed by Debroni et al. seems to be out of reach. Similarly, the cases $|L(T,x,i)|\geq 6$ seems to be beyond the power of modern computers (see e.g. \cite{De}). 

%In \cite{De} Debroni  et al. showed that the maximal clique in a $7$-dimensional Keller graph, where $S=\{0,1,2,3\}$, contains 124 vertices.    

\end{document}